\documentclass[reqno,twoside,11pt,english]{amsart}
\usepackage{amsmath,amsfonts,amssymb,amsthm,epsfig}
\newcommand{\RNum}[1]{\uppercase\expandafter{\romannumeral #1\relax}}
\usepackage{comment}

\usepackage{color}
\usepackage[final,colorlinks,linkcolor=blue,anchorcolor=red,citecolor=blue]{hyperref}
\usepackage{graphicx}
\usepackage{titletoc,epsf}
\usepackage{amsmath,amsfonts,latexsym,amsthm,amsxtra,amssymb,bbm}
\allowdisplaybreaks
\usepackage{graphicx,float}
\usepackage{tikz}
\usetikzlibrary{intersections,patterns}
\usetikzlibrary{calc,spy}
\usepackage{ifthen}
\usepackage{float}

\allowdisplaybreaks

\newtheorem{thm}{Theorem}[section]
\newtheorem{lem}[thm]{Lemma}
\newtheorem{op}[thm]{Open Problem}
\newtheorem{cor}[thm]{Corollary}
\newtheorem{prop}[thm]{Proposition}
\newtheorem{assertion}{Assertion}[section]

\newtheorem{step}{Step}[section]

\newtheorem{cl}{Claim}[section]

\newtheorem{ca}{Case}
\newtheorem{sca}[section]{Subcase}
\newtheorem{scl}[section]{Subclaim}
\newtheorem{conj}[equation]{Conjecture}

\theoremstyle{definition}
\newtheorem{defn}[thm]{Definition}

\newtheorem{ques}[equation]{Question}
\newtheorem{rem}[thm]{Remark}
\newtheorem{exam}[thm]{Example}

\newcounter {own}
\def\theown {\thesection       .\arabic{own}}
\numberwithin{equation}{section}

\newenvironment{pf}[1][]{%
	\vskip 3mm
	\noindent
	\ifthenelse{\equal{#1}{}}%
	{{\slshape Proof. }}%
	{{\slshape #1.} }%
}%
{\qed\bigskip}

\newtheorem{Thm}{Theorem}

\newtheorem{Ques}[Thm]{Question}


\newcommand{\IR}{{\mathbb R}}

\newcommand{\IB}{{\mathbb B}}

\newcommand{\diam}{{\operatorname{diam}}}

\newcommand{\Card}{{\operatorname{Card}}}

\newcommand{\bas}{\begin{assertion}}
	\newcommand{\eas}{\end{assertion}}
\newcommand{\ben}{\begin{enumerate}}
	\newcommand{\een}{\end{enumerate}}
\newcommand{\bst}{\begin{step}}
	\newcommand{\est}{\end{step}}

\def\be{\begin{equation}}
	\def\ee{\end{equation}}

\newcommand{\bee}{\begin{enumerate}}
	\newcommand{\eee}{\end{enumerate}}

\newcommand{\blem}{\begin{lem}}
	\newcommand{\elem}{\end{lem}}
\newcommand{\bthm}{\begin{thm}}
	\newcommand{\ethm}{\end{thm}}
\newcommand{\bcor}{\begin{cor}}
	\newcommand{\ecor}{\end{cor}}
\newcommand{\beg}{\begin{exam}}
	\newcommand{\eeg}{\end{exam}}
\newcommand{\begs}{\begin{examples}}
	\newcommand{\eegs}{\end{examples}}
\newcommand{\bdefe}{\begin{defn}}
	\newcommand{\edefe}{\end{defn}}
\newcommand{\bprob}{\begin{prob}}
	\newcommand{\eprob}{\end{prob}}
\newcommand{\bques}{\begin{ques}}
	\newcommand{\eques}{\end{ques}}
\newcommand{\bei}{\begin{itemize}}
	\newcommand{\eei}{\end{itemize}}
\newcommand{\bcon}{\begin{conj}}
	\newcommand{\econ}{\end{conj}}
\newcommand{\bop}{\begin{op}}
	\newcommand{\eop}{\end{op}}

\newcommand{\bstep}{\begin{step}}
	\newcommand{\estep}{\end{step}}

\newcommand{\bca}{\begin{ca}}
	\newcommand{\eca}{\end{ca}}
\newcommand{\bsca}{\begin{sca}}
	\newcommand{\esca}{\end{sca}}

\newcommand{\bcl}{\begin{cl}}
	\newcommand{\ecl}{\end{cl}}

\newcommand{\bscl}{\begin{scl}}
	\newcommand{\escl}{\end{scl}}

\newcommand{\bcons}{\begin{conjs}}
	\newcommand{\econs}{\end{conjs}}
\newcommand{\bprop}{\begin{prop}}
	\newcommand{\eprop}{\end{prop}}
\newcommand{\br}{\begin{rem}}
	\newcommand{\er}{\end{rem}}
\newcommand{\brs}{\begin{rems}}
	\newcommand{\ers}{\end{rems}}
\newcommand{\bo}{\begin{obser}}
	\newcommand{\eo}{\end{obser}}
\newcommand{\bos}{\begin{obsers}}
	\newcommand{\eos}{\end{obsers}}
\newcommand{\bpf}{\begin{pf}}
	\newcommand{\epf}{\end{pf}}
\newcommand{\ba}{\begin{array}}
	\newcommand{\ea}{\end{array}}
\newcommand{\beq}{\begin{eqnarray}}
	\newcommand{\beqq}{\begin{eqnarray*}}
		\newcommand{\eeq}{\end{eqnarray}}
	\newcommand{\eeqq}{\end{eqnarray*}}

\newcommand{\ds}{\displaystyle}

\newcounter{minutes}
\divide\time by 60
\newcounter{hours}
\multiply\time by 60 \addtocounter{minutes}{-\time}

\begin{document}

	\title[Gromov hyperbolicity III: An improved geometric characterization]{Gromov hyperbolicity III: An improved geometric characterization and its applications}
	
	\author[C.-Y. Guo, M. Huang and X. Wang]{Chang-Yu Guo, Manzi Huang$^*$ and Xiantao Wang}

	\address[C.-Y. Guo]{Research Center for Mathematics and Interdisciplinary Sciences, Shandong University, 266237, Qingdao, P. R. China, and Department of Physics and Mathematics, University of Eastern Finland, 80101, Joensuu, Finland}
	\email{changyu.guo@sdu.edu.cn}
	
	\address[M. Huang]{MOE-LCSM, School of Mathematics and Statistics, Hunan Normal University, Changsha, Hunan 410081, P. R. China} \email{mzhuang@hunnu.edu.cn}
	
	\address[X. Wang]{MOE-LCSM, School of Mathematics and Statistics, Hunan Normal University, Changsha, Hunan 410081, P. R. China} \email{xtwang@hunnu.edu.cn}
	
	\date{}
	\dedicatory{Dedicated to Professor~Pekka Koskela on the occasion of his 65th birthday}
	\thanks{$^*$Corresponding author: Manzi Huang}
	\subjclass[2020]{Primary: 51F30, 30C65; Secondary: 30F45, 51F99}
	\keywords{Gromov hyperbolicity, Gehring-Hayman inequality, Ball separation condition,
		Quasihyperbolic geodesic, Inner uniform domain.}

\begin{abstract}
	
	In the seminal work of Balogh-Buckley [Invent. Math. 2003], the authors asked the following fundamental open problem: For a proper subdomain in Euclidean space, does the ball separation condition alone imply the Gehring-Hayman inequality? 
	
	In this paper, via a completely new measure-independent approach, we establish the following  geometric characterization of Gromov hyperbolicity in a fairly general setting: The Gromov hyperbolicity of a proper subdomain in a doubling metric space is quantitatively equivalent to the geometric ball separation condition, with explicit dependence on the coefficients. 
	In the special case of Euclidean spaces, it affirmatively solves the above Balogh-Buckely problem. Our result also significantly improves the main result of Koskela-Lammi-Manojlovi\'c [Ann. Sci. \'Ec. Norm. Sup\'er. 2014]. {As applications, we obtain the quasiconformal invariance of ball separation condition, a geometric characterization of inner uniformity in terms of ball separation condition, and the Gromov hyperbolicity of quasihyperbolic John length spaces. }
\end{abstract}

\date{\today}


\maketitle

\tableofcontents

\section{Introduction}\label{sec-1}

\subsection{Background}

Recall that the classical uniformization theorem states that the class of simply connected proper subdomains in $\IR^2$ can arise as conformal images of the unit disk $\mathbb{D}\subset \IR^2$, that is, Riemann mapping theorem. Searching for a suitable extension of such a beautiful theory in a higher dimensional Euclidean space $\IR^n$ or even an abstract metric space $X$, we may formally formulate the uniformization problem in the form of
\begin{equation}\label{eq:uniformization}
	\{\text{a class of good domains in } \IR^n \text{ or } X\}=\mathcal{F}\left(\{\text{a class of model domains in } \IR^n \text{ or } X \}\right),\tag{UP}
\end{equation}
where $\mathcal{F}$ consists of  homeomorphisms with good geometric properties (such as conformal maps).

In a seminal work \cite{BHK}, Bonk, Heinonen and Koskela have successfully developed a rich uniformization theory for \eqref{eq:uniformization}. To record their theory, we need to recall a couple of basic definitions. The first notion is the class of (inner) uniform domains in a general metric space.

\begin{defn}\label{def:uniform domain}
	A domain $\Omega$ in a metric space $X=(X,d)$ is called  {\it $c$-uniform}, $c\geq 1$, if each pair of points $z_{1},z_{2}$ in $\Omega$ can be joined by a rectifiable curve $\gamma$ in $\Omega$ satisfying
	\begin{enumerate}
		\item\label{con1} $\ds\min_{j=1,2}\{\ell_d (\gamma [z_j, z])\}\leq c\, d_\Omega(z)$ for all $z\in \gamma$, and
		\item\label{con2} $\ell_d(\gamma)\leq c\,d(z_{1},z_{2})$,
	\end{enumerate}
	where $\ell_d(\gamma)$ denotes the arc-length of $\gamma$ with respect to the metric $d$,
	$\gamma[z_{j},z]$ the subcurve of $\gamma$ between $z_{j}$ and $z$, and $d_\Omega(z):=d(z,\partial \Omega)$, the distance from $z$ to the boundary $\partial \Omega$ of $\Omega$. In a $c$-uniform domain $\Omega$, any curve $\gamma\subset \Omega$, which satisfies the conditions (1) and (2) above, is called a {\it $c$-uniform curve}. 
\end{defn}

If the condition \eqref{con2} in Definition \ref{def:uniform domain} is replaced by the weaker inequality
\begin{equation}\label{eq:def for inner quasiconvexity}
	\ell_d(\gamma)\leq c\,\sigma_\Omega( z_{1}, z_{2}),
\end{equation}
where $\sigma_\Omega$ is the {\it inner metric} in $\Omega$ defined by
$$\sigma_\Omega(z_1,z_2)=\inf \{\ell(\alpha):\; \alpha\subset \Omega\;
\mbox{is a rectifiable curve joining}\; z_1\; \mbox{and}\; z_2 \},$$
then $\Omega$ is said to be $c$-{\it inner uniform}, and the corresponding curve $\gamma$ is called a $c$-\emph{inner uniform curve}. When the context is clear, we often drop the subscript $\Omega$ from $\sigma_{\Omega}$ and simply write $\sigma$.

If $\Omega$ only satisfies the condition \eqref{con1} in Definition \ref{def:uniform domain},  then it is said to be a {\it $c$-John domain}, and the corresponding curve $\gamma$ is called a \emph{$c$-John curve} or a {\it double $c$-cone curve}.

The class of John domains was initially introduced by John in his study of elasticity \cite{Jo}, and the name was coined by Martio and Sarvas in \cite{MS}, where they also introduced the class of uniform domains. These classes of domains are central in the modern geometric function theory in $\mathbb{R}^n$ or more general metric spaces, and also, they have wide connections with many other mathematical subjects related to analysis and geometry; see, for instance, \cite{GM,GO,GePa,GGGKN-2024,Guo-2015,Hajlasz-Koskela-2000,H,Hei-Book-2001,HeKo,HLPW,Jo81,Martin-1985}.

The second notion is the class of Gromov hyperbolic spaces introduced by Gromov in his celebrated work \cite{Gr1}. Recall that
\begin{defn}
	A geodesic metric space $X=(X,d)$ is called \emph{$\delta$-Gromov hyperbolic}, $\delta>0$, if each side of a geodesic triangle in $X$ lies in the $\delta$-neighborhood of the other two sides.
\end{defn}
The Gromov hyperbolicity is a large-scale property, which generalizes the metric properties of classical hyperbolic geometry and trees, and it turns out to be very useful in geometric group theory and metric geometry \cite{Bonk-Kleiner-2005,BS,Bridson-H-1999,Gr2}. For more about the Gromov hyperbolicity and its connection to geometric mapping theory, see, for instance, \cite{BB2,BB4,BHK,Herron-S-Xie-2008,HRWZ-2025,Ntalampekos-K-2025-Duke,Sidler-Wenger-2021,Zhou-R-2022-IJM}.

To define the class of Gromov hyperbolic domains, we first recall the definition of quasihyperbolic metric, which was initially introduced by Gehring and Palka \cite{GePa} for domains in $\IR^n$, and then, has been extensively studied in \cite{GO}.
The {\it quasihyperbolic length} of a rectifiable curve
$\gamma$
in a proper subdomain $\Omega\subsetneq (X,d)$ is defined as
$$
\ell_{k}(\gamma):=\int_{\gamma}\frac{|dz|}{d_\Omega(z)}.
$$

For any $z_1$, $z_2$ in $\Omega$, the {\it quasihyperbolic distance}
$k_\Omega(z_1,z_2)$ between $z_1$ and $z_2$ is set to be
$$k_\Omega(z_1,z_2)=\inf_{\gamma}\{\ell_k(\gamma)\},
$$
where the infimum is taken over all rectifiable curves $\gamma$
joining $z_1$ and $z_2$ in $\Omega$. A rectifiable curve $\gamma$ from $z_1$ to $z_2$ is called a {\it quasihyperbolic geodesic} if
$\ell_k(\gamma)=k_\Omega(z_1,z_2)$. Clearly, each subcurve of a quasihyperbolic
geodesic is a quasihyperbolic geodesic.

\begin{defn}\label{def:Gromov hyperbolic domains}
	Let $(X,d)$ be a metric space and $k$ the quasihyperbolic metric  {induced by $d$}. A domain $\Omega\subsetneq X$ is called {\it $\delta$-Gromov hyperbolic} if the metric space $(\Omega,k)$ is $\delta$-Gromov hyperbolic.
\end{defn}

Now, we are able to present the uniformization theory of Bonk-Heinonen-Koskela \cite{BHK}: On the left hand side of \eqref{eq:uniformization}, one takes the class of Gromov hyperbolic domains, while on the right hand side, one uses the well-known class of uniform domains. Then Bonk, Heinonen and Koskela have established in their main result \cite[Theorem 1.1]{BHK} a rather general uniformization theory for Gromov hyperbolic spaces:
\[
\{\text{some class of Gromov hyperbolic spaces}\}=\mathcal{F}\left(\{\text{bounded locally compact uniform spaces}\}\right),
\] 	
where $\mathcal{F}$ is a class of good homeomorphisms (i.e., quasiisometries). 

In their second main result \cite[Theorem 1.11]{BHK}, Bonk, Heinonen and Koskela proved that all (inner) uniform domains in $\IR^n$ are Gromov hyperbolic, which provides a large class of (nontrivial) examples of Gromov hyperbolic domains in $\IR^n$. It is well known that each quasiconformal image of a Gromov hyperbolic domain is again Gromov hyperbolic in $\IR^n$ {\cite[Chapter 7, $\S$4]{GhHa}. } Then Riemann mapping theorem implies that all simply connected domains in $\IR^2$ are Gromov hyperbolic. Furthermore, it was shown that a Gromov hyperbolic domain $\Omega$ in $\overline{\IR}^n$ equipped with the spherical metric $d_s$ necessarily satisfies both \emph{the Gehring-Hayman inequality} and \emph{the ball separation condition}.
For later references, we recall these two fundamental concepts here. 
\begin{defn}[Gehring-Hayman inequality]\label{def:Gehring-Hayman}
	A domain $\Omega$ in a metric space $(X,d)$ is said to satisfy the \emph{$C$-Gehring-Hayman inequality}, $C\geq 1$, if for all $x,y\in \Omega$ and each quasihyperbolic geodesic $\gamma_{xy}\subset \Omega$ with end points $x$ and $y$, it holds
	\begin{equation}\label{eq:def for Gehring-Hayman}
		\ell_d(\gamma_{xy})\leq C\sigma(x,y).
	\end{equation}
	
	To emphasize the metric $d$, we also say that $(\Omega,d)$ satisfies the $C$-Gehring-Hayman inequality. 
\end{defn}

\begin{defn}[Ball separation condition]\label{def:ball separation}
	A domain $\Omega$ in a metric space $(X,d)$ is said to satisfy the \emph{$C$-ball separation condition}, $C>0$, if for each quasihyperbolic geodesic $\gamma_{xy}\subset \Omega$, every $z\in \gamma_{xy}$ and every rectifiable curve $\gamma\subset \Omega$ joining $x$ and $y$, it holds
	\begin{equation}\label{eq:def for ball separation}
		\IB_{\sigma}\left(z,Cd_{\Omega}(z)\right)\cap \gamma\neq \emptyset.
	\end{equation}
	Here, $\IB_{\sigma}\left(z,r\right)=\{y\in X:\; \sigma(z,y)<r\}$,  the open ball centered at $z$ with radius $r$ in the inner metric $\sigma$. To emphasize the inner metric $\sigma$, we also say that $(\Omega,\sigma)$ satisfies the $C$-ball separation condition.
\end{defn}

The Gehring-Hayman inequality first appeared in the pioneer work \cite{Geh}, where the authors proved a similar inequality for hyperbolic geodesics in simply connected planar domains and the ball separation condition was first introduced by Buckley and Koskela \cite{BK1} in their study of Sobolev-Poincar\'e inequalites in Euclidean domains. 
These two properties are potentially much easier to verify than the more complicated Gromov hyperbolicity, and thus, it is natural to ask \emph{whether these two properties actually fully characterize the Gromov hyperbolicity}. In fact, Bonk, Heinonen and Koskela made the following conjecture.
\vspace{0.1cm}

\textbf{Conjecture} ({\cite[Page 75]{BHK}}): If $\Omega\subset (\overline{\IR}^n,d_s)$ satisfies both the Gehring-Hayman inequality and the ball separation condition, then $(\Omega,k)$ is Gromov hyperbolic.
\medskip

In another seminal work, built on the fundamental works in \cite{BHK} and \cite{HeKo}, Balogh and Buckley \cite[Theorem 0.1]{BB4} have confirmed affirmatively the above conjecture of Bonk-Heinonen-Koskela. As a direct application of their main result \cite[Theorem 0.1]{BB4}, they obtained the following interesting geometric characterization of Gromov hyperbolicity for \emph{proper Euclidean domains}.
\begin{Thm}[{Balogh-Buckley, \cite{BB4}}]\label{Thm-A}
	Let $\Omega\subset \IR^n$ be a proper subdomain. Then the following two statements are equivalent:
	\begin{enumerate}
		\item $(\Omega, k)$ is $\delta$-Gromov hyperbolic.
		\item $\Omega$ satisfies both the $C$-Gehring-Hayman inequality and the $C$-ball separation condition.
	\end{enumerate}
	Moreover, the constants $\delta$ and $C$ depend on each other and on the dimension $n$.
\end{Thm}

A fundamental problem after this work is the relationship between the Gehring-Hayman inequality and the ball separation condition. In \cite{BB4}, Balogh and Buckley constructed a planar domain that satisfies the Gehring-Hayman inequality, but does not have the ball separation condition, showing that the Gehring-Hayman inequality alone does not imply the ball separation condition. Then they asked for the reverse implication, which has become a fundamental open problem in the field.

\begin{Ques}[{\cite[Page 272]{BB4}}]\label{Ques:BB}
	Does the ball separation condition alone imply the Gehring-Hayman inequality for proper Euclidean domains?
\end{Ques}

Due to its importance in the study of Gromov hyperbolic quasihyperbolization of quasihyperbolic John length spaces, very recently, this question was proposed again by Zhou and Ponnusamy \cite[the third paragraph after Question 1.9]{ZP-2024-Pisa}.  
Notice that a positive answer to Question \ref{Ques:BB} would lead to the following concrete geometric characterization of Gromov hyperbolicity: \emph{A proper subdomain $\Omega\subset \IR^n$ is Gromov hyperbolic if and only if it satisfies the ball separation condition}. As the ball separation condition is significantly much easier to verify than the Gromov hyperbolicity, this would be very valuable for the theory of Gromov hyperbolic domains. For instance, it provides a lot of nontrivial interesting examples of Gromov hyperbolic domains.

Notice that all the three conditions in Theorem \ref{Thm-A}  are based only on purely metric concepts, and thus, it is natural to ask for an extension of such a beautiful characterization to abstract metric spaces. Such an extension was initially given in \cite{BB4}, which relies essentially on the existence of a suitable Poincar\'e inequality (or being Loewner) for the underlying spaces (see \cite[the last paragraph on page 265]{BB4}). As was pointed out by Koskela, Lammi and Manojlovi\'c \cite{KLM}, supporting an abstract Poincar\'e inequality (or being Loewner) is a somewhat restrictive assumption for the underlying spaces. In their main result \cite[Theorem 1.2]{KLM}, they have successfully extended Theorem \ref{Thm-A} to the setting of locally compact \emph{$Q$-regular} length spaces that are additionally \emph{annularly quasiconvex}. More precisely, they proved the following result.

\begin{Thm}[{\cite[Theorem 1.2]{KLM}}]\label{Thm-B}
	Let $Q>1$, $C_0\geq 1$, and let $(X,d,\mu)$ be a $(Q,C_0)$-regular metric measure space with $(X,d)$ a locally compact and annularly $C_1$-quasiconvex length space. Let $\Omega\subset X$
	be a bounded proper subdomain.
	Then the following statements are equivalent:
	\begin{enumerate}
		\item $(\Omega, k)$ is $\delta$-Gromov hyperbolic.
		\item $\Omega$ satisfies both the $C$-Gehring-Hayman inequality and the $C$-ball separation condition.
	\end{enumerate}
\end{Thm}

Recall that for $Q\geq 1$ and $C_0\geq 1$, we say that a metric measure space
$X=(X,d,\mu)$ is \emph{(Ahlfors) $(Q,C_0)$-regular} if for each $x\in X$ and $0<r\leq \diam_d(X)$,
$${C_0}^{-1}r^Q\leq {\mu(B_d(x,r))}\leq {C_0}r^Q.$$
In the following, sometimes, we shall ignore the constant $C_0$ and simply say that $X$ is \emph{(Ahlfors) $Q$-regular}.

Theorem \ref{Thm-B} improves the corresponding (metric space) result of Balogh-Buckley by weakening the requirement on the underlying spaces, from supporting an abstract Poincar\'e inequality (or being Loewner) to a weaker geometric assumption, that is, the annular quasiconvexity. But the assumption on the existence of an Ahlfors regular measure is somewhat unnatural as all the three conditions are purely metric. It is thus natural to ask
\begin{Ques}\label{Ques:metric version}
	Does the geometric characterization of Gromov hyperbolicty hold in general metric spaces without any reference measure?
\end{Ques}
Note that the constants involved  in Theorem \ref{Thm-B} depend not only on each other and the dimension $Q$, but also on \emph{$\diam(\Omega)$} and \emph{the constant associated with the annular quasiconvexity}. One would wonder whether the dependence on $\diam(\Omega)$ and $C_1$ are really necessary. These basic questions are the direct motivation of the present paper.

\subsection{Main results}

As was briefly pointed out in the previous subsection, all the known proofs for the geometric characterization of Gromov hyperbolicity  (e.g., Theorems \ref{Thm-A} and \ref{Thm-B}) depend heavily on the uniformization theory of Bonk-Heinonen-Koskela \cite{BHK}, in particular, on the Lebesgue measure and integration in $\IR^n$ or abstract metric measure spaces. The dependence on the dimension $n$ (or $Q$) comes from the Ahlfors $n$-regularity of the Lebesgue $n$-measure in $\IR^n$ (or the Ahlfors $Q$-regularity of the measure in metric measure spaces).

In this paper, we shall provide a new elementary measure-independent proof of Theorem \ref{thm:new main metric space} below ({and thus, is completely different from the proofs of Theorems \ref{Thm-A} and \ref{Thm-B}),} using only the metric doubling property of the underlying spaces. Recall that
\begin{defn}\label{def:Q doubling ms}
	A metric space $X=(X,d)$ is called \emph{$Q$-doubling}, if there is an integer $Q\geq 2$ such that for each ball $\IB(x,r)$, every $r/2$-separated {subset} of $\IB(x,r)$ has at most $Q$ points.
\end{defn}
Here and hereafter,
$\IB(x,r):=\{y\in X: d(x,y)<r\}$ represents the open ball centered at $x$ with radius $r$ in the metric $d$, and a {\it $\nu$-separated} set, $\nu>0$, means a set in $X$ such that every two distinct points in the set have distance
at least $\nu$.
A simple volume comparison implies that $\IR^n$ (equipped with the standard Euclidean metric) is $2^n$-doubling, according to Definition \ref{def:Q doubling ms}.

Our main result of this paper shows surprisingly that the Gromov hyperbolicity of a proper subdomain in a doubling metric length space is \emph{completely characterized by the {geometric} ball separation condition}. Note that no measure is involved in such a characterization.

\begin{thm}\label{thm:new main metric space}
	Let $X=(X, d)$ be a $Q$-doubling length space and $\Omega\subset X$
	a proper subdomain such that $(\Omega,k)$ is geodesic. Then $(\Omega,k)$ is $\delta$-Gromov hyperbolic if and only if $\Omega$ satisfies the $C$-ball separation condition. The statement is quantitative in the sense that the coefficients involved depend only on each other and $Q$, and can be expressed explicitly in terms of these data.
\end{thm}

Note that the assumption of $(\Omega,k)$ being geodesic is quite natural as it is part of the definition of a Gromov hyperbolic domain and also the definition of the ball separation condition. If $X$ is a locally compact length space, then the identity map $id:(\Omega,d)\to (\Omega,\sigma_{\Omega})$ is a homeomorphism, and so, by \cite[Proposition 2.8]{BHK}, $(\Omega,k)$ is geodesic. {This property remains valid for nice domains in certain infinite dimensional spaces as well; see for instance the work of Martio and V\"{a}is\"{a}l\"{a} \cite{Martio-Vaisala-2011}.  }

Besides getting rid of the Gehring-Hayman inequality for a geometric characterization of Gromov hyperbolicity, comparing with Theorem \ref{Thm-B}, there are two other major improvements in Theorem \ref{thm:new main metric space}:
\begin{itemize}
	\item The statement is purely metric and no measure is needed, which in particular provides an affirmative answer to Question \ref{Ques:metric version}.
	\item It removes the annular quasiconvexity assumption on $X$ and the boundedness assumption on $\Omega$ as required in Theorem \ref{Thm-B}.
	The constants $\delta$ and $C$ depend only on each other and on $Q$, not additionally on $\diam(\Omega)$ nor the annular quasiconvexity constant.
\end{itemize}

Theorem \ref{thm:new main metric space} is new even in the case of Euclidean spaces: \emph{For a proper subdomain $\Omega\subset \IR^n$,  $(\Omega, k)$ is Gromov hyperbolic if and only if $\Omega$ satisfies the ball separation condition, quantitatively}. This not only improves the main result of Balogh-Buckley, Theorem \ref{Thm-A}, but also solves affirmatively their fundamental open problem, Question \ref{Ques:BB}.

It is a well-known result, due to Buckley and Koskela {\cite[Lemma 3.3]{BK1}, } that quasiconformal images of uniform domains in $\IR^n$ satisfy the ball separation condition. 
	As a corollary of Theorem \ref{thm:new main metric space} and the quasiconformal invariance of Gromov hyperbolicity in Euclidean spaces, we obtain an important extension of this result: Quasiconformal images of domains with the ball separation condition in Euclidean spaces again satisfy the ball separation condition.
	\begin{cor}[Quasiconformal invariance of ball separation condition]\label{coro:quasiconformal invariance of ball sp}
		Let $\Omega$ and $\Omega'$ be two proper subdomains in $\IR^n$ with $n\geq 2$. If $\Omega$ is $K$-quasiconformally equivalent to $\Omega'$ and satisfies the $\theta$-ball separation condition, then $\Omega'$ satisfies the $\theta'$-ball separation condition, quantitatively.
	\end{cor}
	
	In the final section, Section \ref{sec:6}, we shall give two more geometric applications of our main results. The first one (see Theorem \ref{thm:ball separation+LLC2 implies GH metric space}) gives a concrete geometric characterization of inner uniformity in terms of  ball separation condition, and the second one (see Theorem \ref{thm:Zhou problem}) provides a sufficient condition for a quasihyperbolic John length space to be Gromov hyperbolic (with respect to the quasihyperbolic metric), which provides an affirmative answer to a question raised by Zhou and Ponnusamy \cite[Question 1.9]{ZP-2024-Pisa} in the setting of $Q$-doubling length spaces.

	\subsection{Outline of the proofs}
	In this subsection, we briefly outline the proofs of our results. {It is known that the proofs of Theorem \ref{Thm-A} and/or Theorem \ref{Thm-B} depend on the uniformization theory of Bonk-Heinonen-Koskela \cite{BHK}}, which relies crucially on the Lebesgue measure and integration. Thus these proofs cannot be extended to show our main result, Theorem \ref{thm:new main metric space}, as a measure is needed. Our starting point is to find a new purely metric (and thus, measure-independent) proof of Theorem \ref{thm:new main metric space}.
	
	 The new method developed in this paper seems to be quite robust and we believe that it will have potential applications in many other problems related to quasiconformal maps and Gromov hyperbolic spaces. We shall explore these aspects  in our future works. 
	
	\subsubsection*{On the proof of Theorem \ref{thm:new main metric space}, if part}
	
	The proof of ``if part" in Theorem \ref{thm:new main metric space} consists of two main steps. In the first step, we prove that the ball separation condition alone implies the Gehring-Hayman inequality. In the special case of Euclidean spaces, it affirmatively answers Question \ref{Ques:BB}.
	\begin{thm}\label{positive-answer}
		Fix $Q>1$. Let $X=(X, d)$ be a $Q$-doubling  length space and $\Omega\subset X$
		a proper subdomain such that $(\Omega,k)$ is geodesic. If $\Omega$ satisfies the $\theta$-ball separation condition, then it has the $\theta_1$-Gehring-Hayman inequality, where $\theta_1$ can be expressed explicitly in terms of $Q$ and $\theta$.
	\end{thm}

	At the level of idea, our proof for this step is partially inspired by our recent work \cite{GHW}, where we proved the Gehring-Hayman inequality for certain special domains in general Banach spaces. 
	However, there are essential difficulties to adapt the argument of \cite{GHW} to the current setting: In \cite{GHW}, $\Omega\subsetneq X$ is a proper subdomain that is quasihyperbolically homeomorphic to a uniform domain. This class of domains enjoy much better analytic/geometric properties than merely having the ball separation property. In fact, most of the analysis in \cite[Sections 3, 4 and 5]{GHW} relies heavily on good properties of uniform domains and quasihyperbolic homeomorphisms. In particular, it is not even clear how to extend the proofs to the ``better" class of Gromov hyperbolic domains, let alone the class of ``weaker" domains considered in this paper. Furthermore, in \cite{GHW}, $X$ is assumed to be a Banach space, which has nice linear structure. In our setting, $X$ is a general metric space. Some of the arguments in \cite[Sections 2 and 3]{GHW} relies on this linear structure and does not extend to nonlinear spaces. 
%
	Thus new ideas are really necessary for the proof of Theorem \ref{positive-answer}. There are two main new ingredients in our approach:
	\begin{enumerate}
		\item A version of diameter Gehring-Hayman inequality in Theorem \ref{2017-10-11-1}, which asserts that the diameter (with respect to the inner metric $\sigma$) of a quasihyperbolic geodesic $\gamma$ is quantitatively bounded by the inner distance between end points of $\gamma$. 
		
		\item New constructions of partitions on the quasihyperbolic geodesic $\gamma$ and the ``almost shortest curve" $\alpha$ in Proposition \ref{2017010-20-1}, which provides quantitative control on the quasihyperbolic distance between consecutive points on associated curves $\gamma$ and $\alpha$.
	\end{enumerate}
	In both (1) and (2) above, we shall do some novel constructions to select points with good controls on quasihyperbolic/inner distances, with the aid of $Q$-doubling assumption. We believe this kind of constructions will be of independent interest. 
	

	In the second step, we show that the ball separation condition, together with the Gehring-Hayman inequality, implies the Gromov hyperbolicity. Note that in this step, the $Q$-doubling assumption for $X$ is not needed.
	\begin{thm}\label{thm:sufficient for Gromov hyperbolic}
		Suppose that $X=(X,d)$ is a metric space, and let $\Omega\subset X$ be a proper subdomain such that $(\Omega,k)$ is geodesic. If $\Omega$ satisfies both the $C$-Gehring-Hayman inequality and the $C$-ball separation condition, then $(\Omega, k)$ is $\delta$-Gromov hyperbolic with $\delta=\max\{50 C^6(3+C)^2,100^{10}\}$.
	\end{thm}
	
	
	We shall prove Theorem \ref{thm:sufficient for Gromov hyperbolic}, again by using a contradiction argument. To be more precise, fix an arbitrary geodesic triangle $\Delta_{x_1x_2x_3}$ in $\Omega$ and a point $x_0\in \gamma_{x_1x_2}$, one side of the triangle $\Delta_{x_1x_2x_3}$. Suppose,  on the contrary, that
	$$k_{\Omega}(x_0,\gamma_{x_1 x_3}\cup \gamma_{x_2 x_3})> 50 C^6(3+C)^2.$$
	
	Let $y_0 \in \gamma_{x_1x_3} \cup \gamma_{x_2 x_3}$ be such that
	$$
	\sigma(x_0,y_0)\approx \inf\limits_{y\in\gamma_{x_1x_3} \cup \gamma_{x_2 x_3}}\sigma(x_0,y).
	$$
	Then the key step is to find points with good control on geodesic triangles: There exist two points  $z_{1} \in \gamma_{x_{1} x_{2}}[x_{1}, x_0]$  and  $z_{2} \in \gamma_{x_{1}x_{2}}[x_0, x_{2}]$  such that
	\[
	\sigma(y_0, z_{1}) \leq \frac{3}{5(3+C)^{2}} \sigma(x_0, y_0)\ \text{ and }\  \sigma(y_0, z_{2}) \leq \frac{3}{5 (3+C)^2} \sigma(x_0, y_0).
	\]
	
	Once the above estimates are proved, namely Claim \ref{claim:find z1 and z2} in Section \ref{sec-3}, it follows that
	\begin{align*}
		\ell(\gamma_{x_1 x_2}[z_1, z_2])& \geq \sigma(z_1, x_0)+\sigma(z_2, x_0) \\
		& \geq \sigma(y_0, x_0)-\sigma(y_0, z_1)+\sigma(y_0, x_0)-\sigma(y_0, z_2) \\
		& \geq\left(2-\frac{6}{5(3+C)^{2}}\right) \sigma(y_0, x_0).
	\end{align*}
	
	On the other hand, the $C$-Gehring-Hayman inequality gives
	\begin{align*}
		\ell(\gamma_{x_1 x_2}[z_1, z_2])  \leq C \sigma(z_1, z_2) \leq C(\sigma(z_1, y_0)+\sigma(y_0, z_2))  \leq \frac{6C}{5(3+C)^2} \sigma(y_0, x_0),
	\end{align*}
	which clearly contradicts with the previous estimate.
	
	A novel point for the proof of Claim \ref{claim:find z1 and z2} is that it is simply based on repeated applications of ball separation condition and Gehring-Hayman inequality. Comparing with the proofs of Theorems 2.4 and 6.1 in \cite{BB4}, our proof here is more elementary and direct.
	
	\subsubsection*{On the proof of Theorem \ref{thm:new main metric space}, only if part}
	
	The proof of ``only if part" in Theorem \ref{thm:new main metric space} is given in Section \ref{sec-4-3}; see Theorem \ref{thm-24-5.1} below. As before, we shall apply nontrivial contradiction arguments.

	Suppose that the ball separation condition fails. Then for a quasihyperbolic geodesic $\gamma$ and another rectifiable curve $\alpha$ with the same end points, we may find a point $x_{0,0}\in \gamma$ so that for sufficiently large $\tau$, it holds
	\begin{equation}\label{eq:explain}
		\sigma(x_{0,0},\alpha)>\tau d_{\Omega}(x_{0,0}). \tag{*}
	\end{equation}
	Then a great effort is made to find a sequence of points on $\alpha$ with prescribed control on the quasihyperbolic distances between successive points. Then a contradiction occurs by comparing the quasihyperbolic distance between the final pair of points if $\tau$ in \eqref{eq:explain} is large enough.
	
	More precisely, as $(\Omega,k)$ is $\delta$-Gromov hyperbolic, we may find a sequence of points on $\alpha$ as in the following claim (see Claim \ref{H25-084} below):
	\medskip
	
	\textbf{Claim A:}  Let $N=[\frac{\ell_k(\alpha_{xy})}{3C+1}]+1$. Then for each positive integer $\varsigma\in\{1,\cdots, N\}$, there are $y_{0,\varsigma}^1\in\alpha_{x y}$ and  $y_{0,\varsigma}^2\in\alpha_{x y}[y_{0,\varsigma}^1,y]$ which satisfy the following:
	\ben
	\item
	For each $\varsigma\in\{2,\cdots, N\}$, $y_{0,\varsigma}^1\in\alpha_{x y}[y_{0,\varsigma-1}^1,y]$ and $k_{\Omega}(y_{0,\varsigma-1}^1,y_{0,\varsigma}^1)\geq1+3C$.
	
	\item
	For each $\varsigma\in\{1,\cdots, N\}$ and every $\gamma_{y_{0,\varsigma}^1y_{0,\varsigma}^2}\in \Lambda_{y_{0,\varsigma}^1y_{0,\varsigma}^2}(\Omega)$, there exists $x_{0,\varsigma}\in\gamma_{y_{0,\varsigma}^1y_{0,\varsigma}^2}$ such that
	$$\sigma(x_{0,\varsigma},\alpha[y_{0,\varsigma}^1, y_{0,\varsigma}^2])>\tau d_{\Omega}(x_{0,\varsigma}).$$
	
	\item 
	For each $\varsigma\in\{1,\cdots, N-1\}$, $k_{\Omega}(y_{0,\varsigma}^1,y_{0,\varsigma}^2)\geq 1+3C$ and $k_{\Omega}(y_{0,N}^1,y_{0,N}^2)<1+3C$.
	\een
	\medskip
	
	Here and hereafter, for two points $u$ and $v$ in $\Omega$, we always use $\Lambda_{uv}(\Omega)$ to stand for the set of all quasihyperbolic geodesics in $\Omega$ with end points $u$ and $v$ and $\gamma_{uv}$ an element in $\Lambda_{uv}(\Omega)$.
	
	Suppose that \textbf{Claim A} holds. Then
	$$
	\begin{aligned}
		k_{\Omega}(y_{0,N}^1,y_{0,N}^2)&\geq k_{\Omega}(y_{0,N}^1,x_{0,N}) \stackrel{\eqref{(2.1)}}{\geq} \log\left(1+\frac{\ell(\gamma_{y_{0,N}^1y_{0,N}^2}[y_{0,N}^1,x_{0,N}])}{\min\{d_{\Omega}(y_{0,N}^1), d_{\Omega}(x_{0,N})\}}\right)\\
		&\geq \log\left(1+\frac{\sigma(y_{0,N}^1,x_{0,N})}{d_{\Omega}(x_{0,N})}\right)
		\geq \log \frac{\sigma(x_{0,N},\alpha[y_{0,N}^1,y_{0,N}^2])}{d_{\Omega}(x_{0,N})}>\log\tau,
	\end{aligned}
	$$
	which clearly contradicts with \textbf{Claim A}(3) if $\tau>e^{1+3\delta}$.

	The key ingredient for the proof of \textbf{Claim A} is the construction of point sequences given by Lemma \ref{Lem4.12-1} (or its iterated version Lemma \ref{cl25-047}).  A novel point in our proof of this key lemma is to introduce three new classes $P_{\alpha}^\gamma$, $O_{\alpha}^\gamma$ and $Q_{\alpha}^\gamma$ of curves based on $\gamma$ and $\alpha$ (see Definitions \ref{def:class P}, \ref{def:Class O} and \ref{def:Class AC} below), on which there are certain special points with controlled quasihyperbolic distances. Some basic properties for these curve families in Gromov hyperbolic domains are then developed in Subsection \ref{sec-4-1} and the appendix. \textbf{Claim A} will follows from these basic estimates, together with a smart contradiction argument using the $Q$-doubling property of $X$.

	\medskip
	\textbf{Structure.} The structure of this paper is as follows. In Section \ref{sec-2},
	we recall some basic facts about the quasihyperbolic metric, Gromov hyperbolic spaces and doubling metric spaces. In Section \ref{sec-4-3}, we prove the ``only if part" of Theorem \ref{thm:new main metric space}: The Gromov hyperbolicity implies the ball separation condition. The proof of ``if part" of Theorem \ref{thm:new main metric space} is divided into two sections. In Section \ref{sec-4-2}, we show that the ball separation condition implies the Gehring-Hayman inequality, namely Theorem \ref{positive-answer}, while in Section \ref{sec-3}, we demonstrate Theorem \ref{thm:sufficient for Gromov hyperbolic}: The ball separation condition, together with the Gehring-Hayman inequality, implies the Gromov hyperbolicity. Section \ref{sec:6} contains two more geometric applications of our results.
	\smallskip
	
	\textbf{Notations.} Throughout this paper, 
	for each proper subdomain $\Omega\subset X$, we use $\Lambda_{xy}(\Omega)$ to represent the set of all quasihyperbolic geodesics in $\Omega$ with end points $x$ and $y$, and use $\gamma_{xy}$ to denote some quasihyperbolic geodesic in $\Lambda_{xy}$.
	Meanwhile, we use  $\Gamma_{xy}(\Omega)$ to represent the set of all curves in $\Omega$ with end points $x$ and $y$.

	\section{Preliminaries}\label{sec-2}

	\subsection{Estimates on quasihyperbolic distance}\label{subsec 2-0}

	For any $z_1$, $z_2$ in $\Omega$, let $\gamma\in \Gamma_{z_1z_2}(\Omega)$, where $\Gamma_{z_1z_2}(\Omega)$ denotes the set of all rectifiable curves in $\Omega$ with end points $z_1$ and $z_2$.
	Then we have the following elementary estimates (see, for instance, \cite[Section 2]{Vai3}):
	\beq\label{(2.1)}
	\ell_{k}(\gamma)\geq
	\log\Big(1+\frac{\ell(\gamma)}{\min\{d_\Omega(z_1), d_\Omega(z_2)\}}\Big)
	\eeq
	and
	\beq\label{(2.2)}
	\begin{aligned}
		k_{\Omega}(z_1, z_2) &\geq  \log\Big(1+\frac{\sigma_\Omega(z_1,z_2)}{\min\{d_\Omega(z_1), d_\Omega(z_2)\}}\Big)
		\\
		&\geq
		\log\Big(1+\frac{|z_1-z_2|}{\min\{d_\Omega(z_1), d_\Omega(z_2)\}}\Big)
		\geq
		\Big|\log \frac{d_\Omega(z_2)}{d_\Omega(z_1)}\Big|.
	\end{aligned}
	\eeq
	
	The following two estimates are useful in our later proofs.
	\begin{lem}\label{lem-3-1}
		Suppose that $\Omega\subset X$ is a domain, $u$, $v\in \Omega$, and $\alpha\in \Gamma_{uv}(\Omega)$. Fix $c\geq 1$. If for each $w\in\alpha$,  $\ell(\alpha[u,w])\leq c d_\Omega(w)$, then
		$$k(u,w)\leq 2c \log\Big(1+\frac{2\ell(\alpha[u,w])}{d_\Omega(u)}\Big)\leq 4c \log\Big(1+\frac{\ell(\alpha[u,w])}{d_\Omega(u)}\Big).$$
	\end{lem}
	\bpf
	We first show that for each $x\in \alpha$, it holds
	\be\label{eq:basic on k} 
	d_\Omega(x)\geq \max \left\{\frac{2\ell(\alpha[u,x])+d_\Omega(x)}{4c},\frac{1}{2c}d_\Omega(u)\right\}.
	\ee 
	
	Indeed, if $\ell(\alpha[u,x])\geq \frac{1}{2}d_{\Omega}(u)$, then $$d_\Omega(u)\leq 2\ell(\alpha[u,x])\leq 2c d_\Omega(x),$$
	and so,
	\[
	\max \left\{\frac{2\ell(\alpha[u,x])+d_\Omega(u)}{4c},\frac{1}{2c}d_\Omega(u)\right\}=\frac{2\ell(\alpha[u,x])+d_\Omega(u)}{4c}\leq d_\Omega(x). 
	\] 
	
	If $\ell(\alpha[u,x])< \frac{1}{2}d_{\Omega}(u)$, then $$d_\Omega(x)\geq d_\Omega(u)-\ell(\alpha[u,x])>\frac{1}{2}d_\Omega(u).$$
	Since $c\geq 1$, we infer from the above estimate that
	$$d_\Omega(x)\geq \frac{1}{2c}d_\Omega(u)=\max \left\{\frac{2\ell(\alpha[u,x])+d_\Omega(u)}{4c},\frac{1}{2c}d_\Omega(u)\right\}.$$
	
	In either case, we have proved \eqref{eq:basic on k}.
	
	For each $w\in \alpha$, by \eqref{eq:basic on k}, we have 
	\beqq
	\begin{aligned}
		k_{\Omega}(u,w)&\leq \ell_k(\alpha[u,w])=\int_{x\in\alpha[u,w]}\frac{ds}{d_{\Omega}(x)}\leq 4c\int_{x\in\alpha[u,w]}\frac{ds}{2\ell(\alpha[u,x])+d_\Omega(u)}\\ &\leq 2c \log\Big(1+\frac{2\ell(\alpha[u,w])}{d_\Omega(u)}\Big)\leq 4c \log\Big(1+\frac{\ell(\alpha[u,w])}{d_\Omega(u)}\Big).
	\end{aligned}
	\eeqq
	\epf

	\begin{lem}\label{Lemma-2.1}
		Let $(X, d)$ be a length space, $\Omega$ a proper subdomain of $X$ and $a>1$ a constant. For any points $x_1$ and $x_2$ in $\Omega$, if $d(x_1,x_2)\leq a^{-1}d_{\Omega}(x_1)$, then
		\beq\label{2025-11-12-1}
		k_{\Omega}(x_1,x_2)\leq \frac{9a}{10(a-1)}\frac{d(x_1,x_2)}{d_{\Omega}(x_1)}\leq \frac{10}{9}(a-1)^{-1}
		\eeq
		and
		\beq\label{2025-11-12-2}
		\ell([x_1,x_2])\leq \frac{10a}{9(a-1)}e^{\frac{10}{9}(a-1)^{-1}}d(x_1,x_2).
		\eeq
	\end{lem}
	
	\bpf Since $(X,d)$ is a length space,  for each $\varepsilon\in (0,(9a+1)^{-1})$, there exists some rectifiable curve $\alpha=\alpha_\varepsilon$ in $X$ connecting $x_1$ and $x_2$ such that
	\beq\label{20-7-25-1} \ell(\alpha)\leq \big(1+(a-1)\varepsilon\big)d(x_1,x_2).
	\eeq

	We claim that $\alpha\subset \Omega.$ Indeed, if not, then there exists some point $z\in \alpha\cap \partial \Omega$, and thus, it follows from \eqref{20-7-25-1} and the assumption $d(x_1,x_2)\leq a^{-1}d_{\Omega}(x_1)$ that
	$$d_{\Omega}(x_1)\leq  d(z,x_1)\leq \ell(\alpha)\leq \big(1+(a-1)\varepsilon\big)d(x_1,x_2)\leq \frac{\big(1+(a-1)\varepsilon\big)}{a}d_{\Omega}(x_1).$$
	It follows that $\varepsilon\geq 1$, which clearly contradicts with our choice of $\varepsilon$.
	
	Let $x\in \alpha$. Since $d(x_1,x_2)\leq a^{-1}d_{\Omega}(x_1)$, \eqref{20-7-25-1} gives
	$$
	d_{\Omega}(x) \geq d_{\Omega}(x_1)-\ell(\alpha)\geq a^{-1}(a-1)(1-\varepsilon)d_{\Omega}(x_1).
	$$
	Then it follows from the above estimate and our assumption $d(x_1,x_2)\leq a^{-1}d_{\Omega}(x_1)$ that
	
	\beq\label{20-7-25-2}
	\begin{aligned}
		\log \Big(1+\frac{\ell([x_1,x_2])}{d_{\Omega}(x_1)}\Big) &\stackrel{\eqref{(2.1)}}{\leq}  k_{\Omega}(x_1,x_2)\leq \int_{\alpha}\frac{|dx|}{d_{\Omega}(x)}\\
		&\leq \int_{\alpha}\frac{a|dx|}{(a-1)(1-\varepsilon)d_{\Omega}(x_1)}
		\stackrel{\eqref{20-7-25-1}}{\leq} \frac{a\big(1+(a-1)\varepsilon\big)}{(a-1)(1-\varepsilon)}\cdot\frac{d(x_1,x_2)}{d_{\Omega}(x_1)}\\
		&\leq \frac{10a}{9(a-1)}\cdot \frac{d(x_1,x_2)}{d_{\Omega}(x_1)}\leq \frac{10}{9}(a-1)^{-1},
	\end{aligned}	
	\eeq which proves the estimates in \eqref{2025-11-12-1}.
	
	Since the inequality $\lambda t\leq \log(1+t)$ holds in $[0, e^{-\frac{10}{9}(a-1)^{-1}}]$, where $\lambda=e^{-\frac{10}{9}(a-1)^{-1}}$, and since \eqref{20-7-25-2} gives
	$$\frac{\ell([x_1,x_2])}{d_{\Omega}(x_1)}\leq e^{\frac{10}{9}(a-1)^{-1}}-1,$$
	we know that
	$$\frac{\ell([x_1,x_2])}{d_{\Omega}(x_1)}\leq e^{\frac{10}{9}(a-1)^{-1}}\log \Big(1+\frac{\ell([x_1,x_2])}{d_{\Omega}(x_1)}\Big).$$
	
	Also, \eqref{20-7-25-2} gives $\log \Big(1+\frac{\ell([x_1,x_2])}{d_{\Omega}(x_1)}\Big)\leq \frac{10a}{9(a-1)}\cdot \frac{d(x_1,x_2)}{d_{\Omega}(x_1)}$, and thus, we get $$\ell([x_1,x_2])\leq \frac{10a}{9(a-1)}e^{\frac{10}{9}(a-1)^{-1}}d(x_1,x_2).$$ This shows that the estimate in \eqref{2025-11-12-2} is true, and hence, the proof of lemma is complete.
	\epf

	\subsection{Some elementary estimates}\label{subsec 2-2}
	
		Recall that $\Lambda_{xy}(\Omega)$ denotes the set of all quasihyperbolic geodesics in $\Omega$ with end points $x$ and $y$ and $\gamma_{xy}$ an element in $\Lambda_{xy}(\Omega)$. 
	
	The following simple estimate in Gromov hyperbolic domains will be frequently used in our later proofs.
	\blem\label{lem-2-3.0}
	Suppose that $\Omega$ is $\delta$-Gromov hyperbolic, and  $x, y, z \in \Omega$  are distinct points. For $w \in \gamma_{x y}$, if  $k_{\Omega}(w, y) \geq  2\delta+k_{\Omega}(y, z)$, then there must exist some point $u \in \gamma_{x z}$ such that $k_{\Omega}(w, u) \leq \delta$.
	\elem
	\bpf
	Since $\Omega$ is $\delta$-Gromov hyperbolic, there exists some $u \in \gamma_{xz} \cup \gamma_{yz}$  such that
	$$
	k_{\Omega}(w, u) \leq \delta.
	$$
	If $u \in \gamma_{yz}$, then $$k_{\Omega}(y, u) \geq k_{\Omega}(w, y)-k_{\Omega}(w, u)\geq k_{\Omega}(y, z)+\delta\geq k_{\Omega}(y,u)+\delta,$$
	which is impossible. Thus $u \in \gamma_{x z}$, and hence, the proof of lemma is complete.
	\epf
	
	We shall also need the following elementary estimate in $Q$-doubling metric spaces.
	\begin{lem}[{\cite[Lemma 4.1.11]{HKST-2015}}]\label{qs-5}
		Let $X=(X,d)$ be a $Q$-doubling metric space and $\Omega\subset X$ a domain.	Fix $R>0$ and $a\geq 1$, and let $r=\frac{R}{a}$. Then for any $x\in \Omega$, the ball $\mathbb{B}(x,R)$ contains
		at most $b$ balls with radius $r$ such that they are disjoint from each other, where $b\leq Q^{[\log_2 a]}$. Here and hereafter, $[\cdot]$ means the greatest integer part.
	\end{lem}

	\section{Gromov hyperbolicity implies ball separation condition}\label{sec-4-3}
	In this section, we shall prove the only if part of Theorem \ref{thm:new main metric space}, and thus, we assume throughout this section that $(X,d)$ is a $Q$-doubling length space, and $(\Omega, k)$ is $\delta$-Gromov hyperbolic with $\delta=C\geq 8$. The main result of this section is given as follows.
	
	\bthm\label{thm-24-5.1}
	Suppose that $(\Omega,k)$  is $C$-Gromov hyperbolic. Then $(\Omega,\sigma)$ satisfies the $\tau$-ball separation condition with $\tau=e^{(2CQ)^{192(1+C)}}$.
	\ethm

		\subsection{Three new classes of curves}\label{sec-4-1}
	In this section, we use the notation $\gamma_{xx_0y}\in \Lambda_{xy}(\Omega)$ to denote  a rectifiable curve $\gamma_{xy}\in \Lambda_{xy}(\Omega)$ with $x_0\in \gamma_{xy}$.  An important technical step towards {the proof of Theorem \ref{thm-24-5.1} is to }introduce three new classes of curves and derive a couple of fundamental lemmas about quasihyperbolic distances related to points on these curves.
	The first class is as follows.
	\begin{defn}[Class $P_{\alpha}^\gamma$]\label{def:class P}
		For given $\gamma_{xx_0y}\in \Lambda_{xy}(\Omega)$, $\alpha\in \Gamma_{xy}(\Omega)$ and $\theta>0$, the class $P_{\alpha}^{\gamma_{xx_0y}}(\theta)$ consists of all rectifiable curves $\gamma$ in $\Omega$ satisfying the following properties:
		\begin{enumerate}
			\item
			There is a point $z\in \alpha$ such that $\gamma\in\Lambda_{xz}(\Omega)$.
			\item
			There exists $z_0\in \gamma$ such that $k_{\Omega}(x_0,z_0)=k_{\Omega}(x_0,\gamma)\geq \theta$.
			\item
			For each $w\in \alpha(z,y]$ and for each $\gamma_{xw}\in \Lambda_{xw}(\Omega)$, it holds $k_{\Omega}(x_0,\gamma_{xw})< \theta$.
		\end{enumerate}
		
		To emphasize the points $z_0$ and $z$ on $\gamma$, we shall write $\gamma=\gamma_{xz_0z}$ for a general curve in $ P_{\alpha}^{\gamma_{xx_0y}}(\theta)$.
	\end{defn}

	The following basic result gives a sufficient condition for $P_{\alpha}^{\gamma_{xx_0y}}(\theta)$ to be nonempty.
	
	\blem\label{Lemma4-1.0}
	Suppose that $\gamma_{xx_0y}\in \Lambda_{xy}(\Omega)$, $\alpha\in \Gamma_{xy}(\Omega)$ and $\theta>1+C$. If $k_{\Omega}(x,x_0)>1+\theta$, then there exist a point $z\in\alpha$, a rectifiable curve $\gamma\in\Lambda_{xz}(\Omega)$ and a point $z_0\in \gamma$ such that $\gamma=\gamma_{xz_0z}\in P_{\alpha}^{\gamma_{xx_0y}}(\theta)$.
	\elem
	\bpf
	Take $v_1\in\alpha$ with $k_{\Omega}(v_1,y)\leq 1$, and then, fix $\gamma_{xv_1}\in\Lambda_{xv_1}(\Omega)$ and $\gamma_{v_1y}\in\Lambda_{v_1y}(\Omega)$. Since $(\Omega,k)$ is $C$-Gromov hyperbolic, there exists some point $v_2\in\gamma_{xv_1}\cup\gamma_{v_1y}$
	such that $$k_{\Omega}(x_0, v_2)\leq C.$$
	It follows from the above estimate and the triangle inequality that
	\be\label{eq-19-5}
	\begin{aligned}
		k_{\Omega}(x_0,\gamma_{xv_1} )&\leq  \max\{k_{\Omega}(x_0,v_2), k_{\Omega}(v_1,v_2)+k_{\Omega}(x_0, v_2)\}\\
		&\leq k_{\Omega}(v_1,y)+k_{\Omega}(x_0, v_2)\leq 1+C.
	\end{aligned}
	\ee
	
	Next, select $u_1\in\alpha$ with $k_{\Omega}(u_1,x)\leq 1$, and then, let $\gamma_{xu_1}\in\Lambda_{xu_1}(\Omega)$.
	Since $k_{\Omega}(x,x_0)>1+\theta$, we obtain from our choice of $u_1$ that
	$$k_{\Omega}(x_0,\gamma_{xu_1})\geq k_{\Omega}(x_0,x)-k_{\Omega}(x,u_1)>\theta>1+C.$$
	
	Finally, based on the previous estimate and \eqref{eq-19-5}, we may choose $z$ to be the last point  on $\alpha$ along the direction from $x$ to $y$ such that there exists some $\gamma=\gamma_{xz}\in \Lambda_{xz}(\Omega)$ with $k_{\Omega}(x_0,\gamma){\geq}\theta$. Clearly, there is $z_0\in \gamma$ such that $k_{\Omega}(x_0,z_0)=k_{\Omega}(x_0,\gamma)$. The above discussions show that $\gamma=\gamma_{xz_0z}\in P_{\alpha}^{\gamma_{xx_0y}}(\theta)$, and hence, the lemma is proved.
	\epf
	
	\blem\label{Lemma4-1.1} Suppose that  $\gamma_{xx_0y}\in \Lambda_{xy}(\Omega)$, $\alpha\in \Gamma_{xy}(\Omega)$ and  $\gamma_{xz_0z}\in P_{\alpha}^{\gamma_{xx_0y}}(\theta)$. Then
	$$\theta\leq k_{\Omega}(x_0,z_0)<1+\theta+C.$$
	\elem
	\bpf Let $w_1\in\alpha[z,y]$ be such that \be\label{H25-03-0}k_{\Omega}(w_1,z)\leq 1.\ee
	
	Since $\gamma_{xz_0z}\in P_{\alpha}^{\gamma_{xx_0y}}(\theta)$ and $w_1\in\alpha[z,y]$, by Definition \ref{def:class P}(3), for each $\gamma_{xw_1}\in \Lambda_{xw_1}(\Omega)$, there exists some point
	$y_0\in\gamma_{xw_1}$ such that
	\be\label{H25-03-1}k_{\Omega}(x_0, y_0)< \theta.\ee
	
	As $(\Omega,k)$ is $C$-Gromov hyperbolic, there exists some point $v_1\in\gamma_{xz}\cup\gamma_{zw_1}$
	such that $$k_{\Omega}(y_0, v_1)\leq C.$$
	Moreover, by the triangle inequality and the above estimate, we have
	$$
	k_{\Omega}(y_0,\gamma_{xz} )\leq k_{\Omega}(w_1,z)+k_{\Omega}(y_0, v_1)\stackrel{\eqref{H25-03-0}}{\leq}1+C.
	$$
	Then it follows that
	$$\theta\stackrel{\text{Definition }\ref{def:class P}(2)}{\leq}k_{\Omega}(x_0,z_0)=k_{\Omega}(x_0,\gamma_{xz})\leq k_{\Omega}(x_0, y_0)+k_{\Omega}(y_0,\gamma_{xz} )\stackrel{\eqref{H25-03-1}}{<}1+\theta+C,$$
	which is what we need.
	\epf
	
	Next, we introduce the second new class of curves.
	\begin{defn}[Class $O_{\alpha}^{\gamma}$]\label{def:Class O}
		For given $\gamma_{xx_0y}\in \Lambda_{xy}(\Omega)$, $\alpha\in \Gamma_{xy}(\Omega)$, $z\in\alpha$, $w\in\alpha[z,y]$ and $\vartheta>0$, the class $O_{\alpha[z,w]}^{\gamma_{xx_0y}}(\vartheta)$ consists of all rectifiable curves {$\gamma\in \Lambda_{zw}(\Omega)$} so that
		there exists $y_0\in\gamma$ satisfying
		$$k_{\Omega}(x_0,y_0)=k_{\Omega}(x_0,\gamma)\leq \vartheta.$$
		
		To emphasize the points $z$, $y_0$ and $w$ on $\gamma$, we shall write $\gamma=\gamma_{zy_0w}$ for a general curve in $O_{\alpha[z,w]}^{\gamma_{xx_0y}}(\vartheta)$.
	\end{defn}
	
	The following lemma is very fundamental in our later proofs.
	
	\blem\label{lem-23-4.1}
	Fix  $\gamma_{xz_1y}\in \Lambda_{xy}(\Omega)$, $\alpha\in \Gamma_{xy}(\Omega)$, $\gamma_{x z_{1,1} y_{1}} \in P^{\gamma_{ xz_{1}y}}_{\alpha}(\theta)$, $z_2\in \gamma_{xy}[z_1,y]$ and $y_2\in \alpha[y_1,y]$ with $\gamma_{x z_{1,2} y_{2}} \in P_{\alpha}^{\gamma_{xz_{2}y}}(\theta)$. Suppose that $\theta>2C$, $\gamma_{y_1z_{2,1}y_2}\in O_{\alpha[y_1,y_2]}^{\gamma_{xz_1y}}(2C)$   and $\gamma_{y_{1}z_{1,3}y_{3}} \in P^{\gamma_{ y_{1}z_{3}y_{2}}}_{\alpha[y_{1},y_{2}]}(\theta)$  for some $z_3\in \gamma_{y_1y_2}[y_1,z_{2,1}]$. If $k_{\Omega}(z_{1,1},z_{3})\geq 2(1+\theta+4C)$, then there exists some  $z_{2,3}\in \gamma_{x y_{1}}[y_{1}, z_{1,1}]$ such that
	$$\theta-C\leq k_{\Omega}(z_{2,3}, \gamma_{y_{1}y_{3}}) < 1+\theta+2C\;\mbox{ and }\; k_{\Omega}(z_{2,3}, z_{1,3})<1+\theta+2C.$$
	\elem
	\bpf It follows from Lemma \ref{Lemma4-1.1} and $\gamma_{y_1z_{2,1}y_2}\in O_{\alpha[y_1,y_2]}^{\gamma_{xz_1y}}(2C)$ that
	$$
	k_{\Omega}(z_{1,1},z_{2,1}) \leq k_{\Omega}(z_{1,1}, z_{1})+k_{\Omega}(z_{1}, z_{2,1})<1+\theta+3C.
	$$
	This, together with our assumption $k_{\Omega}(z_{1,1},z_{3})\geq 2(1+\theta+4C)$, gives
	\be\label{Heq-24-5}
	k_{\Omega}(z_{3}, z_{2,1}) \geq k_{\Omega}(z_{1,1}, z_{3})-k_{\Omega}(z_{1,1}, z_{2,1})> 1+\theta+5C>2C+k_{\Omega}(z_{1,1}, z_{2,1}).\ee
	
	Since $z_{2,1}\in\gamma_{y_1y_2}[y_2,z_3]$, by \eqref{Heq-24-5}, we may apply Lemma \ref{lem-2-3.0} (with $x=y_1,y=z_{2,1}, z=z_{1,1}$ and $w=z_3$) to find a point $z_{2,3}\in \gamma_{xy_{1}}[y_{1}, z_{1,1}]$  such that
	$$
	k_{\Omega}(z_{3}, z_{2,3}) \leq C.
	$$
	Moreover, by Lemma \ref{Lemma4-1.1}, we have
	$$\theta \leq k_{\Omega}(z_{3}, z_{1,3})=k_{\Omega}(z_{3}, \gamma_{y_{1}y_{3}}) < 1+\theta+C.$$
	
	Combining the above two estimates gives
	$$
	1+\theta+2C>k_{\Omega}(z_{3}, \gamma_{y_{1}y_{3}})+k_{\Omega}(z_{3}, z_{2,3})\geq k_{\Omega}(z_{2,3}, \gamma_{y_{1}y_{3}}) \geq k_{\Omega}(z_{3}, \gamma_{y_{1}y_{3}})-k_{\Omega}(z_{3}, z_{2,3})\geq\theta-C
	$$
	and
	$$
	k_{\Omega}(z_{2,3}, z_{1,3}) \leq k_{\Omega}(z_{3}, z_{1,3})+k_{\Omega}(z_{3}, z_{2,3})<1+\theta+2C.
	$$
	The proof of the lemma is complete.
	\epf
	\medskip

	Finally, we introduce the third new class of curves.
	\begin{defn}[Class $Q_{\alpha}^\gamma$]\label{def:Class AC}
		Fix $\gamma_{xx_1y}\in \Lambda_{xy}(\Omega)$ and $\alpha\in\Gamma_{xy}(\Omega)$. For  $z\in\alpha$ and $\gamma_{xz}\in \Lambda_{xz}(\Omega)$, if there exists $z_1\in\gamma_{xz}$ such that $2C\leq k_{\Omega}(x_1,\gamma_{xz_1z})=k_{\Omega}(x_1,z_1)\leq 7C$, then we write   $\gamma_{xz_1z}\in Q_{\alpha}^{\gamma_{xx_1y}}$.
	\end{defn}
	
	\blem\label{lem-22-H} Suppose that $\gamma_{xx_1y}\in \Lambda_{xy}(\Omega)$ and $\alpha\in \Gamma_{xy}(\Omega)$. If $\gamma_{xz_1z}\in P_{\alpha}^{\gamma_{xx_1y}}(3C)$, then
	$\gamma_{xz_1z}\in Q_{\alpha}^{\gamma_{xx_1y}}$.\elem
	\bpf This follows directly from Lemma \ref{Lemma4-1.1}.
	\epf
	\medskip
	
	We remark that several elementary estimates related to these new curve families are given in Appendix \ref{appendix-A}.
	
	\subsection{Proof of Theorem \ref{thm-24-5.1}}
	
	In this subsection, we shall prove Theorem \ref{thm-24-5.1}. Namely, we show that for any $\gamma_{xy}\in \Lambda_{xy}(\Omega)$, $z\in \gamma_{xy}$ and $\alpha\in \Gamma_{xy}(\Omega)$,
	$$B_{\sigma}(z,\tau d_{\Omega}(z))\cap \alpha\neq \emptyset,$$
	where $\tau=e^{(2CQ)^{192(1+C)}}$. For notational simplicity, we write $\gamma=\gamma_{xy}$, and fix $\alpha\in \Gamma_{xy}(\Omega)$.
	
	We shall prove it via a contradiction argument. Suppose, on the contrary, that there exists some point $x_{0,0}\in \gamma$ such that
	\be\label{eq-24-1.1}
	\sigma(x_{0,0}, \alpha)>\tau d_{\Omega}(x_{0,0}).
	\ee
	Set $N_1=[\frac{1}{32C^3} \log\tau]$. Then 
	\beqq
	e^{32CN_1}\leq\tau^{\frac{1}{C^2}}<\tau.
	\eeqq
	Let $x_{1,0}=y_{1,0}=x$, $x_{1,N_1+1}=x_{0,0}$, $M_0=2N_1$ and $M_1=[e^{-4-32C}\cdot  e^{\frac{1}{2}\log_Q \frac{N_1}{4}}]$. 
	
	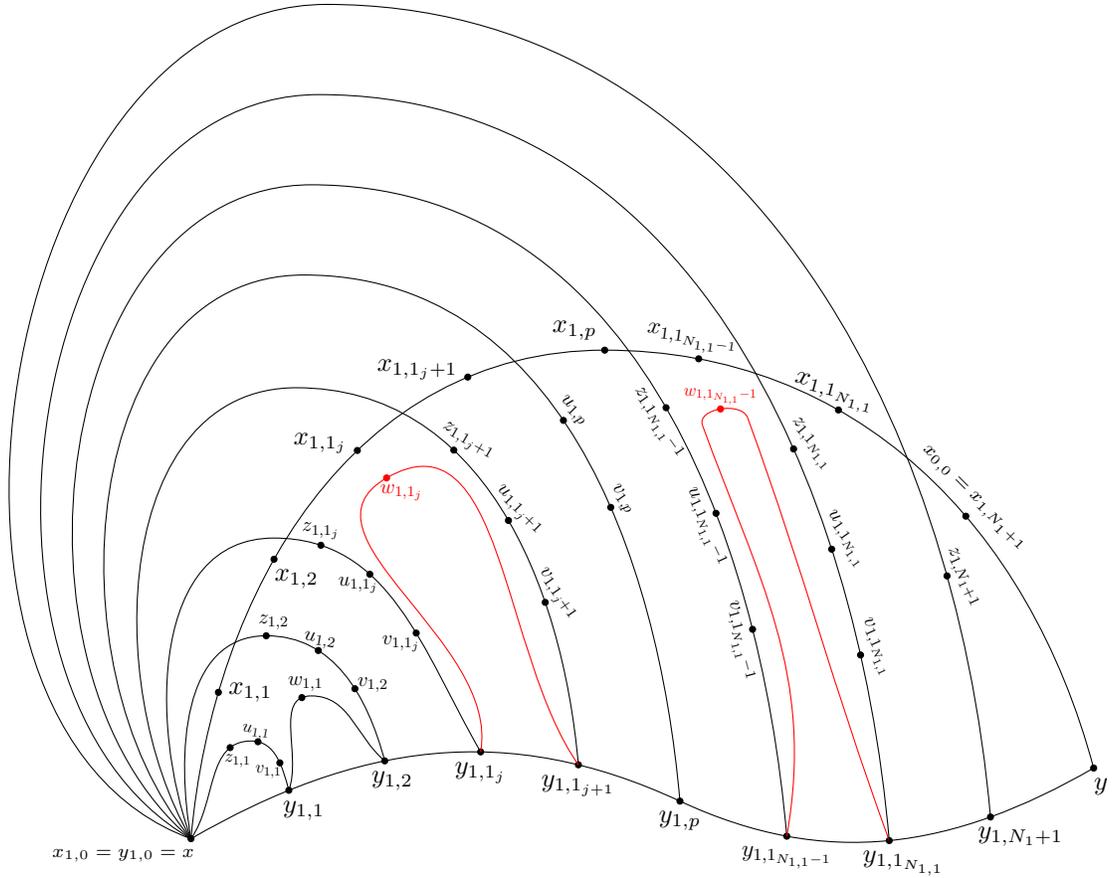
\begin{figure}[htbp]
		\begin{center}
			\begin{tikzpicture}[scale=1]
				
				\coordinate (x) at (-6,-0.5);
				\coordinate (y) at (6,0.44);
				
				\draw (x) to [out=30,in=155]
				coordinate[pos=0.2] (y11)
				coordinate[pos=0.4] (y12)
				coordinate[pos=0.6] (y11j)
				coordinate[pos=0.8] (y11j+1)
				coordinate[pos=1] (y1p)
				(0.5,0) to [out=-25,in=210]
				coordinate[pos=0.25] (y11n11-1)
				coordinate[pos=0.5] (y11n11)
				coordinate[pos=0.75] (y1n1+1)
				(y);
				
				\filldraw  (x)node[below,xshift=0.4cm,xshift=-1.3cm] {\tiny $x_{1,0}=y_{1,0}=x$} circle (0.04);
				\filldraw  (y11)node[below,xshift=0.2cm] {\small $y_{1,1}$} circle (0.04);
				\filldraw  (y12)node[below,xshift=0.1cm] {\small $y_{1,2}$} circle (0.04);
				\filldraw  (y11j)node[below,xshift=0cm] {\small $y_{1,1_j}$} circle (0.04);
				\filldraw  (y11j+1)node[below,xshift=0cm] {\small $y_{1,1_{j+1}}$} circle (0.04);
				\filldraw  (y1p)node[below] {\small $y_{1,p}$} circle (0.04);
				
				\filldraw  (y11n11-1)node[below,xshift=0cm] {\scalebox{0.8}[0.8]{$y_{1,1_{N_{1,1}-1}}$}} circle (0.04);
				\filldraw  (y11n11)node[below,xshift=0.2cm] {\small $y_{1,1_{N_{1,1}}}$} circle (0.04);
				\filldraw  (y1n1+1)node[below,xshift=0.4cm] {\small $y_{1,N_{1}+1}$} circle (0.04);
				\filldraw  (y)node[below,xshift=0.1cm] {\small $y$} circle (0.04);
				
				\draw  (x) to [out=85,in=180]
				coordinate[pos=0.2] (x11)
				coordinate[pos=0.4] (x12)
				coordinate[pos=0.6] (x11j)
				coordinate[pos=0.8] (x11j+1)
				coordinate[pos=1] (x1p)
				(-0.5,6) to [out=0,in=105]
				coordinate[pos=0.13] (x11n11-1)
				coordinate[pos=0.35] (x11n11)
				coordinate[pos=0.6] (x1n1+1)
				(y);
				
				\filldraw  (x11)node[right,xshift=0cm] {\small $x_{1,1}$}  circle (0.04);
				\filldraw  (x12)node[below,xshift=0.3cm] {\small $x_{1,2}$}  circle (0.04);
				\filldraw  (x11j)node[left,yshift=0.1cm] {\small $x_{1,1_j}$}  circle (0.04);
				\filldraw  (x11j+1)node[left,yshift=0.1cm] {\small $x_{1,1_j+1}$}  circle (0.04);
				\filldraw  (x1p)node[above,xshift=-0.4cm] {\small $x_{1,p}$}   circle (0.04);
				
				\filldraw  (x11n11-1)node[above,rotate=-8,xshift=-0.12cm,yshift=-0.05cm] {\scalebox{0.8}[0.8]{$x_{1,1_{N_{1,1}-1}}$}}  circle (0.04);
				\filldraw  (x11n11)node[above,rotate=-15,xshift=-0.1cm,yshift=-0.1cm] {\small $x_{1,1_{N_{1,1}}}$}  circle (0.04);
				\filldraw  (x1n1+1)node[above,rotate=-45,xshift=-0.1cm,yshift=0cm] {\tiny $x_{0,0}=x_{1,N_{1}+1}$}  circle (0.04);
				
				\draw  (x) to [out=55,in=180]
				coordinate[pos=0.8] (z11)
				(-5.2,0.8) to [out=0,in=105]
				coordinate[pos=0.1] (u11)
				coordinate[pos=0.6] (v11)
				(y11);
				
				\filldraw  (z11)node[below,rotate=0,xshift=0.1cm,yshift=0.05cm] {\scalebox{0.6}[0.6]{$z_{1,1}$}}  circle (0.04);
				\filldraw  (u11)node[above,rotate=-10,xshift=-0.05cm,yshift=-0.1cm] {\scalebox{0.6}[0.6]{$u_{1,1}$}}  circle (0.04);
				\filldraw  (v11)node[below,rotate=0,xshift=-0.15cm,yshift=0.1cm] {\scalebox{0.6}[0.6]{$v_{1,1}$}}  circle (0.04);
				
				\draw  (x) to [out=100,in=180]
				coordinate[pos=1] (z12)
				(-5,2.2) to [out=0,in=105]
				coordinate[pos=0.3] (u12)
				coordinate[pos=0.6] (v12)
				(y12);
				
				\filldraw  (z12)node[above,rotate=0,xshift=0.1cm,yshift=-0.05cm] {\scalebox{0.7}[0.7]{$z_{1,2}$}}  circle (0.04);
				\filldraw  (u12)node[above,rotate=-10,xshift=0cm,yshift=-0.1cm] {\scalebox{0.7}[0.7]{$u_{1,2}$}}  circle (0.04);
				\filldraw  (v12)node[right,rotate=0,xshift=-0.1cm,yshift=0.05cm] {\scalebox{0.7}[0.7]{$v_{1,2}$}}  circle (0.04);
				
				\draw  (y11) to [out=70,in=180]
				coordinate[pos=0.9] (w11)
				(-4.4,1.4) to [out=0,in=140]
				(y12);
				
				\filldraw  (w11)node[above,rotate=0,xshift=0.05cm,yshift=-0.05cm] {\scalebox{0.7}[0.7]{$w_{1,1}$}}  circle (0.04);
				
				\draw  (x) to [out=110,in=180]
				(-4.9,3.5) to [out=0,in=120]
				coordinate[pos=0.15] (z11j)
				coordinate[pos=0.35] (u11j)
				coordinate[pos=0.6] (v11j)
				(y11j);
				
				\filldraw  (z11j)node[above,rotate=0,xshift=0cm,yshift=-0.05cm] {\scalebox{0.7}[0.7]{$z_{1,1_j}$}}  circle (0.04);
				\filldraw  (u11j)node[below,rotate=0,xshift=-0.15cm,yshift=0.1cm] {\scalebox{0.7}[0.7]{$u_{1,1_j}$}}  circle (0.04);
				\filldraw  (v11j)node[below,rotate=0,xshift=-0.2cm,yshift=0.1cm] {\scalebox{0.7}[0.7]{$v_{1,1_j}$}}  circle (0.04);
				
				\draw  (x) to [out=120,in=180]
				(-4.6,5.5) to [out=0,in=95]
				coordinate[pos=0.35] (z11j+1)
				coordinate[pos=0.53] (u11j+1)
				coordinate[pos=0.7] (v11j+1)
				(y11j+1);
				
				\filldraw  (z11j+1)node[above,rotate=-30,xshift=0.1cm,yshift=-0.05cm] {\scalebox{0.7}[0.7]{$z_{1,1_j+1}$}}  circle (0.04);
				\filldraw  (u11j+1)node[above,rotate=-45,xshift=0cm,yshift=-0.05cm] {\scalebox{0.7}[0.7]{$u_{1,1_j+1}$}}  circle (0.04);
				\filldraw  (v11j+1)node[above,rotate=-55,xshift=0cm,yshift=-0.05cm] {\scalebox{0.7}[0.7]{$v_{1,1_j+1}$}}  circle (0.04);
				
				\draw[red]  (y11j) to [out=80,in=210]
				coordinate[pos=1] (w11j)
				(-3.4,4.3) to [out=30,in=125]
				(y11j+1);
				
				\filldraw[red]  (w11j)node[below,rotate=0,xshift=0.2cm,yshift=0.05cm] {\scalebox{0.7}[0.7]{$w_{1,1_j}$}}  circle (0.04);
				
				\draw  (x) to [out=130,in=180]
				(-4.5,7) to [out=0,in=95]
				coordinate[pos=0.46] (u1p)
				coordinate[pos=0.6] (v1p)
				(y1p);
				
				\filldraw  (u1p)node[above,rotate=-45,xshift=0cm,yshift=-0.05cm] {\scalebox{0.7}[0.7]{$u_{1,p}$}}  circle (0.04);
				\filldraw  (v1p)node[above,rotate=-55,xshift=0cm,yshift=-0.05cm] {\scalebox{0.7}[0.7]{$v_{1,p}$}}  circle (0.04);
				
				\draw  (x) to [out=140,in=180]
				(-4.4,8.2) to [out=0,in=95]
				coordinate[pos=0.52] (z11n11-1)
				coordinate[pos=0.65] (u11n11-1)
				coordinate[pos=0.78] (v11n11-1)
				(y11n11-1);
				
				\filldraw  (z11n11-1)node[below,rotate=-55,xshift=0.1cm,yshift=0.05cm] {\scalebox{0.7}[0.7]{$z_{1,1_{N_{1,1}}-1}$}}  circle (0.04);
				\filldraw  (u11n11-1)node[below,rotate=-65,xshift=0.1cm,yshift=0.05cm] {\scalebox{0.7}[0.7]{$u_{1,1_{N_{1,1}}-1}$}}  circle (0.04);
				\filldraw  (v11n11-1)node[below,rotate=-75,xshift=0.1cm,yshift=0.05cm] {\scalebox{0.7}[0.7]{$v_{1,1_{N_{1,1}}-1}$}}  circle (0.04);

				\draw[red]  (y11n11-1) to [out=80,in=290]
				(0.8,5);
				\draw[red]  (0.8,5) to [out=110,in=110] coordinate[pos=0.5] (w11n11-1) (1.4,5.1);
				\draw[red]  (1.4,5.1) to [out=290,in=110] (y11n11);
				
				\filldraw[red]  (w11n11-1)node[above,rotate=0,xshift=0cm,yshift=-0.08cm] {\scalebox{0.6}[0.6]{$w_{1,1_{N_{1,1}}-1}$}}  circle (0.04);
				
				\draw  (x) to [out=150,in=180]
				(-4.3,9.4) to [out=0,in=95]
				coordinate[pos=0.63] (z11n11)
				coordinate[pos=0.73] (u11n11)
				coordinate[pos=0.83] (v11n11)
				(y11n11);
				
				\filldraw  (z11n11)node[above,rotate=-50,xshift=0.1cm,yshift=-0.05cm] {\scalebox{0.7}[0.7]{$z_{1,1_{N_{1,1}}}$}}  circle (0.04);
				\filldraw  (u11n11)node[above,rotate=-55,xshift=0cm,yshift=-0.05cm] {\scalebox{0.7}[0.7]{$u_{1,1_{N_{1,1}}}$}}  circle (0.04);
				\filldraw  (v11n11)node[above,rotate=-65,xshift=0cm,yshift=-0.05cm] {\scalebox{0.7}[0.7]{$v_{1,1_{N_{1,1}}}$}}  circle (0.04);
				
				\draw  (x) to [out=160,in=180]
				(-4.2,10.6) to [out=0,in=95]
				coordinate[pos=0.8] (z1n1+1)
				(y1n1+1);
				
				\filldraw  (z1n1+1)node[above,rotate=-65,xshift=0.1cm,yshift=-0.05cm] {\scalebox{0.7}[0.7]{$z_{1,N_1+1}$}}  circle (0.04);
				
			\end{tikzpicture}
		\end{center}
		\caption{Illustration for the proof of Lemma \ref{Lem4.12-1}} \label{fig-25-8-7}
	\end{figure}
	
	The following lemma will be crucial for the proof of Theorem \ref{thm-24-5.1}.
	\blem\label{Lem4.12-1}  For each $p\in\{1,\cdots,N_1\}$, there exists a sequence of successive points  $\{x_{1,p}\}\subset\gamma[x_{1,p-1},x_{0,0}]$ along the direction from $x_{1,p-1}$ to  $x_{0,0}$ such that $\sigma(x_{0,0}, x_{1,p})\leq e^{32C(1+N_1)} d_{\Omega}(x_{0,0})$. Moreover, the following conclusions hold:
	\ben
	\item\label{Lem4.12-1-1}
	For each $p\in\{0,\cdots,N_1\}$, $k_{\Omega}(x_{1,p},x_{1,p+1})>30C$.
	\item\label{Lem4.12-1-2}
	For each $p\in\{1,\cdots,N_1+1\}$, there exist $y_{1,p}\in\alpha[y_{1,p-1},y]$, $\gamma_{xy_{1,p}}\in \Lambda_{xy_{1,p}}(\Omega)$ and $z_{1,p}\in\gamma_{xy_{1,p}}$ such that
	$\gamma_{xz_{1,p}y_{1,p}}\in P_{\alpha}^{\gamma_{xx_{1,p}y}}(3C)$; and for each $p\in\{1,\cdots,N_1\}$ and $\gamma_{y_{1,p}y_{1,p+1}}\in \Lambda_{y_{1,p}y_{1,p+1}}(\Omega)$, there exists $w_{1,p}\in\gamma_{y_{1,p}y_{1,p+1}}$ such that $\gamma_{y_{1,p}w_{1,p}y_{1,p+1}}\in O_{\alpha[y_{1,p},y_{1,p+1}]}^{\gamma_{xx_{1,p}y}}(2C)$.
	
	\item\label{Lem4.12-1-3}
	There exists some integer $N_{1,1}>\frac{N_1}{2}$ such that
	\begin{itemize}
		\item there exists a sequence of integers  $\{1_j\}_{j=1}^{N_{1,1}}\subset \{1,\cdots,N_1-1\}$ with $1_j<1_{j+1}$ for each $j\in \{1,\cdots,N_{1,1}-1\}$.
		
		\item there exists a sequence $\{x_{2,jM_1-(t-1)}\}_{j,t}$, indexed with $j\in\{1,\cdots,N_{1,1}\}$ and $t\in\{1,\cdots,M_1\}$, of successive points on $\gamma_{y_{1,1_j}y_{1,1_j+1}}[y_{1,1_j},w_{1,1_j}]$ along the direction from $y_{1,1_j}$ to $w_{1,1_j}$ satisfying
		$$d_{\Omega}(x_{2,jM_1-(t-1)})\leq r_1,$$ where $r_1=e^{-\frac{1}{2}\log_Q \frac{N_1}{4}} \cdot e^{32C(N_1+2)+2}d_{\Omega}(x_{0,0}).$
		Moreover,  it holds
		$$\frac{3}{4}e^{32C(N_1+2)}d_{\Omega}(x_{0,0})< \sigma(x_{0,0},x_{2,jM_1-(t-1)})<\frac{3}{4}e^{32C(N_1+2)+2}d_{\Omega}(x_{0,0}).$$
		\item for each $j\in\{1,\cdots,N_{1,1}\}$ and $t\in\{2,\cdots,M_1\}$, $$k_{\Omega}(x_{2,jM_1-(t-1)}, x_{2,jM_1-(t-2)})>30C.$$
	\end{itemize}
	\een
	\elem
	
	\bpf  (1) For each $p \in\{1, \cdots, N_1\}$, let $x_{1,p} \in \gamma[x, x_{0,0}]$ be such that
	\be\label{eq-24-2}
	\sigma(x_{0,0}, x_{1,p})=e^{32C(N_1+1-p)} d_{\Omega}(x_{0,0}).
	\ee
	
	As $x_{1,N_1+1}=x_{0,0}$, it follows from \eqref{eq-24-2} that
	\be\label{H2504-02-1}
	d_{\Omega}(x_{1,p})\leq d_{\Omega}(x_{0,0})+\sigma(x_{0,0}, x_{1,p})=(1+e^{32C(N_1+1-p)}) d_{\Omega}(x_{0,0}).
	\ee
	
	For each $p\in\{1, \ldots, N_1\}$, we have
	\[
	\sigma(x_{1,p},x_{1,p+1})\geq \sigma(x_{1,p},x_{0,0})-\sigma(x_{1,p+1},x_{0,0})\stackrel{\eqref{eq-24-2}}{=}e^{32C(N_1-p)}(e^{32C}-1)d_{\Omega}(x_{0,0}),
	\]
	and then,
	\beq
	\begin{aligned}\label{eq-24-3}
		k_{\Omega}(x_{1,p}, x_{1,p+1}) &\stackrel{\eqref{(2.2)}}{\geq} \log\left(1+\frac{\sigma(x_{1,p}, x_{1,p+1})}{d_{\Omega}(x_{1,p+1})}\right) \\ &\stackrel{\eqref{eq-24-1.1}+\eqref{H2504-02-1}}{\geq} \log\left(1+\frac{e^{32C(N_1-p)}(e^{32C}-1)}{1+e^{32C(N_1-p)}}\right){>}30C.
	\end{aligned}
	\eeq
	This proves \eqref{Lem4.12-1-1}.
	
	(2) For each $p \in \{1, \dots, N_1+1\}$, we note that
	\[
	\sigma(x,x_{1,p})\geq \sigma(x,x_{0,0})-\sigma(x_{0,0},x_{1,p})\stackrel{\eqref{eq-24-1.1}+\eqref{eq-24-2}}{\geq} \tau d_{\Omega}(x_{0,0})-e^{32C(N_1+1-p)}d_{\Omega}(x_{0,0}),
	\]
	and so,
	\be\label{H25-0507-1}
	k_{\Omega}(x,x_{1,p})\stackrel{\eqref{(2.2)}}{\geq} \log\left(1+\frac{\sigma(x,x_{1,p})}{d_{\Omega}(x_{1,p})}\right)\stackrel{\eqref{H2504-02-1}}{\geq} \log\left(1+\frac{\tau-e^{32C(N_1+1-p)}}{1+e^{32C(N_1+1-p)}}\right)>\frac{1}{2}\log \tau>1+3C.
	\ee
	Then Lemma \ref{Lemma4-1.0} and \eqref{H25-0507-1} imply that  for each $p \in \{1, \dots, N_1+1\}$, there exist $y_{1,p}\in \alpha[y_{1,p-1},y]$, $\gamma_{xy_{1,p}}\in\Lambda_{xy_{1,p}}(\Omega)$ and $z_{1,p}\in \gamma_{xy_{1,p}}$ such that
	\be\label{eq-24-4}
	\gamma_{xz_{1,p}y_{1,p}}\in P_{\alpha}^{\gamma_{x x_{1,p}y}}(3C).
	\ee
	
	Moreover, for each $p \in \{1, \dots, N_1\}$,  {by \eqref{eq-24-3} and \eqref{eq-24-4}, we may apply Lemmas \ref{lem-22-H} and  \ref{lem-22-3.1}(\ref{Lemma4-1}) from Appendix \ref{appendix-A} (with $z=y_{1,p}$ and $w=y_{1,p+1}$) to infer that for each $\gamma_{y_{1,p}y_{1,p+1}}\in\Lambda_{y_{1,p}y_{1,p+1}}(\Omega)$, there exists $w_{1,p}\in \gamma_{y_{1,p}y_{1,p+1}}$ such that}
	\beqq
	\gamma_{y_{1,p}w_{1,p}y_{1,p+1}}\in O_{\alpha[y_{1,p}, y_{1,p+1}]}^{\gamma_{xx_{1,p}y}}(2C).
	\eeqq
	These prove \eqref{Lem4.12-1-2}.
	
	(3) For each $p \in \{1, \dots, N_1\}$, it follows from Lemma \ref{Lemma4-1.1} and \eqref{(2.2)} that
	$$
	\log \Big(1+\frac{\sigma(x_{1,p},z_{1,p})}{d_{\Omega}(x_{1,p})}\Big)\leq k_{\Omega}(x_{1,p},z_{1,p})<5C,
	$$
	and so,
	\[
	\sigma(x_{1,p},z_{1,p})\leq (e^{5C}-1)d_{\Omega}(x_{1,p})\stackrel{\eqref{H2504-02-1}}{\leq}(e^{5C}-1)(1+e^{32C(N_1+1-p)})d_{\Omega}(x_{0,0}).
	\]
	It follows from the above estimate and \eqref{eq-24-2} that
	\be\label{H25-04-06-1}\sigma(x_{0,0},z_{1,p})\leq \sigma(x_{0,0},x_{1,p})+\sigma(x_{1,p},z_{1,p})\leq e^{5C}(1+e^{32C(N_1+1-p)})d_{\Omega}(x_{0,0}).
	\ee
	
	Based on \eqref{eq-24-1.1} and \eqref{H25-04-06-1}, we may choose
	\beq\label{eq:selection of u1p}
	u_{1,p}\in \gamma_{xz_{1,p}y_{1,p}}[z_{1,p}, y_{1,p}]\cap \mathbb{S}_{\sigma}(x_{0,0}, e^{32C(N_1+2)}d_{\Omega}(x_{0,0})),
	\eeq
	and then, let \be\label{447}v_{1,p}\in \gamma_{xz_{1,p}y_{1,p}} [u_{1,p}, y_{1,p}] \cap \mathbb{S}_{\sigma}(x_{0,0}, e^{32C(N_1+2)+1} d_{\Omega}(x_{0,0}))\ee be the first point along the direction from $u_{1,p}$ to $y_{1,p}$.
	
	Next, we claim that there exists an integer
	\beqq
	N_{1,1}> \frac{N_1}{2}
	\eeqq
	such that
	\begin{itemize}
		\item there is a sequence of integers $\{1_j\}_{j=1}^{N_{1,1}}\subset \{1,\cdots,N_1-1\}$ with $1_j<1_{j+1}$ for each $j\in \{1,\cdots,N_{1,1}-1\}$.
		
		\item for each $j \in \{1, \dots, N_{1,1}\}$ and each $u \in \gamma_{xz_{1,1_{j}}y_{1,1_{j}}}[u_{1,1_{j}}, v_{1,1_{j}}]$,
		\be\label{eq-24-5}
		d_{\Omega}(u)\leq r_1=e^{-\frac{1}{2}\log_Q \frac{N_1}{4}}\cdot e^{32C(N_1+2)+2}d_{\Omega}(x_{0,0}).
		\ee
	\end{itemize}
	
	Indeed, suppose, on the contrary, that our claim fails. Then there exist $N_{1,2}(> \frac{N_1}{3})$ integers  $\rho_1, \dots, \rho_{N_{1,2}}$ in $\{1, \dots, N_1\}$ such that
	\begin{itemize}
		\item for each $t\in\{1, \dots, N_{1,2}-1\}$, $\rho_t<\rho_{t+1}$.
		\item  for each $t\in \{1, \dots, N_{1,2}\}$, there exists a point $u_{1,\rho_{t}}^1\in\gamma_{xz_{1,\rho_{t}}y_{1,\rho_{t}}}[u_{1,\rho_{t}}, v_{1,\rho_{t}}]$ so that
		\be\label{eq-24-7}
		d_{\Omega}(u_{1,\rho_{t}}^1) > r_1.
		\ee
	\end{itemize}
	
	Let $B_{0,0}=\mathbb{B}\big(x_{0,0},e^{32C(N_1+2)+2}d_{\Omega}(x_{0,0})\big)$. For each $t \in \{1, \dots, N_{1,2}\}$, we take
	$$B_{\rho_t}=\mathbb{B}\left(u_{1,\rho_{t}}^1, \frac{1}{3}r_1\right).$$
	Then by the choice of $u_{1,\rho_t}$ and $v_{1,\rho_t}$ in \eqref{eq:selection of u1p} and \eqref{447}, for each $u \in \overline{B_{\rho_t}}$, we have
	\begin{align*}
		d(x_{0,0}, u) &\leq d(x_{0,0}, u_{1,\rho_{t}}^1)+d(u,u_{1,\rho_{t}}^1)\leq e^{32C(N_1+2)+1}d_{\Omega}(x_{0,0})+\frac{1}{3}r_1\\
		&< e^{32C(N_1+2)+2}d_{\Omega}(x_{0,0}),
	\end{align*}
	and so, $$\overline{B_{\rho_t}}\subset B_{0,0}.$$
	
	If all these balls are disjoint, then applying Lemma \ref{qs-5} with $R=e^{32C(N_1+2)+2}d_{\Omega}(x_{0,0})$ and $r=\frac{1}{3}e^{-\frac{1}{2}\log_Q \frac{N_1}{4}} R$ gives
	\[
	N_{1,2}< \frac{N_1}{3},
	\]
	which is a contradiction.
	
	For the remaining case, there exist two integers $q_1<q_2\in \{1, \dots, N_{1,2}\}$ such that $B_{\rho_{q_1}}\cap B_{\rho_{q_2}}\not=\emptyset$. 
	It follows that
	$$d(u_{1,\rho_{q_1}},u_{1,\rho_{q_2}})\leq \frac{2}{3}r_1\stackrel{\eqref{eq-24-7}}{\leq}\frac{2}{3}\min\left\{d_{\Omega}(u_{1,\rho_{q_1}}),d_{\Omega}(u_{1,\rho_{q_2}})\right\},$$ and thus, by Lemma \ref{Lemma-2.1}, we have
	\be\label{H25-086-1}k_{\Omega}(u_{1,\rho_{q_1}},u_{1,\rho_{q_2}})\leq\frac{20}{9}<3C.\ee
	
	Note that by \eqref{eq:selection of u1p}, we have $$u_{1,\rho_{q_1}}\in\gamma_{xz_{1,\rho_{q_1}}y_{1,\rho_{q_1}}}[z_{1,\rho_{q_1}},y_{1,\rho_{q_1}}]\quad \text{and}\quad  u_{1,\rho_{q_2}}\in\gamma_{xz_{1,\rho_{q_2}}y_{1,\rho_{q_2}}}[z_{1,\rho_{q_2}},y_{1,\rho_{q_2}}].$$ Then it follows from the assertions  (\ref{Lem4.12-1-1}) and (\ref{Lem4.12-1-2}) of Lemma \ref{Lem4.12-1}, together with Lemma \ref{lem-22-3.1} from Appendix \ref{appendix-A}, that
	$$k_{\Omega}(u_{1,\rho_{q_1}},u_{1,\rho_{q_2}})\geq 3C,$$
	which clearly contradicts with \eqref{H25-086-1}. The proof of claim is thus complete.
	
	Let us continue the proof based on the above claim (i.e., \eqref{eq-24-5}). For each $j\in\{1,\dots,N_{1,1}\}$, let $x_{1,j}^1=u_{1,1_{j}}$, $y_{1,j}^1=y_{1,1_{j}}$ and $y_{1,j}^2=y_{1,1_{j}+1}$. For each $t\in\{2,\cdots,M_1\}$, let $x_{1,j}^t\in\gamma_{y_{1,1_{j}}z_{1,1_{j}}y_{1,1_{j}-1}}$ be such that
	\be\label{H25-04-06-2}
	\sigma(x_{1,j}^t,x_{1,j}^{t-1})=e^{32C}r_1.
	\ee
	Then we have
	\be\label{H25-04-06-3}
	k_{\Omega}(x_{1,j}^t,x_{1,j}^{t-1})\stackrel{\eqref{(2.1)}}{\geq} \log\Big(1+\frac{\sigma(x_{1,j}^t,x_{1,j}^{t-1})}{d_{\Omega}(x_{1,j}^{t-1})}\Big)\stackrel{\eqref{eq-24-5}}{\geq} \log\Big(1+\frac{\sigma(x_{1,j}^t,x_{1,j}^{t-1})}{r_1}\Big)\stackrel{\eqref{H25-04-06-2}}{>} 32C.
	\ee
	
	Note that by \eqref{H25-04-06-1} and \eqref{eq:selection of u1p}, we have $$\sigma(x_{1,j}^1,z_{1,1_{j}})\geq \sigma(x_{0,0},x_{1,j}^1)-\sigma(x_{0,0},z_{1,1_{j}})>e^{32C(N_1+2)-1}d_{\Omega}(x_{0,0}),$$
	from which it follows that
	$$k_{\Omega}(x_{1,j}^1,z_{1,1_{j}})\stackrel{\eqref{(2.1)}}{\geq}\log\Big(1+\frac{\sigma(x_{1,j}^1,z_{1,1_{j}})}{d_{\Omega}(x_{1,j}^1)}\Big)\stackrel{\eqref{eq-24-5}}{>}\frac{1}{2}\log_Q \frac{N_1}{4}-3{\geq 10C}.$$
	This, together with Lemmas \ref{lem-24-5.3}(\ref{lem-24-5.3-1}) and  \ref{Lem4.12-1}(\ref{Lem4.12-1-2}), shows that there exists some point $x_{1,j}^{1,1}\in \gamma_{y_{1,j}^1w_{1,1_{j}}y_{1,j}^2}[y_{1,j}^1,w_{1,1_{j}}]$
	such that
	\be\label{H25-04-06-4}
	k_{\Omega}(x_{1,j}^1,x_{1,j}^{1,1})\leq C.
	\ee
	
	Next, we observe the following iteration: For each $t\in\{1,\cdots,M_1-1\}$, if there exists some point $x_{1,j}^{1,t}\in  \gamma_{y_{1,j}^1w_{1,1_{j}}y_{1,j}^2}[y_{1,j}^1,w_{1,1_{j}}]$
	such that $$k_{\Omega}(x_{1,j}^t, x_{1,j}^{1,t})\leq C,$$
	then there exists some point $x_{1,j}^{1,t+1} \in\gamma_{y_{1,j}^1w_{1,1_{j}}y_{1,j}^2}[y_{1,j}^1,x_{1,j}^{1,t}]$ such that
	$$ k_{\Omega}(x_{1,j}^{t+1},x_{1,j}^{1,t+1})\leq C.$$
	
	Indeed, note that
	$$k_{\Omega}(x_{1,j}^{t+1},x_{1,j}^t)\stackrel{\eqref{H25-04-06-3}}>32C\geq 31C+k_{\Omega}(x_{1,j}^t, x_{1,j}^{1,t}).$$
	The observation follows directly from  Lemma \ref{lem-2-3.0} (with $x=y_{1,j}^1, y=x_{1,j}^t, z=x_{1,j}^{1,t}$ and $w=x_{1,j}^{t+1}$).
	
	The above observation, together with \eqref{H25-04-06-4}, implies that for each $j\in\{1,\dots,N_{1,1}\}$ and $t\in\{1,\cdots,M_1\}$, there exists some point $x_{1,j}^{1,t+1} \in\gamma_{y_{1,j}^1w_{1,1_{j}}y_{1,j}^2}[y_{1,j}^1,x_{1,j}^{1,t}]$ such that
	\be\label{H25-04-06-5} k_{\Omega}(x_{1,j}^{t+1},x_{1,j}^{1,t+1})\leq C.\ee
	
	Note that
	\beq\label{H25-0512-1}
	\begin{aligned}
		\frac{3}{4}e^{32C(N_1+2)}d_{\Omega}(x_{0,0})&\stackrel{\eqref{eq:selection of u1p}+\eqref{H25-04-06-2}}{<}\sigma(x_{0,0},x_{1,j}^t)-\sigma(x_{1,j}^t,x_{1,j}^{1,t})\leq\sigma(x_{0,0},x_{1,j}^{1,t})\\ &\leq  \sigma(x_{0,0},x_{1,j}^t)+\sigma(x_{1,j}^t,x_{1,j}^{1,t})\stackrel{\eqref{eq:selection of u1p}+\eqref{H25-04-06-2}}{<}\frac{3}{4}e^{32C(N_1+2)+2}d_{\Omega}(x_{0,0}),
	\end{aligned}
	\eeq
	and thus,
	$$\max\left\{\log \frac{d_{\Omega}(x_{1,j}^{1,t+1})}{d_{\Omega}(x_{1,j}^{t+1})},\log \frac{\sigma(x_{1,j}^{1,t},x_{1,j}^{1,t+1})}{d_{\Omega}(x_{1,j}^{t+1})}\right\}\stackrel{\eqref{(2.1)}}{\leq} k_{\Omega}(x_{1,j}^{t+1},x_{1,j}^{1,t+1})\stackrel{\eqref{H25-04-06-5}}{\leq} C.$$
	It follows from the above estimate and \eqref{eq-24-5} that
	\be\label{H25-04-07-1} \max\{\sigma(x_{1,j}^{1,t},x_{1,j}^{1,t+1}),d_{\Omega}(x_{1,j}^{1,t+1})\}\leq e^C r_1.\ee
	
	Moreover, for each $j\in\{1,\dots,N_{1,1}\}$ and $t\in\{1,\cdots,M_1-1\}$, we have
	\be\label{H25-0508-2}
	k_{\Omega}(x_{1,j}^{1,t}, x_{1,j}^{1,t+1})\geq  k_{\Omega}(x_{1,j}^{t}, x_{1,j}^{t+1})-k_{\Omega}(x_{1,j}^{t}, x_{1,j}^{1,t})-k_{\Omega}(x_{1,j}^{t+1}, x_{1,j}^{1,t+1})\stackrel{\eqref{H25-04-06-3}+\eqref{H25-04-06-5}}{>}30C.
	\ee
	
	Finally, for each $j\in\{1,\dots,N_{1,1}\}$ and $t\in\{1,\cdots,M_1\}$, if we take $x_{2,jM_1-(t-1)}=x_{1,j}^{1,t}$, then it follows from \eqref{H25-0512-1} $\sim$ \eqref{H25-0508-2} that Lemma \ref{Lem4.12-1}(\ref{Lem4.12-1-3}) holds. The proof of the lemma is thus complete.
	\epf
	\medskip

	For each $i\in\{2,\cdots,M_0\}$, set $r_{i-1}=e^{-\frac{1}{2}\log_Q \frac{N_{i-1}}{4}}\cdot e^{20C(N_1+2)+2i-2}d_{\Omega}(x_{0,0})$ and  $N_i=(M_1-1) N_{i-1,1}$. Then an iteration of Lemma \ref{Lem4.12-1}  gives the following more general result.
	
	\begin{figure}[htbp]
		\begin{center}
			\begin{tikzpicture}[scale=1]
				
				\coordinate (x) at (-6,-0.5);
				\coordinate (y) at (6,0.44);
				
				\draw (x) to [out=-20,in=155]
				coordinate[pos=0.25] (y11j)
				coordinate[pos=0.5] (y22j)
				coordinate[pos=0.75] (yiij)
				coordinate[pos=1] (yiij+1)
				(0.5,0) to [out=-25,in=210]
				coordinate[pos=0.3] (y22j+1)
				coordinate[pos=0.6] (y11j+1)
				(y);
				
				\filldraw  (x)node[below,xshift=-0.05cm] {\small $x$} circle (0.04);
				\filldraw  (y11j)node[below,xshift=0.2cm] {\small $y_{1,1_j}$} circle (0.04);
				
				\filldraw  (y22j)node[below,xshift=0.1cm] {\small $y_{2,2_j}$} circle (0.04);
				\node[below,rotate=90, xshift=-0.57cm, yshift=0.15cm] at(y22j) {\scalebox{0.8}[0.8]{$=$}};
				\node[below, xshift=0.05cm, yshift=-0.55cm] at(y22j) {\scalebox{0.8}[0.8]{$y_{2,jM_1-(t-1)}$}};
				
				\filldraw  (yiij)node[below,xshift=0cm] {\small $y_{i,i_j}$} circle (0.04);
				\filldraw  (yiij+1)node[below,xshift=0cm] {\small $y_{i,i_{j}+1}$} circle (0.04);
				
				\filldraw  (y22j+1)node[below,xshift=0cm] {\small $y_{2,2_{j}+1}$} circle (0.04);
				\node[below,rotate=90, xshift=-0.65cm, yshift=0.15cm] at(y22j+1) {\scalebox{0.8}[0.8]{$=$}};
				\node[below, xshift=0.05cm, yshift=-0.6cm] at(y22j+1) {\scalebox{0.8}[0.8]{$y_{2,jM_1-(t-2)}$}};
				
				\filldraw  (y11j+1)node[below,xshift=0.2cm] {\small $y_{1,1_{j}+1}$} circle (0.04);
				\filldraw  (y)node[right,xshift=0cm] {\small $y$} circle (0.04);

				\draw  (x) to [out=85,in=180]
				coordinate[pos=0.3] (x11j)
				(-0.5,6) to [out=0,in=105]
				coordinate[pos=0.13] (x11j+1)
				coordinate[pos=0.4] (x00)
				(y);
				
				\filldraw  (x11j)node[left,xshift=0cm] {\small $x_{1,1_{j}}$} circle (0.04);
				\filldraw  (x11j+1)node[below,xshift=-0.3cm] {\small $x_{1,1_{j}+1}$} circle (0.04);
				\filldraw  (x00)node[right,xshift=0cm] {\small $x_{0,0}$} circle (0.04);

				\draw  (x) to [out=110,in=180]
				(-2,7) to [out=0,in=100]
				coordinate[pos=0.5] (z11j+1)
				(y11j+1);
				
				\filldraw  (z11j+1)node[above,rotate=-50,xshift=0.1cm,yshift=-0.05cm] {\scalebox{0.8}[0.8]{$z_{1,1_{j}+1}$}}  circle (0.04);
				
				\draw  (x) to [out=60,in=180]
				coordinate[pos=1] (z11j)
				(-5,1) to [out=0,in=110]
				(y11j);
				
				\filldraw  (z11j)node[above,rotate=0,xshift=0cm,yshift=-0.05cm] {\scalebox{0.8}[0.8]{$z_{1,1_{j}}$}}  circle (0.04);

				\draw  (y11j) to [out=60,in=180]
				coordinate[pos=0.3] (x11j)
				(-3.5,0.6) to [out=0,in=110]
				coordinate[pos=0.13] (x11j+1)
				coordinate[pos=0.35] (x00)
				(y22j);

				\draw  (y11j) to [out=80,in=180]
				coordinate[pos=0.9] (x2jm1-t-2)
				(-1,4) to [out=0,in=80]
				(y22j+1);
				
				\filldraw  (x2jm1-t-2)node[above,rotate=10,xshift=-0.1cm,yshift=-0.05cm] {\scalebox{0.7}[0.7]{$x_{2,jM_1-(t-2)}$}}  circle (0.04);
				
				\draw[red]  (y11j) to [out=90,in=180]
				coordinate[pos=0.25] (x2j-1m1+1)
				coordinate[pos=0.46] (x2j-1m-t-1)
				coordinate[pos=0.7] (x2jm1)
				coordinate[pos=0.9] (w11j)
				(-0.5,5.3) to [out=0,in=110]
				(y11j+1);
				
				\filldraw  (x2j-1m1+1)node[above,rotate=75,xshift=0.1cm,yshift=-0.05cm] {\scalebox{0.7}[0.7]{$x_{2,(j-1)M_1+1}$}}  circle (0.04);
				\filldraw  (x2j-1m-t-1)node[above,rotate=65,xshift=0.1cm,yshift=-0.05cm] {\scalebox{0.7}[0.7]{$x_{2,jM_1-(t-1)}$}}  circle (0.04);
				\filldraw  (x2jm1)node[left,rotate=0,xshift=0.1cm,yshift=0.05cm] {\scalebox{0.7}[0.7]{$x_{2,jM_1}$}}  circle (0.04);
				\filldraw  (w11j)node[above,rotate=0,xshift=-0.1cm,yshift=-0.05cm] {\scalebox{0.7}[0.7]{$w_{1,1_j}$}}  circle (0.04);
				
				\draw[red]  (y22j) to [out=80,in=180]
				coordinate[pos=0.8] (x2jm1-t-1)
				(0,3) to [out=0,in=100]
				(y22j+1);
				
				\filldraw[red]  (x2jm1-t-1)node[above,rotate=20,xshift=-0.1cm,yshift=-0.05cm] {\scalebox{0.7}[0.7]{$x_{2,jM_1-(t-1)}$}}  circle (0.04);
				
				\draw  (yiij) to [out=80,in=180]
				coordinate[pos=0.9] (wiij)
				(-0.2,1.5) to [out=0,in=90]
				(yiij+1);
				
				\filldraw  (wiij)node[above,rotate=0,xshift=0cm,yshift=-0.05cm] {\scalebox{0.7}[0.7]{$w_{i,i_j}$}}  circle (0.04);

			\end{tikzpicture}
		\end{center}
		\caption{Illustration for the proof of Lemma \ref{cl25-047}} \label{fig-25-8-8}
	\end{figure}
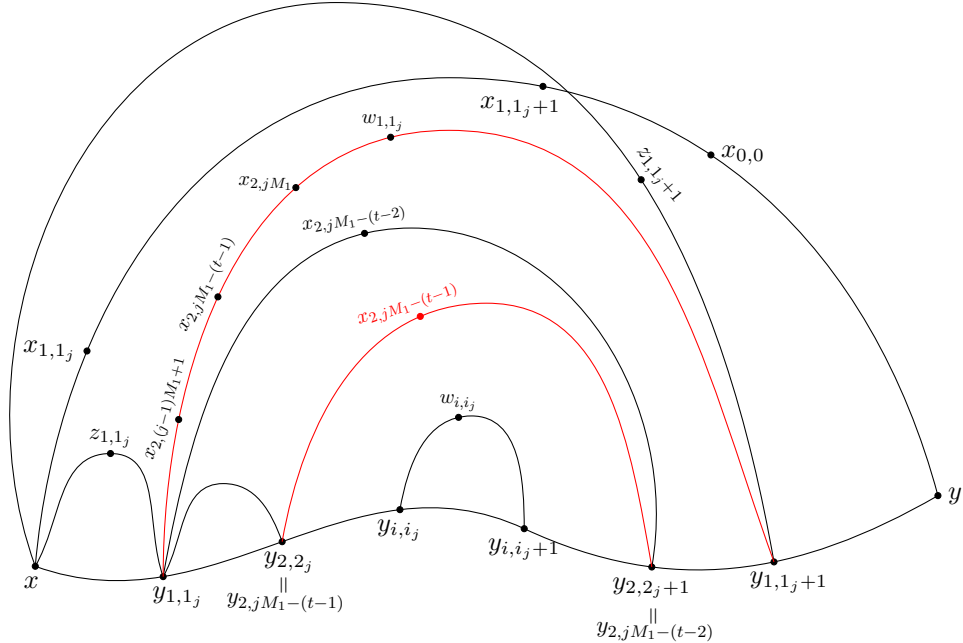
	\blem\label{cl25-047}Given $r_i$ and $N_i$ as above,  there exists some integer  $N_{i-1,1}>\frac{N_{i-1}}{2}$ such that the following conclusions hold:
	\ben
	
	\item\label{cl25-047-3}
	For each $i\in\{3,\cdots, M_0\}$, there exists $\{(i-1)_j\}_{j=1}^{N_{i-1,1}}\subset \bigcup\limits_{t=2}^{M_1}\bigcup\limits_{p=1}^{N_{i-2,1}}\{pM_1-t+1\}$ $($a sequence of integers$)$ with  $(i-1)_j<(i-1)_{j+1}$ for each $j\in \{1,\cdots,N_{i-1,1}-1\}$.

	\item\label{cl25-047-1}
	For each $i\in\{3,\cdots, M_0\}$,  there exists a sequence $\{x_{i,jM_1-(t-1)}\}$, indexed with $j\in\{1,\cdots, N_{i-1,1}\}$ and $t\in\{1,\cdots, M_1\}$, of successive points on $\gamma_{y_{i-1, (i-1)_j}y_{i-1, (i-1)_j+1}}\in \Lambda_{y_{i-1, (i-1)_j}y_{i-1, (i-1)_j+1}}(\Omega)$ so that
	$$d_{\Omega}(x_{i,jM_1-t+1})\leq  r_{i-1}.$$ Moreover, it holds
	$$\frac{3}{4}e^{32C(N_1+2)+2(i-2)}d_{\Omega}(x_{0,0})<\sigma(x_{0,0}, x_{i,jM_1-t+1})<\frac{3}{4}e^{32C(N_1+2)+2(i-1)}d_{\Omega}(x_{0,0}).$$

	\item\label{cl25-047-2}
	For each $i\in\{3,\cdots, M_0\}$, $j\in\{1,\cdots, N_{i-1,1}\}$ and $t\in\{2,\cdots, M_1\}$,
	$$k_{\Omega}(x_{i,jM_1-t+1}, x_{i,jM_1-t+2})>30C.$$
	Moreover, there exist $y_{i,jM_1-t+1}\in \alpha[y_{i-1, (i-1)_j}, y_{i-1, (i-1)_j+1}]$, $\gamma_{y_{i-1, (i-1)_j}y_{i,jM_1-t+1}}$ from $\Lambda_{y_{i-1, (i-1)_j}y_{i,jM_1-t+1}}(\Omega)$ and $z_{i,jM_1-t+1}\in \gamma_{y_{i-1, (i-1)_j}y_{i,jM_1-t+1}}$ such that $$\gamma_{y_{i-1, (i-1)_j}z_{i,jM_1-t+1}y_{i,jM_1-t+1}}\in P_{\alpha[y_{i-1, (i-1)_j}, y_{i-1, (i-1)_j+1}]}^{\gamma_{y_{i-1, (i-1)_j}x_{i,jM_1-t+1}y_{i-1, (i-1)_j+1}}}(3C),$$ and for each $\gamma_{y_{i,jM_1-t+1}y_{i,jM_1-t+2}}\in \Lambda_{y_{i,jM_1-t+1}y_{i,jM_1-t+2}}(\Omega)$, there exists  $w_{i,jM_1-t+1}\in \gamma_{y_{i,jM_1-t+1}y_{i,jM_1-t+2}}$
	such that $$\gamma_{y_{i,jM_1-t+1}w_{i,jM_1-t+1}y_{i,jM_1-t+2}}\in O_{\alpha[y_{i,jM_1-t+1},y_{i,jM_1-t+2}]}^{\gamma_{y_{i-1, (i-1)_j}x_{i,jM_1-t+1}y_{i-1, (i-1)_j+1}}}(2C).$$
	\een
	\elem
	
	\bpf
	We shall prove the following iteration claim:
	\bcl\label{Cl25-5.12} If Lemma \ref{cl25-047} holds for all $i\leq k$, then it holds when $i=k+1$.
	\ecl
	
	By Lemma \ref{Lem4.12-1}, we see that Lemma \ref{cl25-047} holds when $k=2$. To prove the claim, we may assume that Lemma \ref{cl25-047} holds for all $i\leq k$, where $k\geq 2$.
	Since for each $s\in\{1,\cdots,N_{k-1,1}\}$ and $t\in\{1,\cdots,M_1\}$, Lemma \ref{cl25-047}(\ref{cl25-047-2}) yields
	$$\log\Big(1+\frac{\sigma(z_{k,sM_1-t+1},x_{k,sM_1-t+1})}{d_{\Omega}(x_{k,sM_1-t+1})}\Big)\stackrel{\eqref{(2.2)}}{\leq} k_{\Omega}(z_{k,sM_1-t+1},x_{k,sM_1-t+1})\leq 3C,$$
	and so, $$\sigma(z_{k,sM_1-t+1},x_{k,sM_1-t+1})\leq (e^{3C}-1)d_{\Omega}(x_{k,sM_1-t+1})\leq (e^{3C}-1)r_{k-1},$$
	we infer from Lemma \ref{cl25-047}(\ref{cl25-047-1}) that
	$$\sigma(x_{0,0},z_{k,sM_1-t+1})< \sigma(x_{0,0},x_{k,sM_1-t+1})+\sigma(z_{k,sM_1-t+1},x_{k,sM_1-t+1})<e^{32C(N_1+2)+2(k-1)}d_{\Omega}(x_{0,0}).$$
	Thus we may take
	\beqq\label{H25-0510-0}
	u_{k,s}^t\in \gamma_{y_{k-1, (k-1)_s}z_{k,sM_1-t+1}y_{k,sM_1-t+1}}[z_{k,sM_1-t+1},y_{k,sM_1-t+1}]\cap \mathbb{S}_\sigma(x_{0,0}, e^{32C(N_1+2)+2k-2}d_{\Omega}(x_{0,0})).\eeqq

Based on the positions of points on $\gamma_{y_{k-1, (k-1)_s}z_{k,sM_1-t+1}y_{k,sM_1-t+1}}[u_{k,s}^t,y_{k,sM_1-t+1}]$, along the direction from $u_{k,s}^t$ to $y_{k,sM_1-t+1}$,
we arrange the points in the intersection  $\gamma_{y_{k-1, (k-1)_s}z_{k,sM_1-t+1}y_{k,sM_1-t+1}}$ $[u_{k,s}^t,y_{k,sM_1-t+1}]\cap \mathbb{S}_\sigma(x_{0,0}, e^{32C(N_1+2)+2k-2}d_{\Omega}(x_{0,0}))$. 
Let $v_{k,s}^t\not= u_{k,s}^t$ be the first point.
	Now, we define $u_{k,sM_1-t+1}=u_{k,s}^t$ and $v_{k,sM_1-t+1}=v_{k,s}^t$.
	
	As in the proof of Lemma \ref{Lem4.12-1}, we next assert that there exists an integer
	\be\label{H25-0510-1}
	N_{k,1}>\frac{N_{k}}{2}
	\ee
	such that
	\begin{itemize}
		\item there is a sequence of integers $\{k_{1,j}\}_{j=1}^{N_{k,1}}\subset \{1,\cdots,N_k\}$ with $k_{1,j}<k_{1,j+1}$ for each $j \in \{1, \dots, N_{k,1}-1\}$.
		
		\item  for each $j \in \{1, \dots, N_{k,1}\}$, there exists some $s_j\in\{1,\cdots,N_{k-1,1}\}$ so that for each $u \in \gamma_{y_{k-1,(k-1)_{s_j}}z_{k,k_j}y_{k,k_j}} [u_{k,k_j}, y_{k,k_j}]$,
		\be\label{H25-049-1}
		d_{\Omega}(u)\leq r_k=e^{-\frac{1}{2}\log_Q \frac{N_k}{4}}\cdot e^{32C(N_1+2)+2k}d_{\Omega}(x_{0,0}).
		\ee
	\end{itemize}
	
	Indeed, suppose, on the contrary, that our assertion fails. Then there exist $N_{k,2}$ integers  $k_{1,1},\dots,k_{1,N_{k,2}}$ in $\{1,\dots, N_{k}\}$, with $N_{k,2}\geq \frac{N_{k}}{3}$, such that
	\begin{itemize}
		\item for each $s\in\{1, \dots, N_{k,2}-1\}$, $k_{1,s}<k_{1,s+1}$.
		
		\item for each $s\in \{1, \dots, N_{k,2}\}$, there exist an integer $t_s\in\{1,\cdots,N_{k-1,1}\}$ and some point $u_{k,k_{1,s}}^0\in \gamma_{y_{k-1,(k-1)_{t_s}}z_{k,k_{1,s}}y_{k,k_{1,s}}} [u_{k,k_{1,s}}, y_{k,k_{1,s}}]$ such that
		\be\label{H25-0510-2}
		d_{\Omega}(u_{k,k_{1,s}}^0) > r_{k}.
		\ee
	\end{itemize}
	
	Let $B_{0,k}=\mathbb{B}\big(x_{0,0},e^{32C(N_1+2)+2k}d_{\Omega}(x_{0,0})\big)$. For each $s \in \{1, \dots, N_{k,2}\}$, we take
	$$B_{k_{1,s}}=\mathbb{B}\left(u_{k,k_{1,s}}^0, \frac{1}{3}r_k\right).$$
	Then we know from the choice of $u_{k,sM_1-t+1}(=u_{k,s}^t)$ and $v_{k,sM_1-t+1}(=v_{k,s}^t)$ that
	for each $u \in \overline{B_{k_{1,s}}}$,
	\begin{align*}
		d(x_{0,0}, u) &\leq d(x_{0,0}, u_{k,k_{1,s}}^0)+d(u,u_{k,k_{1,s}}^0)\leq e^{32C(N_1+2)+2k-1}d_{\Omega}(x_{0,0})+\frac{1}{3}r_k\\
		&< e^{32C(N_1+2)+2k}d_{\Omega}(x_{0,0}),
	\end{align*}
	and so, $$\overline{B_{k_{1,s}}}\subset B_{0,k}.$$
	
	If all these balls are disjoint, then applying Lemma \ref{qs-5} with $R_k=e^{20C(N_1+2)+2k}d_{\Omega}(x_{0,0})$ and $r=\frac{1}{3}e^{-\frac{1}{2}\log_Q \frac{N_k}{4}}R_k$ gives
	\[
	N_{k,2}<\frac{N_{k}}{3},
	\]
	which is a contradiction.
	
	In the other case, there exist two integers $s_1<s_2\in \{1, \dots, N_{k,2}\}$ such that $B_{k_{1,s_1}}\cap B_{k_{1,s_2}}\not=\emptyset$. It follows that
	$$d(u_{k,k_{1,s_1}}^0,u_{k,k_{1,s_2}}^0)\leq \frac{2}{3}r_k\stackrel{\eqref{H25-0510-2}}{<}\frac{2}{3}\min\{d_{\Omega}(u_{k,k_{1,s_1}}^0),d_{\Omega}(u_{1,k_{1,s_2}}^0)\},$$
	and thus, by Lemma \ref{Lemma-2.1}, we have
	\be\label{H25-0513-1}
	k_{\Omega}(u_{k,k_{1,s_1}}^0,u_{k,k_{1,s_2}}^0)<\frac{20}{9}.
	\ee
	
	For each $i\in\{1,\cdots,k-1\}$, $j\in\{1,\cdots, N_{i,1}\}$ and $t\in\{2,\cdots,M_1\}$, note that $x_{i+1,jM_1}\in\gamma_{y_{i,i_j}y_{i,i_j+1}}[x_{i+1,jM_1}, y_{i,i_j+1}]$,
	and thus, it follows from our induction assumption that
	\be\label{H25-0513-3}
	k_{\Omega}(w_{i,i_j}, x_{i+1,jM_1-t+1})>k_{\Omega}(x_{i+1,jM_1}, x_{i+1,jM_1-t+1})>30C.
	\ee
	Moreover, by {Lemmas \ref{lem-23-4.2}} and  \ref{cl25-047}(\ref{cl25-047-2}),
	there exist  $q\in \{1,\cdots,k-1\}$, $j\in\{1,\cdots, N_{q,1}\}$ and $t_1<t_2\in\{2,\cdots,M_1\}$ such that
	$$\gamma_1:=\gamma_{y_{q,q_j}v_{q+1,jM_1-t_1+1}y_{k-1,(k-1)_{1,s_1}}}\in Q_{\alpha[y_{q,q_j}y_{q,q_j+1}]}^{\gamma_{y_{q,q_j}x_{q+1,jM_1-t_1+1}y_{q,q_j+1}}},$$
	$$\gamma_2:=\gamma_{y_{q,q_j}v_{q+1,jM_1-t_1+1}^1y_{k,k_{1,s_1}}}\in Q_{\alpha[y_{q,q_j}y_{q,q_j+1}]}^{\gamma_{y_{q,q_j}x_{q+1,jM_1-t_1+1}y_{q,q_j+1}}},$$
	$$\gamma_3:=\gamma_{y_{q,q_j}v_{q+1,jM_1-t_2+1}y_{k-1,(k-1)_{1,s_2}}}\in Q_{\alpha[y_{q,q_j}y_{q,q_j+1}]}^{\gamma_{y_{q,q_j}x_{q+1,jM_1-t_1+1}y_{q,q_j+1}}}$$
	and
	$$\gamma_4:=\gamma_{y_{q,q_j}v_{q+1,jM_1-t_2+1}^1y_{k,k_{1,s_2}}}\in Q_{\alpha[y_{q,q_j}y_{q,q_j+1}]}^{\gamma_{y_{q,q_j}x_{q+1,jM_1-t_2+1}y_{q,q_j+1}}}.$$
	Then combining Lemma \ref{lem-24-5.3} with \eqref{H25-0513-3}, we infer that there exist two points
	$$v_{k,k_{1,s_1}}^0\in \gamma_1[v_{q+1,jM_1-t_1+1},y_{k-1,(k-1)_{1,s_1}}]\bigcup \gamma_2[v_{q+1,jM_1-t_1+1}^1, y_{k,k_{1,s_1}}]$$
	and
	$$v_{k,k_{1,s_2}}^0\in\gamma_{3}[v_{q+1,jM_1-t_2+1},y_{k-1,(k-1)_{1,s_2}}]\bigcup \gamma_{4}[v_{q+1,jM_1-t_2+1}^1,y_{k,k_{1,s_2}}]$$
	such that
	$$k_{\Omega}(u_{k,k_{1,s_1}}^0,v_{k,k_{1,s_1}}^0)\leq C\;\quad\mbox{ and }\quad\;k_{\Omega}(u_{k,k_{1,s_2}}^0,v_{k,k_{1,s_2}}^0)\leq C.$$
	This, together with \eqref{H25-0513-1}, shows that
	$$k_{\Omega}(v_{k,k_{1,s_1}}^0,v_{k,k_{1,s_2}}^0)\leq k_{\Omega}(u_{k,k_{1,s_1}}^0,v_{k,k_{1,s_1}}^0)+k_{\Omega}(u_{k,k_{1,s_2}}^0,v_{k,k_{1,s_2}}^0)+k_{\Omega}(u_{k,k_{1,s_1}}^0,u_{k,k_{1,s_2}}^0)<3C,$$
	which, together with Lemma \ref{cl25-047}(\ref{cl25-047-2}), clearly contradicts with Lemma \ref{lem-22-3.1}(\ref{Lemma4-2}). The proof of the assertion is thus complete.

	Let us continue the proof based on the assertion. Let $r_{k}=e^{-\frac{1}{2}\log_Q \frac{N_k}{4}} R_k$. For each $j\in\{1,\dots,N_{k,1}\}$, let $x_{k,j}^1=u_{k,k_j}$, and for each $t\in\{2,\cdots,M_1\}$, let $x_{k,j}^t\in\gamma_{y_{k-1,(k-1)_{s_j}}z_{k,k_{j}}y_{k,k_{j}}}$ be such that
	\be\label{H25-0510-4}
	\sigma(x_{k,k_j}^t,x_{k,k_j}^{t-1})=e^{32C}r_{k}.
	\ee
	Then we have
	\beqq
	k_{\Omega}(x_{k,k_j}^t,x_{k,k_j}^{t-1})\stackrel{\eqref{(2.1)}}{\geq} \log\Big(1+\frac{\sigma(x_{k,k_j}^t,x_{k,k_j}^{t-1})}{d_{\Omega}(x_{k,k_j}^{t-1})}\Big)\stackrel{\eqref{H25-049-1}}{\geq}
	\log\Big(1+\frac{\sigma(x_{k,k_j}^t,x_{k,k_j}^{t-1})}{r_{k}}\Big)\stackrel{\eqref{H25-0510-4}}{>}32C.
	\eeqq
	
	Furthermore, by Lemmas \ref{lem-24-5.3}(\ref{lem-24-5.3-1}) and  \ref{lem-23-4.1}, there exists some point $x_{k,k_j}^{1,1}\in \gamma_{y_{k,k_j}w_{k,k_{j}}y_{k,k_j+1}}[y_{k,k_j},$ $w_{k,k_j}]$
	such that \be\label{H25-0510-6}k_{\Omega}(x_{k,k_j}^1,x_{k,k_j}^{1,1})\leq C.\ee
	
	Based on \eqref{H25-0510-6}, we may argue similarly as in \eqref{H25-04-06-5} to conclude that for each $t\in\{2,\cdots,M_1\}$, there exists some point $x_{k,k_j}^{1,t}\in  \gamma_{y_{k,k_j}w_{k,k_j}y_{k,k_j+1}}[y_{k,k_j},$  $x_{k,k_j}^{1,t-1}]$ satisfying
	\be\label{H25-0510-7}
	k_{\Omega}(x_{k,k_j}^t, x_{k,k_j}^{1,t})\leq C.
	\ee
	
	Finally, for each $j\in\{1,\dots,N_{k,1}\}$ and $t\in\{1,\cdots,M_1\}$, if we take $x_{k+1,jM_1-(t-1)}=x_{k,j}^{1,t}$, then a similar discussion as in the proof of \eqref{H25-0512-1} $\sim$ \eqref{H25-0508-2} in Lemma \ref{Lem4.12-1}(\ref{Lem4.12-1-3}) (using \eqref{H25-0510-4} $\sim$ \eqref{H25-0510-7}), together with the assertion (i.e., \eqref{H25-0510-1}), shows that the claim holds when $i=k+1$. This implies that the claim is true, and thus, the lemma is proved.
	\epf
	
	\bpf[Proof of Theorem \ref{thm-24-5.1}]
	Towards a contradiction, suppose that \eqref{eq-24-1.1} holds, and then, we shall prove the following claim.
	
	\bcl\label{H25-084}
	Let $N=[\frac{\ell_k(\alpha_{xy})}{3C+1}]+1$. Then for each positive integer $\varsigma\in\{1,\cdots, N\}$, there are $y_{0,\varsigma}^1\in\alpha_{x y}$ and  $y_{0,\varsigma}^2\in\alpha_{x y}[y_{0,\varsigma}^1,y]$ which satisfy the following:
	\ben
	\item\label{H5-08-04-1}
	For each $\varsigma\in\{2,\cdots, N\}$, $y_{0,\varsigma}^1\in\alpha_{x y}[y_{0,\varsigma-1}^1,y]$ and $k_{\Omega}(y_{0,\varsigma-1}^1,y_{0,\varsigma}^1)\geq1+3C$.
	
	\item \label{H5-08-04-2}
	For each $\varsigma\in\{1,\cdots, N\}$ and every $\gamma_{y_{0,\varsigma}^1y_{0,\varsigma}^2}\in \Lambda_{y_{0,\varsigma}^1y_{0,\varsigma}^2}(\Omega)$, there exists $x_{0,\varsigma}\in\gamma_{y_{0,\varsigma}^1y_{0,\varsigma}^2}$ such that
	$$\sigma(x_{0,\varsigma},\alpha[y_{0,\varsigma}^1, y_{0,\varsigma}^2])>\tau d_{\Omega}(x_{0,\varsigma}).$$
	
	\item \label{H5-08-04-3}
	For each $\varsigma\in\{1,\cdots, N-1\}$, $k_{\Omega}(y_{0,\varsigma}^1,y_{0,\varsigma}^2)\geq 1+3C$ and $k_{\Omega}(y_{0,N}^1,y_{0,N}^2)<1+3C$.
	\een
	\ecl
	
	Suppose that Claim \ref{H25-084} holds. Then by Claim \ref{H25-084}(\ref{H5-08-04-2}), we have
	$$
	\begin{aligned}
		k_{\Omega}(y_{0,N}^1,y_{0,N}^2)&\geq k_{\Omega}(y_{0,N}^1,x_{0,N}) \stackrel{\eqref{(2.1)}}{\geq} \log\left(1+\frac{\ell(\gamma_{y_{0,N}^1y_{0,N}^2}[y_{0,N}^1,x_{0,N}])}{\min\{d_{\Omega}(y_{0,N}^1), d_{\Omega}(x_{0,N})\}}\right)\\
		&\geq \log\left(1+\frac{\sigma(y_{0,N}^1,x_{0,N})}{d_{\Omega}(x_{0,N})}\right)
		\geq \log \frac{\sigma(x_{0,N},\alpha[y_{0,N}^1,y_{0,N}^2])}{d_{\Omega}(x_{0,N})}>\log\tau,
	\end{aligned}
	$$
	which clearly contradicts with Claim \ref{H25-084}\eqref{H5-08-04-3}.
	Thus to prove the theorem, it suffices to show the claim.
	
	For the case $\varsigma=1$, we infer from Lemma \ref{cl25-047}(\ref{cl25-047-1})  (with $i=M_0$) that for each $j\in\{1\cdots, N_{M_0-1,1}\}$ and $t\in\{1,\cdots,M\}$, there exists some point $$x_{M_0,jM_1-t+1}\in\gamma_{y_{M_0-1,(M_0-1)_j}y_{M_0-1,(M_0-1)_j+1}}$$
	with $y_{M_0-1,(M_0-1)_j}$ and $y_{M_0-1,(M_0-1)_j+1}$ given by Lemma \ref{cl25-047}(\ref{cl25-047-1}), so that
	$$\sigma(x_{0,0},x_{M_0,jM_1-t+1})< \frac{5}{4}e^{32C(N_1+2)+2M_0+1}d_{\Omega}(x_{0,0})<\tau^{\frac{1}{C}}d_{\Omega}(x_{0,0}).$$
	Then for each $z\in\alpha[y_{M_0-1,(M_0-1)_j},y_{M_0-1,(M_0-1)_j+1}]$, it follows from the triangle inequality that
	\begin{align*}
		\sigma(x_{M_0,jM_1-t+1},z) &\geq\sigma(x_{0,0},z)- \sigma(x_{0,0},x_{M_0,jM_1-t+1})\stackrel{\eqref{eq-24-1.1}}{\geq} (\tau -\tau^{\frac{1}{C}})d_{\Omega}(x_{0,0}).
	\end{align*}
	
	Moreover, by Lemma \ref{cl25-047}(\ref{cl25-047-1}), we have
	\[
	\begin{aligned}
		d_{\Omega}(x_{M_0,jM_1-t+1})\leq r_{M_0-1}=e^{-\frac{1}{2}\log_Q \frac{N_{M_0-1}}{4}}\cdot e^{32C(N_1+2)+2M_0-2}d_{\Omega}(x_{0,0})
	\end{aligned}
	\]
	and
	$$N_{M_0-1}\geq \frac{(M_1-1)^{M_0-2}}{2^{M_0-2}}\cdot N_1>e^{128QC^3 N_1}.$$
	Thus we deduce that
	$$d_{\Omega}(x_{M_0,jM_1-t+1})\leq \tau^{-1}d_{\Omega}(x_{0,0}).$$
	Combining the above two estimates gives
	\[
	\sigma(x_{M_0,jM_1-t+1},\alpha[y_{M_0-1,(M_0-1)_j},y_{M_0-1,(M_0-1)_j+1}])>\tau d_{\Omega}(x_{M_0,jM_1-t+1}).
	\]
	
	Meanwhile, {by Lemma \ref{Lem4.12-1}\eqref{Lem4.12-1-2}, we know that
		$y_{M_0-1,(M_0-1)_j}\in\alpha[y_{1,1},y]$. This, together with
		\eqref{eq-24-4}, shows that there exists $\zeta_{1,1}\in \gamma_{xy_{M_0-1,(M_0-1)_j}}$ such that
		$$k_{\Omega}(\zeta_{1,1}, x_{1,1})<3C.$$
		Hence we get from \eqref{H25-0507-1} that
		$$k_{\Omega}(x,y_{M_0-1,(M_0-1)_j})>k_{\Omega}(x,\zeta_{1,1})\geq k_{\Omega}(x,x_{1,1})-k_{\Omega}(x_{1,1},\zeta_{1,1})>\frac{1}{2}\log \tau-3C>1+3C.$$}
	
	Set $y_{0,1}^1=y_{M_0-1,(M_0-1)_j}\in\alpha_{x y}$,  $y_{0,1}^2=y_{M_0-1,(M_0-1)_j+1}\in\alpha_{x y}[y_{0,1}^1,y]$ and $x_{0,1}=x_{M_0,jM_1-t+1}\in\gamma_{y_{0,1}^1y_{0,1}^2}$. Then Claim \ref{H25-084} holds for the case $\varsigma=1$.
	
	Suppose that Claim \ref{H25-084} holds when $\varsigma=k$ for $k\in\{1,\cdots,N-1\}$. By replacing $x=y_{0,k}^1$, $y=y_{0,k}^2$ and $x_{0,k}=x_{0,0}$, we see that Claim \ref{H25-084} holds when $\varsigma=k+1$.
	This completes the proof of Claim \ref{H25-084}(\ref{H5-08-04-1}) and (\ref{H5-08-04-2}).
	
	Note that  Claim \ref{H25-084}\eqref{H5-08-04-3} is a direct consequence of  Claim \ref{H25-084}(\ref{H5-08-04-1}) and (\ref{H5-08-04-2}). Indeed, by Claim \ref{H25-084}(\ref{H5-08-04-1}), we have
	\beqq
	\begin{aligned}
		(N-1)(1+3C)&+k_{\Omega}(y_{0,N}^1,y_{0,N}^2)\leq \sum\limits_{\varsigma=2}^N k_{\Omega}(y_{0,\varsigma-1}^1,y_{0,\varsigma}^1)+k_{\Omega}(y_{0,N}^1,y_{0,N}^2)\\ &\leq \sum\limits_{\varsigma=2}^N \ell_k(\alpha_{xy}[y_{0,\varsigma-1}^1,y_{0,\varsigma}^1])+k_{\Omega}(y_{0,N}^1,y_{0,N}^2)\leq \ell_k(\alpha_{xy})\leq N(1+3C),
	\end{aligned}
	\eeqq
	which implies $k_{\Omega}(y_{0,N}^1,y_{0,N}^2)<1+3C$, and thus, the proof of theorem is complete.
	\epf

	\section{Ball separation condition implies Gehring-Hayman inequality}\label{sec-4-2}
In this section, we shall prove Theorem \ref{positive-answer}, namely, the ball separation condition implies the Gehring-Hayman inequality. Throughout this section, $X=(X, d)$ is assumed to be a $Q$-doubling length space, and $\Omega$ denotes a proper subdomain in $X$ so that $(\Omega,k)$ is geodesic. We also assume that $\Omega$ satisfies the $\theta$-ball separation condition with $\theta>1$ and shall show that $\Omega$ satisfies the $\theta_1$-Gehring-Hayman inequality with $\theta_1=(2Q)^{4(36\theta\cdot Q^{5\theta})^{(8Q\theta)^8}}$.

\subsection{A version of diameter Gehring-Hayman inequality}

Let $G_{xy}(\Omega)$ denote the collection of all rectifiable curves $\alpha\in\Gamma_{xy}(\Omega)$ such that
$$\ell(\alpha)\leq \sigma(x,y)+e^{-\theta_1}\min\{\sigma(x,y), d_{\Omega}(x),d_{\Omega}(y)\}.$$
For convenience, we shall use the notation $\alpha_{xy}$ to represent a curve in $G_{xy}(\Omega)$.

For the proof of Theorem \ref{positive-answer}, we first prove the following weaker version of Gehring-Hayman inequality.
\begin{thm}\label{2017-10-11-1}
	{For any $x_1,x_2\in \Omega$, let $\gamma_{x_1x_2}\in \Lambda_{x_1x_2}(\Omega)$ and $\alpha_{x_1x_1}\in G_{x_1x_2}(\Omega)$.
		Then $$\diam_{\sigma}(\gamma_{x_1x_2})\leq \theta_0\ell(\alpha_{x_1x_2}),$$
		where $\theta_0=e^{(36\theta\cdot Q^{5\theta})^{(8Q\theta)^8}}$.}
\end{thm}

For the proof of Theorem \ref{2017-10-11-1}, we need several auxiliary lemmas.

\blem\label{H25-09-20} Fix $x_1, x_2\in\Omega$, $\gamma_{x_1x_2}\in \Lambda_{x_1x_2}(\Omega)$ and $\alpha_{x_1x_2}\in G_{x_1x_2}(\Omega)$. Then the following holds.
\ben
\item\label{H25-09-20-2}
Fix $x\in\gamma_{x_1x_2}$. If $\sigma(x_1,x)\geq \varsigma\ell(\alpha_{x_1x_2})$ with $\varsigma>1$, then $\sigma(x_1,x)\leq \frac{\varsigma\theta}{\varsigma-1} d_{\Omega}(x)$.

\item\label{H25-09-20-1} Fix $\varsigma>0$ and $y_1\in \Omega$ with $k_{\Omega}(x_1,y_1)\leq \frac{1}{2}$. Suppose that for each $x\in\gamma_{x_1x_2}$, $\sigma(x_1,x)\leq \varsigma d_{\Omega}(x)$. Then for each {$\gamma_{y_1x_2}\in \Lambda_{y_1x_2}(\Omega)$ and each $y\in \gamma_{y_1x_2}$}, it holds
$$\sigma(y_1,y)\leq 5(\theta+\varsigma(1+\theta)) d_{\Omega}(y).$$

\item\label{H25-09-20-3a}
Fix {$\varsigma>1$}, $y_0\in\gamma_{x_1x_2}$, $y_1\in\gamma_{x_1x_2}[x_1,y_0]$ and $y_2\in\gamma_{x_1x_2}[x_2,y_0]$  with $\sigma(x_1,y_2)=\sigma(x_2,y_1)=\varsigma\ell(\alpha_{x_1x_2})$. Suppose that for each  $i\in\{1,2\}$ and each $x\in\gamma_{x_1x_2}[y_i,y_0]$, $\sigma(x_j,x)\geq \varsigma\ell(\alpha_{x_1x_2})$ for each $j\not=i$ and $j\in\{1,2\}$. Then $\sigma(y_1,x)< \frac{3\varsigma\theta}{\varsigma-1} d_{\Omega}(x)$.
\een

\elem

\bpf (1). Fix $x\in\gamma_{x_1x_2}$ with $\sigma(x_1,x)\geq \varsigma\ell(\alpha_{x_1x_2})$.  Then the $\theta$-ball separation condition implies that there exists some point
$y\in \alpha_{x_1x_2}$ such that $$\sigma(x,y)\leq \theta d_{\Omega}(x).$$
Hence by the choice of $x$, we have $$\theta d_{\Omega}(x)\geq\sigma(x,y)\geq \sigma(x_1,x)-\sigma(x_1,y)\geq (1-\frac{1}{\varsigma})\sigma(x_1,x),$$
which implies that the statement (\ref{H25-09-20-2}) of the lemma holds.

(2). Fix $\gamma_{y_1x_2}\in\Lambda_{y_1x_2}(\Omega)$ and $y\in\gamma_{y_1x_2}$. If $\sigma(y_1,y)<\frac{7}{8}d_{\Omega}(y_1)$, then we have
$$d_{\Omega}(y)\geq d_{\Omega}(y_1)-\sigma(y_1,y)>\frac{1}{7}\sigma(y_1,y)>\frac{1}{5(\theta+\varsigma(1+\theta))} \sigma(y_1,y),$$
which gives our desired estimate.

It remains to consider the case $\sigma(y_1,y)\geq \frac{7}{8}d_{\Omega}(y_1)$. In this case,  we note that
$$\log\Big(1+\frac{\sigma(x_1,y_1)}{\min\{d_{\Omega}(x_1),d_{\Omega}(y_1)\}}\Big)\stackrel{\eqref{(2.2)}}{\leq} k_{\Omega}(x_1,y_1)\leq \frac{1}{2},$$
and so, \be\label{H25-09-20-3}\sigma(x_1,y_1)<\frac{7}{10}\min\{d_{\Omega}(x_1),d_{\Omega}(y_1)\}.\ee
This, together with our assumption, gives
\be\label{H25-09-20-0}
\sigma(y_1,y)>\frac{5}{4}\sigma(y_1,x_1).
\ee

Fix $\gamma_{x_1y_1}\in \Lambda_{x_1y_1}(\Omega)$. Then it follows from the $\theta$-ball separation condition that there exists some point
$w\in \gamma_{x_1x_2}\cup \gamma_{x_1y_1}$ such that \be\label{H25-09-20-4}\sigma(y,w)\leq \theta d_{\Omega}(y),\ee
and thus,
\be\label{H25-09-20-5}
d_{\Omega}(w)\leq \sigma(y,w)+d_{\Omega}(y)\leq (1+\theta)d_{\Omega}(y).
\ee

If $w\in \gamma_{x_1y_1}$, then
$$d_{\Omega}(w)\geq d_{\Omega}(x_1)-\sigma(x_1,y_1)\stackrel{\eqref{H25-09-20-3}}{\geq} \frac{3}{7}\sigma(x_1,y_1),$$
from which it follows that
$$\sigma(y_1,y)\leq \sigma(y_1,w)+\sigma(y,w)\stackrel{\eqref{H25-09-20-4}}{\leq} \frac{7}{3}d_{\Omega}(w)+\theta d_{\Omega}(y)\stackrel{\eqref{H25-09-20-5}}{<}(3+4\theta) d_{\Omega}(y).$$

If $w\in \gamma_{x_1x_2}$, then it follows from the assumption and \eqref{H25-09-20-0} $\sim$ \eqref{H25-09-20-5} that
\beq
\begin{aligned}
	\sigma(y_1,y)\nonumber&\leq \sigma(y_1,x_1)+\sigma(x_1,w)+\sigma(w,y)
	\\ &\leq \frac{4}{5}\sigma(y_1,y)+\varsigma d_{\Omega}(w)+\theta d_{\Omega}(y)
	\\ &\leq \frac{4}{5}\sigma(y_1,y)+ (\theta+\varsigma(1+\theta))d_{\Omega}(y),
\end{aligned}
\eeq
which gives $$\sigma(y_1,y)\leq 5(\theta+\varsigma(1+\theta))d_{\Omega}(y).$$

{In either case}, the statement (\ref{H25-09-20-1}) is proved.

(3). For each $i\in \{1,2\}$, by the assumption, it holds
\be\label{H25-12-02-1}\sigma(x_i,x)\geq \varsigma \ell(\alpha_{x_1x_2})= \sigma(x_1,y_1).\ee
If $x\in\gamma_{x_1x_2}[y_1,y_0]$, then we have
$$\sigma(y_1,x)\leq \sigma(x_2,y_1)+\sigma(x_2,x)\stackrel{ (\ref{H25-12-02-1})}{\leq} 2\sigma(x_2,x)\stackrel{\text{Lemma }\ref{H25-09-20} (\ref{H25-09-20-2})}{\leq} \frac{2\varsigma\theta}{\varsigma-1} d_{\Omega}(x).$$	
If $x\in\gamma_{x_1x_2}[y_2,y_0]$, then
\beqq
\begin{aligned}\sigma(y_1,x) &\leq \sigma(x_1,y_1)+\sigma(x_1,x)\leq\sigma(x_2,y_1)+\sigma(x_1,x_2)+\sigma(x_1,x)\\ & \stackrel{ (\ref{H25-12-02-1})}{\leq} (2+\varsigma^{-1})\sigma(x_2,x)\stackrel{\text{Lemma }\ref{H25-09-20} (\ref{H25-09-20-2})}{<} \frac{3\varsigma\theta}{\varsigma-1} d_{\Omega}(x).\end{aligned}
\eeqq																													
These establish the statement \eqref{H25-09-20-3a}, and hence, the proof of Lemma \ref{H25-09-20} is complete.
\epf

\begin{lem}\label{qs-6}
	
	For any pair of points $x, y\in \Omega$ with $$d(x,y)\geq \frac{1}{2}\max\{d_{\Omega}(x),d_{\Omega}(y)\},$$
	and $\beta\in \Lambda_{xy}(\Omega)$, there exists a finite sequence of balls $\{B_i\}_{i=1}^{k_1}$ in $\Omega$ such that
	\ben
	\item\label{qs-6(1)}
	for each $i\in\{1,\ldots,k_1\}$, $B_i=\mathbb{B}(y_i,r_i)$ with $r_i=\frac{1}{4}d_{\Omega}(y_i)$, where $y_0=x$, $y_i\in \beta[y_{i-1},y]$, but $y_i\notin B_{i-1}$.
	\item\label{qs-6(2)}
	$y_{k_1+1}=y\in B_{k_1}$ $($Possibly, $y_{k_1+1}=y_{k_1}$$)$.
	\item\label{qs-6(4)}
	for any pair $\{i,j\} \subset\{1,\ldots,k_1\}$ with $|j-i|>1$,
	$B_i\cap B_j=\emptyset$	
	provided that $k_1\geq 3$.
	\item for each $i\in\{1,\ldots,k_1-1\}$, we have
	\ben
	\item\label{qs-6(3)}
	$B_i\cap B_{i+1}\not=\emptyset$.
	\item\label{qs-6(7)}
	$\ell(\beta[y_i,y_{i+1}])\leq \frac{11}{9}d_{\Omega}(y_i)\leq \frac{44}{9}d(y_i,y_{i+1}).$
	\item \label{qs-6(6)}
	$\frac{9}{20}d_{\Omega}(z)\leq d_{\Omega}(w)\leq \frac{20}{9}d_{\Omega}(z)$
	for all $z,$ $w\in \beta[y_i,y_{i+1}]$.
	\een
	\item\label{qs-6(0-4)}
	for each $i\in\{1,\ldots,k_1-1\}$, $\log \frac{5}{4}\leq k_{\Omega}(y_i,y_{i+1})\leq \frac{20}{27}$.
	\item\label{qs-6(8)}
	$\ell(\beta[y_{k_1},y_{k_1+1}])\leq \frac{9}{20}d_{\Omega}(y_{k_1})$ and $k_{\Omega}(y_{k_1},y_{k_1+1})\leq \frac{10}{27}$.
	\een
\end{lem}

\bpf The proof consists of a few steps.
\medskip

\textbf{Step 1}. Construct an initial sequence $\{x_j\}_{j=1}^{k_0}$ of points on $\beta$.
\medskip

In this step, we shall prove that there exists a finite sequence $\{x_j\}_{j=1}^{k_0}\subset \beta$ with $x_1=x$ such that
$$d(x_j, x_{j+1})=\frac{1}{4}d_{\Omega}(x_j)\;\;\mbox{and}\;\;y\in \mathbb{B}\left(x_{k_0}, \frac{1}{4}d_{\Omega}(x_{k_0})\right).$$

To this end, set $x_1=x$. According to the positions of points on
$\beta$, along the direction from $x$ to $y$, we arrange the points in the intersection $\beta\cap \mathbb{S}(x, \frac{1}{4}d_{\Omega}(x))$. Let $x_2$ be the last point.

If $y\in \mathbb{B}(x_2, \frac{1}{4}d_{\Omega}(x_2))$, then we take $k_0=2$. If  $y\notin \mathbb{B}(x_2, \frac{1}{4}d_{\Omega}(x_2))$,
according to the positions of points on
$\beta[x_2,y]$, along the direction from $x_2$ to $y$, we arrange the points in the intersection $\beta[x_2,y]\cap \mathbb{S}(x_2, \frac{1}{4}d_{\Omega}(x_2))$. Let $x_3$ be the last point.

If $y\in \mathbb{B}(x_3, \frac{1}{4}d_{\Omega}(x_3))$, then we take $k_0=3$. If  $y\notin \mathbb{B}(x_3, \frac{1}{4}d_{\Omega}(x_3))$, $\ldots$.

Repeating this procedure for $\nu$ times, we find a point $x_{\nu}\in \beta[x_{\nu-1},y]$ such that $y\in \mathbb{B}(x_{\nu}, \frac{1}{4}d_{\Omega}(x_{\nu}))$. We take $k_0=\nu$. It is possible that $x_{k_0}=y$.
\medskip

\textbf{Step 2}. Select a subsequence $\{y_i\}_{i=1}^{k_1+1}$ of  $\{x_j\}_{j=1}^{k_0}$ and, accordingly, get a ball sequence $\{B_i\}_{i=1}^{k_1}$, where
$B_i=\mathbb{B}\left(y_i, \frac{1}{4}d_{\Omega}(y_i)\right)$.
\medskip

Let $\{x_j\}_{j=1}^{k_0}$ be the sequence constructed in Step 1. For each $i\in \{1, \ldots, k_0\}$, set
$$B_{1,j}=\mathbb{B}\left(x_j, \frac{1}{4}d_{\Omega}(x_j)\right).$$
We are going to select the desired sequence of balls from $\{B_{1,j}\}_{j=1}^{k_0}$.

\begin{enumerate}
	\item Set $B_1=B_{1,1}$ and $y_1=x_1$. Then we define
	$$r_1=\max\{r:\;r\in \{2,\ldots, k_0\}\;\;\mbox{and}\;\; B_1\cap B_{1,r}\not=\emptyset\}.$$
	
	\item Set $B_2=B_{1,r_1}$ and $y_2=x_{r_1}$.
	\begin{itemize}
		\item 	If $r_1=k_0$, then we find the needed sequence of balls $\{B_i\}_{i=1}^{k_1}$ with $k_1=2$.
		\item  If $r_1< k_0$, then define
		$$r_2=\max\{r:\;r\in \{r_1+1,\ldots, k_0\}\;\;\mbox{and}\;\; B_2\cap B_{1,r}\not=\emptyset\}.$$
	\end{itemize}
	
	\item Set $B_3=B_{1,r_2}$ and $y_3=x_{r_2}$.
	\begin{itemize}
		\item If $r_2=k_0$, then we find the needed sequence of balls $\{B_i\}_{i=1}^{k_1}$ with $k_1=3$.
		\item If $r_2< k_0$, then define
		$$r_3=\max\{r:\;r\in \{r_2+1,\ldots, k_0\}\;\;\mbox{and}\;\; B_3\cap B_{1,r}\not=\emptyset\}.$$
	\end{itemize}
	
	\item Set $B_4=B_{1,r_3}$ and $y_4=x_{r_3}$. $\cdots$.
\end{enumerate}

By repeating this procedure, we find an integer $k_1\leq k_0$ such that
$$\max\{r:\;r\in \{r_{k_1-1}+1,\ldots, k_0\}\;\;\mbox{and}\;\; B_{k_1-1}\cap B_{1,r}\not=\emptyset\}=k_0.$$
Then set $B_{k_1}=B_{1,k_0}$, $y_{k_1}=x_{k_0}$ and $y_{k_1+1}=y$. It is possible that $y_{k_1+1}=y_{k_1}$.

In this way, we find the desired subsequence $\{y_i\}_{i=1}^{k_1+1}$ of  $\{x_j\}_{j=1}^{k_0}$ and the corresponding ball sequence $\{B_i\}_{i=1}^{k_1}$.
\medskip

\textbf{Step 3}. Verify all the listed properties.
\medskip

In this step, we shall prove that the point sequence $\{y_i\}_{i=1}^{k_1+1}$ and the ball sequence $\{B_i\}_{i=1}^{k_1}$
selected in \text{Step 2} satisfy all listed properties. That $\{B_i\}_{j=1}^{k_1}$ satisfies the properties \eqref{qs-6(1)} $ \sim$ \eqref{qs-6(4)} and \eqref{qs-6(3)} of the lemma is clear from the construction.

For each $i\in\{1,\ldots,k_1-1\}$, by the property \eqref{qs-6(4)} of the lemma, we have
$$\mathbb{S}\left(y_i, \frac{1}{4}d_{\Omega}(y_i)\right)\bigcap \mathbb{B}\left(y_{i+1}, \frac{1}{4}d_{\Omega}(y_{i+1})\right)\not=\emptyset.$$
Select $y_{1,i}\in\mathbb{S}\left(y_i, \frac{1}{4}d_{\Omega}(y_i)\right)\bigcap \mathbb{B}\left(y_{i+1}, \frac{1}{4}d_{\Omega}(y_{i+1})\right)$.
Applying Lemma \ref{Lemma-2.1} with $a=4$, we obtain
$$k_{\Omega}(y_i,y_{1,i})\leq \frac{10}{27}\;\;\mbox{ and }\;\;k_{\Omega}(y_{1,i},y_{i+1})\leq \frac{10}{27},$$ which implies
\be\label{2017-10-09-12} k_{\Omega}(y_i,y_{i+1})\leq \frac{20}{27}.\ee

By \eqref{(2.2)}, for any pair of points $z,$ $w\in \beta[y_i, y_{i+1}]$, we have
$$\max\Big\{\log\Big(1+\frac{\ell(\beta[y_i,y_{i+1}])}{d_{\Omega}(y_i)}\Big),\Big|\log\frac{d_{\Omega}(w)}{d_{\Omega}(z)}\Big|,\Big|\log\frac{d_{\Omega}(z)}{d_{\Omega}(w)}\Big|\Big\}\leq k_{\Omega}(y_i,y_{i+1}),$$
and thus, it follows from \eqref{2017-10-09-12} and the above estimate that
\beq\label{20-7-27-1}
\frac{9}{20}d_{\Omega}(z)\leq d_{\Omega}(w)\leq \frac{20}{9}d_{\Omega}(z)\;\;\mbox{and}\;\;\ell(\beta[y_i,y_{i+1}])\leq \frac{11}{9}d_{\Omega}(y_i).
\eeq
Moreover, the fact that $y_{i+1}\notin \mathbb{B}(y_i, \frac{1}{4}d_{\Omega}(y_i))$ for $i\in\{1,\ldots,k_1-1\}$ implies
\beqq
d(y_i,y_{i+1})\geq \frac{1}{4}d_{\Omega}(y_i),
\eeqq
which, together with \eqref{(2.1)}, shows that
\beqq
k_{\Omega}(y_i,y_{i+1})\geq \log\Big(1+\frac{d(y_i,y_{i+1})}{d_{\Omega}(y_i)}\Big)>\log \frac{5}{4}.\eeqq

Since $y_{k_1+1}=y\in B_{k_1}$, Lemma \ref{Lemma-2.1} with $a=4$ gives
\beq\label{20-07-30-1}k_{\Omega}(y_{k_1},y_{k_1+1})\leq\frac{10}{27}.\eeq
A similar argument as in \eqref{20-7-27-1}, together with \eqref{(2.2)} and \eqref{20-07-30-1}, gives
\beq\label{20-07-30-2}
\ell(\beta[y_{k_1},y_{k_1+1}])\leq \frac{9}{20}d_{\Omega}(y_{k_1}).
\eeq

Now, we conclude from \eqref{2017-10-09-12} $\sim$ \eqref{20-07-30-2} that all the remaining properties listed in the lemma hold, and hence, the proof  is complete.
\epf

\begin{lem}\label{qs-8}
	Let $\beta\in \Lambda_{xy}(\Omega)$. Suppose that there is a point $x_0\in\beta$ satisfying {$d_{\Omega}(x_0)\geq \frac{1}{2}\sup\limits_{w\in\beta}\{d_{\Omega}(w)\}$}. If there exists a constant $\mu_1\geq 1$ such that for any $z\in\beta$,
	\beq\label{20-07-30-3}
	d(x,z)\leq \mu_1 d_{\Omega}(z),
	\eeq
	then
	$$\ell(\beta)\leq  \lambda d_{\Omega}(x_0)\;\;\mbox{ and }\;\;
	k_{\Omega}(x,y)\leq \lambda \log \frac{3d_{\Omega}(x_0)}{d_{\Omega}(x)},$$
	where $\lambda={\frac{22}{9}}e^2\mu_1\big([Q^{\log_2 4e^2\mu_1(1+\mu_1)}]+1\big)$.
\end{lem}

\bpf If $d(x,y)<\frac{1}{2}\max\{d_{\Omega}(x), d_{\Omega}(y)\}$, then Lemma \ref{Lemma-2.1} with $a=2$ implies
\beqq
k_{\Omega}(x,y)\leq \frac{10}{9}\leq \lambda \log \frac{3d_{\Omega}(x_0)}{d_{\Omega}(x)}.
\eeqq
This, together with \eqref{(2.2)}, gives
$$\log\Big(1+\frac{\ell(\beta)}{d_{\Omega}(x)}\Big)\leq k_{\Omega}(x,y)\leq \frac{10}{9},$$ from which it follows that
\beqq
\ell(\beta)\leq \big(e^{\frac{10}{9}}-1\big)d_{\Omega}(x)\leq \big(e^{\frac{10}{9}}-1\big)d_{\Omega}(x)< \lambda d_{\Omega}(x_0).
\eeqq

For the remaining case, that is,
\beqq
d(x,y)\geq \frac{1}{2}\max\{d_{\Omega}(x), d_{\Omega}(y)\},
\eeqq
by Lemma \ref{qs-6}, there exist finite sequences of balls $\big\{B_i=\mathbb{B}(y_i,r_i)\big\}_{i=1}^{k_1}$ and points $\{y_i\}_{i=1}^{k_1+1}$ in $\Omega$ which satisfy all the properties listed in Lemma \ref{qs-6}. Let $w_0\in \beta$ be such that
$$
d(x,w_0)\geq \frac{1}{2}\sup_{z\in\beta}\{d(x,z)\}.
$$
Then there is an integer $k_2$, depending on $x$ and $w_0$, such that
\be\label{qs-11}e^{k_2-1}d_{\Omega}(x)\leq d(x,w_0)<e^{k_2}d_{\Omega}(x).\ee

To continue the proof, we consider two cases.
\bca\label{yako-0-0}
$k_2\leq \log(3\mu_1)$.
\eca

Let $$\mathcal{B}=\mathbb{B}(x,4\mu_1d_{\Omega}(x)).$$
Then for each $i\in\{1,\ldots,k_1\}$, we have
$$d(x,y_i)\leq d(x,w_0)\stackrel{\eqref{qs-11}}{\leq} e^{k_2}d_{\Omega}(x)\leq 3\mu_1 d_{\Omega}(x),$$
and then,
$$4r_i=d_{\Omega}(y_i)\leq d_{\Omega}(x)+d(x,y_i)\leq (1+3\mu_1)d_{\Omega}(x).$$
It follows that for each $i\in\{1,\ldots,k_1\}$, $d(y_i,x)+r_i<4\mu_1d_{\Omega}(x)$, and so,
\beq\label{20-07-28}
B_i\subset \mathcal{B}.
\eeq

Let $a=16\mu_1(1+\mu_1)$ and $R=4\mu_1d_{\Omega}(x)$. Since
for each $i\in\{1,\ldots,k_1\}$,
\beq\label{20-07-29-2}
d_{\Omega}(x)\leq d(x, y_i)+d_{\Omega}(y_i)\stackrel{\eqref{20-07-30-3}}{\leq} (1+\mu_1)d_{\Omega}(y_i),
\eeq
we get from the above estimate that
$$\frac{R}{a}=\frac{d_{\Omega}(x)}{4(1+\mu_1)}\leq \frac{1}{4}d_{\Omega}(y_i).$$
Based on this, \eqref{20-07-28} and Lemma \ref{qs-6}\eqref{qs-6(4)}, we may apply Lemma \ref{qs-5} to conclude
\beqq k_1\leq [Q^{\log_2 16\mu_1(1+\mu_1)}]+1.\eeqq
Then it follows from the above estimate and Lemma \ref{qs-6}\eqref{qs-6(0-4)} that
\beqq
k_{\Omega}(x,y)=\sum_{i=1}^{k_1}k_{\Omega}(y_i,y_{i+1})
\leq  \frac{20}{27}\left([Q^{\log_2 16\mu_1(1+\mu_1)}]+1\right),
\eeqq
and so, the assumption $d_{\Omega}(x_0)\geq \frac{1}{2}\sup\limits_{w\in\beta}\{d_{\Omega}(w)\}$ in the lemma ensures that $$k_{\Omega}(x,y)\leq \lambda \log \frac{3d_{\Omega}(x_0)}{d_{\Omega}(x)}.$$
Moreover, since for all $i\in\{1,\ldots,k_1\}$, $d_{\Omega}(y_i)\leq d_{\Omega}(x_0)$, we obtain
\beqq
\begin{aligned}
	\ell(\beta)&= \sum_{i=1}^{k_1}\ell(\beta[y_i,y_{i+1}])\stackrel{\text{Lemma } \ref{qs-6}}{\leq}  \frac{11}{9}\sum_{i=1}^{k_1}d_{\Omega}(y_i)\\
	& \leq {\frac{22}{9}}\left([Q^{\log_2 16\mu_1(1+\mu_1)}]+1\right)d_{\Omega}(x_0).
\end{aligned}
\eeqq

In this case, we have proved the lemma.

\bca\label{yako-0-1}
$k_2>\log(3\mu_1)$.
\eca
For each $p\in\{1,\ldots,k_2+1\}$, set $$\mathcal{B}_p=\mathbb{B}(x, e^p d_{\Omega}(x)).$$ Then \eqref{qs-11} implies that
\beq\label{20-07-29-3}
\beta\subset \mathcal{B}_{k_2}.
\eeq

For each $z\in \mathcal{B}_p\cap \beta$, we observe that
\beq\label{20-07-29-1}
\mathbb{B}\left(z,\frac{1}{4}d_{\Omega}(z)\right)\subset \mathcal{B}_{p+1}.
\eeq
Indeed, for any $w\in \mathbb{B}\big(z,\frac{1}{4}d_{\Omega}(z)\big)$, we have
$$d(x,w)\leq d(x,z)+\frac{1}{4}d_{\Omega}(z)\leq \frac{5}{4}d(x,z)+\frac{1}{4}d_{\Omega}(x)\leq \frac{1}{4}(1+5e^p)d_{\Omega}(x),$$
which gives \eqref{20-07-29-1}.

It follows from \eqref{20-07-29-3} and \eqref{20-07-29-1} that $$\bigcup_{i=1}^{k_1}B_i\subset \mathcal{B}_{k_2+1}.$$

For $p\in\{1,\ldots,k_2\}$, let $$\beta_{p}=\beta\cap (\mathcal{B}_p\backslash \mathcal{B}_{p-1}),$$ where $\mathcal{B}_0=\emptyset,$
and define
\[
t_{p}=
\begin{cases}
	0, &\text{if }\beta_{p}\cap \{y_i\}_{i=1}^{k_1}=\emptyset,\\
	{\rm card} \big\{\beta_{p}\cap \{y_i\}_{i=1}^{k_1}\big\}, &\text{otherwise}.
\end{cases}
\]

When $\beta_{p}\cap \{y_i\}_{i=1}^{k_1}\not=\emptyset$, let
$$\{y_{p,j}\}_{j=1}^{t_{p}}=\beta_{p}\cap \{y_i\}_{i=1}^{k_1},$$
and denote by $B_{p,j}$ the ball in $\{B_i\}_{i=1}^{k_1}$ with center $y_{p,j}$. Then it follows from \eqref{20-07-29-1} that
$$\bigcup_{j=1}^{t_{p}}B_{p,j}\subset \mathcal{B}_{p+1}.$$

Next, we claim that for all $p\in \{1,\cdots,k_2\}$, it holds
\be\label{qs-14}
t_{p}\leq [Q^{\log_2 4e^2(1+\mu_1)}]+1.
\ee

Indeed, when $p=1$ and $t_1\geq 1$, we set $R_1=e^2d_{\Omega}(x)$ and $a_1=4(1+\mu_1)e^2$.
Then for each $j\in \{1, \ldots, t_1\}$, it holds
$$\frac{R_1}{a_1}=\frac{d_{\Omega}(x)}{4(1+\mu_1)}\stackrel{\eqref{20-07-29-2}}{\leq} \frac{1}{4}d_{\Omega}(y_{1,j}),$$ and thus,
we know from Lemmas \ref{qs-5} and \ref{qs-6}\eqref{qs-6(4)} that
\beqq
t_1\leq \big[Q^{\log_2 4e^2(1+\mu_1)}\big]+1.
\eeqq

When $p\in\{2,\ldots,k_2\}$ and $t_p\geq 1$,
note that for any $u\in \beta_{p}$, we have
\beq\label{20-07-29-4}
d_{\Omega}(u)\stackrel{\eqref{20-07-30-3}}{\geq} \frac{1}{\mu_1}d(x,u)\geq \frac{1}{\mu_1}e^{p-1}d_{\Omega}(x).
\eeq
Set $R_p=e^{p+1}d_{\Omega}(x)$ and $a_p=4\mu_1 e^2$. Then for each $j\in \{1, \ldots, t_p\}$, it holds
$$\frac{R_p}{a_p}=\frac{e^{p-1}d_{\Omega}(x)}{4\mu_1}\stackrel{\eqref{20-07-29-4}}{\leq} \frac{1}{4}d_{\Omega}(y_{p,j}),$$
and thus, again, the desired estimate \eqref{qs-14} follows from Lemmas \ref{qs-5} and \ref{qs-6}\eqref{qs-6(4)}.

Since for each $p\in \{1,\ldots,k_2\}$ and each $j\in \{1,\ldots,t_{p}\}$, $y_{p,j}\in \beta_p\subset \mathcal{B}_p\backslash \mathcal{B}_{p-1}$, and so,
\beq\label{2017-7-7-1}d_{\Omega}(y_{p,j})\leq d(x,y_{p,j})+d_{\Omega}(x)\leq (1+e^p)d_{\Omega}(x).\eeq
This, together with the assertions \eqref{qs-6(7)} and \eqref{qs-6(8)} of Lemma \ref{qs-6}, implies that
\beqq
\ell(\beta)=\sum_{i=1}^{k_1}\ell(\beta[y_i,y_{i+1}]) \leq \frac{11}{9}\sum_{i=1}^{k_1}d_{\Omega}(y_i)\leq \frac{11}{9}\sum_{p=1}^{k_2}\sum_{j=1}^{t_{p}}d_{\Omega}(y_{p,j}).\eeqq

Moreover, since
\[
\sum_{j=1}^{t_{p}}d_{\Omega}(y_{p,j}) \stackrel{\eqref{2017-7-7-1}}{\leq}  (1+e^p)t_pd_{\Omega}(x)\stackrel{\eqref{qs-14}}{\leq} \big([Q^{\log_2 4e^2(1+\mu_1)}]+1\big)(1+e^p)d_{\Omega}(x)
\]
and
\beq\label{20-07-31-1}
d(x,w_0) \stackrel{\eqref{20-07-30-3}}{\leq} \mu_1 d_{\Omega}(w_0)\leq {2}\mu_1 d_{\Omega}(x_0),
\eeq
we get from the estimate of $\ell(\beta)$ that
\beqq
\begin{aligned}
	\ell(\beta)&\leq  \frac{11}{9}\big([Q^{\log_2 4e^2(1+\mu_1)}]+1\big)e^{k_2+1}d_{\Omega}(x)\\
	&\stackrel{\eqref{qs-11}}{\leq} \frac{11}{9}e^2\big([Q^{\log_2 4e^2(1+\mu_1)}]+1\big)d(x,w_0)\\
	&\stackrel{\eqref{20-07-31-1}}{\leq}  {\frac{22}{9}}e^2\mu_1\big([Q^{\log_2 4e^2(1+\mu_1)}]+1\big)d_{\Omega}(x_0)=\lambda d_{\Omega}(x_0).
\end{aligned}
\eeqq

As
\[
\begin{aligned}
	k_1&=\sum_{p=1}^{k_2}t_{p}\stackrel{\eqref{qs-14}}{\leq} \big([Q^{\log_2 4e^2(1+\mu_1)}]+1\big)k_2\\
	&\stackrel{\eqref{qs-11}}{\leq} \big([Q^{\log_2 4e^2(1+\mu_1)}]+1\big)\Big(1+\log \frac{d(x,w_0)}{d_{\Omega}(x)}\Big)\\
	&\stackrel{\eqref{20-07-31-1}}{\leq}  2\mu_1\big([Q^{\log_2 4e^2(1+\mu_1)}]+1\big)\log \frac{3d_{\Omega}(x_0)}{d_{\Omega}(x)},
\end{aligned}
\]
we see from Lemma \ref{qs-6}\eqref{qs-6(0-4)} that
\beqq
k_{\Omega}(x,y) \leq  \sum_{i=1}^{k_1} k_{\Omega}(y_i, y_{i+1})\leq  \frac{40}{27}\mu_1\big([Q^{\log_2 4e^2(1+\mu_1)}]+1\big)\log \frac{3d_{\Omega}(x_0)}{d_{\Omega}(x)}\leq \lambda \log \frac{3d_{\Omega}(x_0)}{d_{\Omega}(x)}.
\eeqq
These prove the lemma in this case, and hence, the proof is complete.
\epf

We shall prove Theorem \ref{2017-10-11-1} by contradiction. Suppose, on the contrary, that
\be\label{H25-09-19-0}
\diam_{\sigma}(\gamma_{x_1x_2})> \theta_0\ell(\alpha_{x_1x_2}).
\ee
Then we know that
\be\label{H25-09-29-0}
\max\{d_{\Omega}(x_1),d_{\Omega}(x_2)\}\leq{\frac{5}{4}}d(x_1,x_2)\leq{\frac{5}{4}}\ell(\alpha_{x_1x_2}).
\ee
Otherwise, by Lemma \ref{Lemma-2.1}, we get
$$\ell(\gamma_{x_1x_2})\leq \frac{50}{9}e^{\frac{40}{9}}d(x_1,x_2)\leq \frac{50}{9}e^{\frac{40}{9}}\ell(\alpha_{x_1x_2}),$$
which contradicts with (\ref{H25-09-19-0}), and so, (\ref{H25-09-29-0}) holds.

Let $x_0\in\gamma_{x_1,x_2}$ be such that
$$\min\{\sigma(x_1,x_0),\sigma(x_2,x_0)\}\geq \frac{1}{3}\diam_{\sigma}(\gamma_{x_1x_2}).$$
Then \be\label{H25-09-21-1}
\min\{\sigma(x_1,x_0),\sigma(x_2,x_0)\}\stackrel{\eqref{H25-09-19-0}}{\geq} \frac{\theta_0}{3}\ell(\alpha_{x_1x_2}).
\ee

Next, by the above two estimates, we may choose $x_{1,0}\in\gamma_{x_1x_2}[x_1,x_0]$ to be the last point along the direction from $x_1$ to $x_0$, and $y_{1,0}\in\gamma_{x_1x_2}[x_2,x_0]$ to be the last point along the direction from $x_2$ to $x_0$, such that
\be\label{H25-09-21-2}
\sigma(x_1,x_{1,0})=\sigma(x_2,y_{1,0})=\frac{3}{2}\ell(\alpha_{x_1x_2}),
\ee
and so,  for each $x\in\gamma_{x_1x_2}[x_{1,0},x_0]$ and each $y\in\gamma_{x_1x_2}[y_{1,0},x_0]$,
\be\label{H25-09-21-3}\sigma(x_2,x)\geq \frac{3}{2}\ell(\alpha_{x_1x_2})\;\quad\mbox{and}\quad\;\sigma(x_1,y)\geq \frac{3}{2}\ell(\alpha_{x_1x_2}).\ee

Now, set $M=16([8Q\theta]^8+1)$, $M_0=2^{\frac{1}{[4Q]}}$, $M_1=6[QM]$, $\theta_{0,1}=(72Q\theta)^4(1+Q^{4\theta})^4$, $\theta_{0,2}=2^M$, $r_1=\ell(\alpha_{x_1x_2})$ and $r_2=5r_1$. Then we take $x_{1,1}(0)={x_{1,0}}$, $y_{1,1}(0)={y_{1,0}}$ and $1_t=1$ for each $t\in\{1,\cdots, M_1\}$.  The following lemma plays a key role in the proof of Theorem \ref{2017-10-11-1}.
\smallskip

\blem\label{25H-09-27} There exists a point sequence $\bigcup\limits_{j=1}^{M-1}\bigcup\limits_{t=j}^{M_1}\bigcup\limits_{i=1}^{j_t}\{x_{j,i}(tr_2)\}$ with the following properties.
\ben
\item\label{CL-09-27-1}  For any $t\in\{1,\cdots, M_1\}$,

{\upshape{(i)}} $x_{1,1}(tr_2)\in\gamma_{x_1x_2}[x_{1,1}(0),x_0]$ is the last point along the direction from $x_{1,1}(0)$ to
$x_0$ such that
$$\sigma(y_{1,0},x_{1,1}(tr_2))=M_0^tr_2;$$

{\upshape{(ii)}} $y_{1,1}(tr_2)\in\gamma_{x_1x_2}[y_{1,0},x_0]$ is the last point along the direction from $y_{1,0}$ to
$x_0$ such that
$$\sigma(x_{1,0},y_{1,1}(tr_2))=\sigma(y_{1,0},x_{1,1}(tr_2));$$

\item\label{CL-09-27-2} For each $j\in\{1,\cdots, M-2)\}$, each $t\in\{j,\cdots, M_1-1\}$ and each $i\in\{1,\cdots,j_t\}$,

{\upshape{(i)}} there exists $x_{j+1,i}^1(tr_2)\in {\alpha_{y_{j,i}(tr_2)x_{j,i}(tr_2)}}$ with $k_{\Omega}(x_{j+1,i}^1(tr_2), x_{j,i}(tr_2))=\frac{1}{2}$;

{\upshape{(ii)}} there exists $x_{j+1,i}^2(tr_2)\in \gamma_{y_{j,i}(tr_2)x_{j+1,i}^1(tr_2)}$ with
$$\min\{\sigma(y_{j,i}(tr_2),x_{j+1,i}^2(tr_2)),\sigma(x_{j+1,i}^1(tr_2),x_{j+1,i}^2(tr_2))\}\geq \frac{1}{3}\diam_{\sigma}\left(\gamma_{y_{1,0}x_{j+1,i}^1(tr_2)}\right);$$

{\upshape{(iii)}} for each $s\in\{t+1,\cdots, M_1\}$, there exist $q\in\{1,\cdots,(j+1)_{s}\}$ and $x_{j+1,q}(sr_2)\in \gamma_{y_{j,i}(tr_2)x_{j+1,i}^1(tr_2)}$ such that $x_{j+1,q}(sr_2)$ is
the last point along the direction from $x_{j+1,i}^1(tr_2)$ to $x_{j+1,i}^2(tr_2)$ with
$$\sigma(y_{j,i}(tr_2),x_{j+1,q}(sr_2))=\sigma(y_{j,i}(tr_2),x_{j,i}(tr_2))+M_0^sr_2,$$
and $y_{j+1,q}(sr_2)\in \gamma_{y_{j,i}(tr_2)x_{j+1,i}^1(tr_2)}$ is
the last point along the direction from $y_{j,i}(tr_2)$ to $x_{j+1,i}^2(tr_2)$ with
$$\sigma(x_{j+1,i}^1(tr_2),y_{j+1,q}(sr_2))=\sigma(y_{j,i}(tr_2),x_{j+1,q}(sr_2)).$$

\item\label{CL-09-27-3}For each $j\in\{2,\cdots,M-1\}$ and $t\in\{j,\cdots, M_1\}$, $j_t=\sum\limits_{s=j-1}^{t-1}(j-1)_s$.
\een
\elem


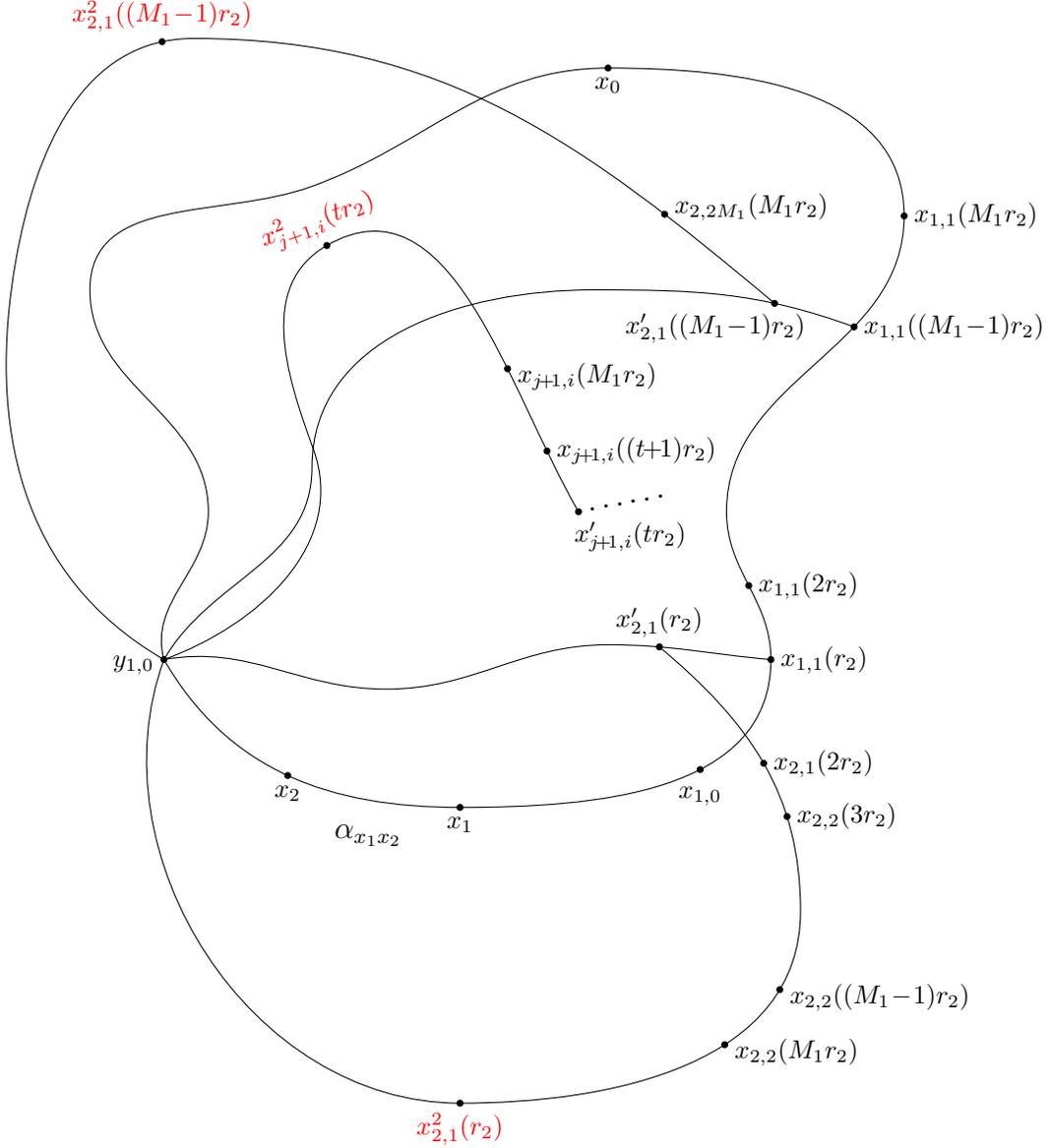
\begin{figure}[htbp]
	\begin{center}
		\begin{tikzpicture}[scale=2]
			
			\coordinate (y10) at (-2,0);
			\filldraw (y10)node[left,yshift=-0.08cm] {\small $y_{1,0}$} circle (0.02);
			\draw (y10)
			to [out=140,in=-110]coordinate[pos=0.7] (y11M1-1r2)   (-2.3,2.5)
			to [out=70,in=200]coordinate[pos=0.15] (y11M1r2)    (-1,3.3)
			to [out=20,in=180]coordinate[pos=1] (x0)   (1,4)
			to [out=0,in=90] coordinate[pos=1] (x11theta02r2)  (3,3)
			to [out=-90,in=90] coordinate[pos=0.35] (x11theta02-1r2)  (1.8,1)
			to [out=-90,in=90] coordinate[pos=0.5] (x112r2) coordinate[pos=1] (x11r2) (2.1,0)
			to [out=-90,in=0] coordinate[pos=0.4] (x10) coordinate[pos=1] (x1)  (0,-1)
			to [out=180,in=-60] node [below, font=\large, pos=0.25,yshift=-0.2cm]{$\alpha_{x_1x_2}$} coordinate[pos=0.5] (x2) (y10);
			
			\filldraw (y11M1-1r2)node[right,xshift=-0.05cm] {\small $y_{1,1}((\!M_1\!\!-\!\!1)r_2)$} circle (0.02);
			\filldraw (y11M1r2)node[below, xshift=0.9cm, yshift=0.1cm] {\small $y_{1,1}(M_1r_2)$} circle (0.02);
			\filldraw (x0)node[below] {\small $x_0$} circle (0.02);
			\filldraw (x11theta02r2)node[right] {\small $x_{1,1}(M_1r_2)$} circle (0.02);
			\filldraw (x11theta02-1r2)node[right] {\small $x_{1,1}((M_1\!-\!1)r_2)$} circle (0.02);
			\filldraw (x112r2)node[right] {\small $x_{1,1}(2r_2)$} circle (0.02);
			\filldraw (x11r2)node[right] {\small $x_{1,1}(r_2)$} circle (0.02);
			\filldraw (x10)node[below,yshift=-0.1cm] {\small $x_{1,0}$} circle (0.02);
			\filldraw (x1)node[below] {\small $x_{1}$} circle (0.02);
			\filldraw (x2)node[below] {\small $x_{2}$} circle (0.02);

			\draw (y10)
			to [out=10,in=180]   (-0.5,-0.2)
			to [out=0,in=180]   (1,0.1)
			to [out=0,in=175] coordinate[pos=0.3] (x21r2') (x11r2);
			
			\filldraw (x21r2')node[above] {\small $x_{2,1}^1(r_2)$} circle (0.02);
			
			\draw (y10)
			to [out=-110,in=180]
			coordinate[pos=0.15] (y212r2)
			coordinate[pos=0.3] (y223r2)
			coordinate[pos=0.65] (y22M1-1r2)
			coordinate[pos=0.8] (y22M1r2)
			coordinate[pos=1] (x21r2^2)  (0,-3)
			to [out=0,in=270] coordinate[pos=0.6] (x22theta02r2) coordinate[pos=0.8] (x22theta02-1r2) (2.3,-1.7)
			to [out=90,in=-40] coordinate[pos=0.3] (x223r2) coordinate[pos=0.5] (x212r2) (x21r2');
			
			\filldraw (y212r2)node[left,xshift=0.05cm,yshift=0cm] {\small $y_{2,1}(2r_2)$} circle (0.02);
			\filldraw (y223r2)node[left,xshift=0.05cm,yshift=0cm] {\small $y_{2,2}(3r_2)$} circle (0.02);
			\filldraw (y22M1-1r2)node[left,xshift=0.1cm,yshift=-0.3cm] {\small $y_{2,2}((\!M_1\!\!-\!\!1)r_2)$} circle (0.02);
			\filldraw (y22M1r2)node[left,xshift=0.1cm,yshift=-0.3cm] {\small $y_{2,2}(M_1r_2)$} circle (0.02);
			\filldraw (x21r2^2)node[below,red] {\small $x_{2,1}^2(r_2)$} circle (0.02);
			\filldraw (x22theta02r2)node[right,yshift=-0.1cm] {\small $x_{2,2}(M_1r_2)$} circle (0.02);
			\filldraw (x22theta02-1r2)node[right,yshift=-0.1cm] {\small $x_{2,2}((M_1\!-\!1)r_2)$} circle (0.02);
			\filldraw (x223r2)node[right] {\small $x_{2,2}(3r_2)$} circle (0.02);
			\filldraw (x212r2)node[right] {\small $x_{2,1}(2r_2)$} circle (0.02);

			\draw (y10)
			to [out=60,in=-90]   (-1,1.3)
			to [out=90,in=180]  coordinate[pos=0.05] (fuzhu)  (0.9,2.5)
			to [out=0,in=160] coordinate[pos=0.7] (x21theta02-1r2')
			(x11theta02-1r2);
			
			\filldraw (x21theta02-1r2')node[below,xshift=-0.8cm] {\small $x_{2,1}^1((M_1\!-\!1)r_2)$} circle (0.02);
			
			\draw (y11M1-1r2)
			to [out=150,in=-100] coordinate[pos=0.3] (y22M1M12^2)  (-3,2.7)
			to [out=80,in=180]  coordinate[pos=0.9] (x21theta02-1r2^2) (-1.8,4.2)
			to [out=0,in=140] coordinate[pos=0.8] (x22theta02r2)
			(x21theta02-1r2');

			\filldraw (y22M1M12^2)node[below,xshift=-1cm] {\small $y_{2,2_{M_1}}(M_1r_2)$} circle (0.02);
			\filldraw (x21theta02-1r2^2)node[above,red,yshift=0cm] {\small $x_{2,1}^2((M_1\!-\!1)r_2)$} circle (0.02);
			\filldraw (x22theta02r2)node[right,yshift=0.1cm] {\small $x_{2,2_{M_1}}(M_1r_2)$} circle (0.02);

			\draw (-1, 0.4)
			to [out=50,in=-70]
			coordinate[pos=0] (yjitr2)
			coordinate[pos=0.7] (yj+1t+1r2)
			(fuzhu)
			to [out=110,in=210] coordinate[pos=0.48] (yj+1M1r2) coordinate[pos=1] (xj+1itr2^2) (-0.9,2.8)
			to [out=30,in=120] coordinate[pos=0.66] (xj+1itheta02r2) coordinate[pos=0.85] (xj+1it+1r2) coordinate[pos=1] (xj+1itr2')
			(0.8,1);
			
			\node[left,rotate=22] at (yjitr2) {\large $\cdots\cdots$};
			\filldraw (yjitr2)node[right,rotate=0,yshift=0cm] {\small $y_{j,i}(tr_2)$} circle (0.02);
			\filldraw (yj+1t+1r2)node[right,rotate=0,yshift=-0.05cm] {\small $y_{j\!+\!1,i}((t+1)r_2)$} circle (0.02);
			\filldraw (yj+1M1r2)node[right,xshift=0.1cm,rotate=20] {\small $y_{j+1,i}(M_1r_2)$} circle (0.02);
			\filldraw (xj+1itr2^2)node[above,xshift=0.1cm, red,rotate=20] {\small $x_{j+1,i}^2(tr_2)$} circle (0.02);
			\filldraw (xj+1itheta02r2)node[right,rotate=0,yshift=-0.1cm] {\small $x_{j\!+\!1,i}(M_1r_2)$} circle (0.02);
			\filldraw (xj+1it+1r2)node[right,rotate=0,yshift=0cm] {\small $x_{j\!+\!1,i}((t\!\!+\!\!1)r_2)$} circle (0.02);
			\filldraw (xj+1itr2')node[below,rotate=0,xshift=0.7cm] {\small $x_{j\!+\!1,i}^1(tr_2)$} circle (0.02);
			\node[right,rotate=10] at (xj+1itr2') {\large $\cdots\cdots$};
			
		\end{tikzpicture}
	\end{center}
	\caption{Illustration for the proof of Lemma \ref{25H-09-27}} \label{Fig-10-5}
\end{figure}

\bpf
For each $t\in\{1,\cdots,M_1\}$, based on \eqref{H25-09-21-1}, Lemma \ref{25H-09-27}(\ref{CL-09-27-1}) holds. 

To prove Lemma \ref{25H-09-27}(\ref{CL-09-27-2}), we shall use an induction argument. Namely, suppose that we have found the points $x_{j,i}(tr_2)$ and $y_{j,i}(tr_2)$ for some $j\in\{1,\cdots, M-2\}$, $t\in\{j,\cdots, M_1-1\}$ and $i\in\{1,\cdots,j_t\}$ as in the lemma, then we may proceed to construct the desired new points $x_{j+1,i}(tr_2)$ and $y_{j+1,i}(tr_2)$ via the following inductive claim.
\bcl\label{Cl25-09-29} For $j\in\{1,\cdots, M-2\}$, $t\in\{j,\cdots, M_1-1\}$ and $i\in\{1,\cdots,j_t\}$, suppose that $x_{j,i}(tr_2)$ and $y_{j,i}(tr_2)$ satisfy the conclusions of Lemma \ref{25H-09-27}. Moreover, we assume that $\gamma_{y_{j,i}(tr_2)x_{j,i}(tr_2)}\in\Lambda_{y_{j,i}(tr_2)x_{j,i}(tr_2)}(\Omega)$  satisfies
$$\diam_{\sigma}(\gamma_{y_{j,i}(tr_2)x_{j,i}(tr_2)})\geq e^{\theta_{0,1}^{8M-\frac{3}{2}j}}r_1,$$
and for each $w\in \gamma_{y_{j,i}(tr_2)x_{j,i}(tr_2)}$, it holds 
$$\sigma(x_{j,i}(tr_2),w)\leq (2^{\frac{1}{[4Q]}}-1)^{-1}\cdot 6\theta d_{\Omega}(w).$$ 
Then there exists some point $x_{j+1,i}^1(tr_2)\in {\alpha_{y_{j,i}(tr_2)x_{j,i}(tr_2)}}$ with the following properties.
\ben
\item\label{Cl25-09-29-1}
$k_{\Omega}(x_{j+1,i}^1(tr_2), x_{j,i}(tr_2))=\frac{1}{2}$.

\item\label{Cl25-09-29-2} For any $\gamma_{y_{j,i}(tr_2)x_{j+1,i}^1(tr_2)}\in\Lambda_{y_{j,i}(tr_2)x_{j+1,i}^1(tr_2)}(\Omega)$,
$$\diam_{\sigma}(\gamma_{y_{j,i}(tr_2)x_{j+1,i}^1(tr_2)})> 2^{-4QM}e^{\theta_{0,1}^{8M-\frac{3}{2}j-1}}r_1,$$
and for any $z\in \gamma_{y_{j,i}(tr_2)x_{j+1,i}^1(tr_2)}$,
$$\sigma(x_{j+1,i}^1(tr_2),z)<5(\theta+(2^{\frac{1}{[4Q]}}-1)^{-1}(6\theta+6\theta^2)) d_{\Omega}(z).$$

\item\label{Cl25-09-29-3} There exists $x_{j+1,i}^2(tr_2)\in \gamma_{y_{j,i}(tr_2)x_{j+1,i}^1(tr_2)}$ such that
$$\sigma(x_{j+1,i}^1(tr_2),x_{j+1,i}^2(tr_2))\geq 2^{-2-4QM}e^{\theta_{0,1}^{8M-\frac{3}{2}j-1}}r_1.$$

\item\label{Cl25-09-29-4} Let  $x_{j+1,i}^2(tr_2)\in \gamma_{y_{j,i}(tr_2)x_{j+1,i}^1(tr_2)}$ satisfy the statement (\ref{Cl25-09-29-3}) above, and select
$$y_1\in \gamma_{y_{j,i}(tr_2)x_{j+1,i}^1(tr_2)}[x_{j+1,i}^1(tr_2),x_{j+1,i}^2(tr_2)].$$
If for each $x\in \gamma_{y_{j,i}(tr_2)x_{j+1,i}^1(tr_2)}[y_1,x_{j+1,i}^2(tr_2)]$, $$\sigma(y_{j,i}(tr_2),x)\geq \sigma(y_{j,i}(tr_2),x_{j,i}(tr_2))+M_0^{t+1}r_2,$$
then for each $y\in\gamma_{y_{j,i}(tr_2)x_{j+1,i}^1(tr_2)}[y_1,x_{j+1,i}^2(tr_2)]$,
$$\sigma(y_1,y)\leq (2^{\frac{1}{[4Q]}}-1)^{-1}\cdot 6\theta d_{\Omega}(y).$$
\een
\ecl

When $j=1$, for each $t\in\{1,\cdots,M_1\}$, we let $j_t=1$, $x_{1,1}^1(r_2)=x_0$, $x_{0,1}(tr_2)=x_{1,0}$ and $y_{0,1}(tr_2)=y_{1,0}$.

For each $j\in\{1,\cdots, M-2)\}$, each $t\in\{j,\cdots, M_1-1\}$ and each $i\in\{1,\cdots,j_t\}$, since the assumptions in this claim implies that $x_{j,i}(tr_2)$ and $y_{j,i}(tr_2)$ satisfy the conclusions of Lemma \ref{25H-09-27}, we obtain from Lemma \ref{25H-09-27}(\ref{CL-09-27-1}) and (\ref{CL-09-27-2})$(iii)$ that there exist $p< t$
and $q\in\{1,\cdots,(j-1)_p\}$ such that
$x_{j,i}(tr_2)\in\gamma_{y_{j-1,q}(pr_2)x_{j,q}^1(pr_2)}[x_{j,q}^1(pr_2),x_{j,q}^2(pr_2)]$ and \be\label{H25-12-08-1}\sigma(y_{j-1,q}(pr_2),x_{j,i}(tr_2))=\sigma(y_{j-1,q}(pr_2),x_{j-1,q}(pr_2))+M_0^tr_2.\ee
Then we have
\be\label{H25-12-08-2}\sigma(x_{j,q}^1(pr_2),x_{j-1,q}(pr_2))\geq d(x_{j,q}^1(pr_2),x_{j-1,q}(pr_2))\geq \frac{1}{4}d_{\Omega}(x_{j-1,q}(pr_2)).
\ee
Otherwise, Lemma \ref{Lemma-2.1} implies that $$k_{\Omega}(x_{j,q}^1(pr_2),x_{j-1,q}(pr_2))\leq \frac{10}{27}<\frac{1}{2},$$
which contradicts with Lemma \ref{25H-09-27} (\ref{CL-09-27-2})$(i)$.

Since $x_{j,q}^1(pr_2)\in \alpha_{y_{j-1,q}(pr_2)x_{j-1,q}(pr_2)}$, we get
$$\ell(\alpha_{y_{j-1,q}(pr_2)x_{j-1,q}(pr_2)})\leq \sigma(y_{j-1,q}(pr_2),x_{j-1,q}(pr_2))+e^{-\theta_1} d_{\Omega}(x_{j-1,q}(pr_2)),$$
and so, by \eqref{H25-12-08-2} and Lemma \ref{25H-09-27} (\ref{CL-09-27-2})$(iii)$, \begin{equation}\label{H25-12-09-1}
	\begin{aligned}
		\sigma(y_{j-1,q}(pr_2),x_{j,q}^1(pr_2)&\leq \ell(\alpha_{y_{j-1,q}(pr_2)x_{j-1,q}(pr_2)})-\sigma(x_{j,q}^1(pr_2),x_{j-1,q}(pr_2))\\ &< \sigma(y_{j-1,q}(pr_2),x_{j-1,q}(pr_2))\leq 2\cdot\sum_{i=1}^p M^i r_2\\ &\leq \frac{2}{M_0-1}(M_0^{p+1}-M_0)r_2.
	\end{aligned}
\end{equation}
Since $p<t$, we have 
$$
\begin{aligned}
	\sigma(y_{j-1,q}(pr_2),x_{j,i}(tr_2)) &\stackrel{\eqref{H25-12-08-1}}{=}\sigma(y_{j-1,q}(pr_2),x_{j-1,q}(pr_2))+M_0^tr_2
	\\ &\geq \sigma(y_{j-1,q}(pr_2),x_{j-1,q}(pr_2))+M_0^{p+1}r_2\\
	&\stackrel{\eqref{H25-12-09-1}}{>}\Big(1+\frac{M_0-1}{2}\Big)\sigma(y_{j-1,q}(pr_2),x_{j-1,q}(pr_2))\\
	&\stackrel{\eqref{H25-12-09-1}}{>} \Big(1+\frac{M_0-1}{2}\Big)\sigma(y_{j-1,q}(pr_2),x_{j,q}^1(pr_2))\\
	&\geq \Big(1+\frac{M_0-1}{2}\Big)(1+e^{-\theta_1})^{-1}\ell(\alpha_{y_{j-1,q}(pr_2),x_{j,q}^1(pr_2)}),
\end{aligned}
$$
and thus, we obtain from Lemma \ref{H25-09-20}(\ref{H25-09-20-2}) and \eqref{H25-12-08-1} that \be\label{H25-12-08-3}M_0^t r_2\leq \sigma(y_{j-1,q}(pr_2),x_{j,i}(tr_2))\leq (M_0-1)^{-1}\cdot 6\theta  d_{\Omega}(x_{j,i}(tr_2)).\ee

Note that based on \eqref{H25-09-21-3}, we may apply  {Lemma \ref{H25-09-20}(\ref{H25-09-20-2})} to obtain
\beqq
\frac{3}{2}\ell(\alpha_{x_1x_2})\stackrel{\eqref{H25-09-21-2}}{=}\sigma(x_2,y_{1,0})\leq 3\theta d_{\Omega}(y_{1,0}),
\eeqq
which, together with (\ref{H25-09-29-0}) and (\ref{H25-09-21-2}), shows
\be\label{H25-12-11-0}\frac{1}{2\theta}\ell(\alpha_{x_1x_2})\leq d_{\Omega}(y_{1,0})\leq \sigma(x_2,y_{1,0})+d_{\Omega}(x_2)\leq \frac{11}{4}\ell(\alpha_{x_1x_2}),\ee
and so, by (\ref{H25-12-08-3}) and Lemma \ref{25H-09-27} (\ref{CL-09-27-2})$(iii)$,
\be\label{H25-09-30-1}
\frac{(2^{\frac{1}{[4Q]}}-1)}{6\theta} r_2\leq d_{\Omega}(x_{j,i}(tr_2))\leq \sigma(y_{1,0},x_{j,i}(tr_2))+d_{\Omega}(y_{1,0})<2(2^{\frac{1}{[4Q]}}-1)^{-1}M_0^{t+1}r_2.
\ee
Then by using the assumption in this claim,

\beq\label{H25-09-29-1}
\begin{aligned}
	k_{\Omega}(y_{j,i}(tr_2),x_{j,i}(tr_2))\stackrel{\eqref{(2.2)}}{\geq} \log\left(1+\frac{\diam_{\sigma}(\gamma_{y_{j,i}(tr_2)x_{j,i}(tr_2)})}{d_{\Omega}(x_{j,i}(tr_2))}\right)>\theta_{0,1}^{8M-\frac{3}{2}j}-4MQ> 1.
\end{aligned}
\eeq
Thus there exists $x_{j+1,i}^1(tr_2)\in \alpha_{y_{j,i}(tr_2)x_{j,i}(tr_2)}$ such that Claim \ref{Cl25-09-29}\eqref{Cl25-09-29-1} holds.

For any $z\in \gamma_{y_{j,i}(tr_2)x_{j+1,i}^1(tr_2)}$, by the conclusion \eqref{Cl25-09-29-1} and assumptions of the claim, we may apply Lemma \ref{H25-09-20}(\ref{H25-09-20-1}) with $x_1=x_{j,i}(tr_2)$, $y_1=x_{j+1,i}^1(tr_2)$ and $\varsigma=(2^{\frac{1}{[4Q]}}-1)^{-1}\cdot 6\theta$ to derive
\be\label{H25-09-29-2}
\sigma(x_{j+1,i}^1(tr_2),z)<5(\theta+(2^{\frac{1}{[4Q]}}-1)^{-1}(6\theta+6\theta^2)) d_{\Omega}(z).
\ee

Let $x_{j+1,i}^0(tr_2)\in \gamma_{y_{j,i}(tr_2)x_{j+1,i}^1(tr_2)}$ be such that
$$
d_{\Omega}(x_{j+1,i}^0(tr_2))\geq \frac{1}{2}\sup _{z\in \gamma_{y_{j,i}(tr_2)x_{j+1,i}^1(tr_2)}} \{d_{\Omega}(z)\}.
$$
Then we obtain from {Lemma \ref{qs-8} (with $\beta=\gamma_{y_{j,i}(tr_2)x_{j+1,i}^1(tr_2)}$) and \eqref{H25-09-29-2}} that
$$k_{\Omega}(y_{j,i}(tr_2),x_{j+1,i}^1(tr_2))\leq \theta_{0,1}\log\Big(1+\frac{3d_{\Omega}(x_{j+1,i}^0(tr_2))}{d_{\Omega}(x_{j+1,i}^1(tr_2))}\Big).$$
Note that by Claim \ref{Cl25-09-29}(\ref{Cl25-09-29-1}) and \eqref{H25-09-29-1}, we have
$$k_{\Omega}(y_{j,i}(tr_2),x_{j+1,i}^1(tr_2))\geq k_{\Omega}(y_{j,i}(tr_2),x_{j,i}(tr_2))-k_{\Omega}(x_{j,i}(tr_2),x_{j+1,i}^1(tr_2))>\theta_{0,1}^{8M-\frac{3}{2}j}-4M Q-\frac{1}{2}.$$
Combining the above two inequalities gives
$$d_{\Omega}(x_{j+1,i}^0(tr_2))>\frac{1}{4}e^{\theta_{0,1}^{8M-\frac{3}{2}j-1}}d_{\Omega}(x_{j+1,i}^1(tr_2)).$$
Notice also that 
$$\log \frac{d_{\Omega}(x_{j,i}(tr_2))}{d_{\Omega}(x_{j+1,i}^1(tr_2))}\stackrel{\eqref{(2.2)}}{\leq} k_{\Omega}(x_{j,i}(tr_2),x_{j+1,i}^1(tr_2))=\frac{1}{2}$$
implies that $$d_{\Omega}(x_{j,i}(tr_2))\leq e\cdot d_{\Omega}(x_{j+1,i}^1(tr_2)).$$
Using the above two inequalities and \eqref{H25-09-30-1}, we obtain that
$$\diam_{\sigma}(\gamma_{y_{j,i}(tr_2)x_{j+1,i}^1(tr_2)})\geq d_{\Omega}(x_{j+1,i}^0(tr_2))-d_{\Omega}(x_{j+1,i}^1(tr_2))\geq 2^{-4QM}e^{\theta_{0,1}^{8M-\frac{3}{2}j-1}}r_1,$$
which, together with (\ref{H25-09-29-2}), shows that Claim \ref{Cl25-09-29}(\ref{Cl25-09-29-2}) holds.

Let $x_{j+1,i}^2(tr_2)\in \gamma_{y_{j,i}(tr_2)x_{j+1,i}^1(tr_2)}$ be such that
$$\min\{\sigma(y_{j,i}(tr_2),x_{j+1,i}^2(tr_2)),\sigma(x_{j+1,i}^1(tr_2),x_{j+1,i}^2(tr_2))\}\geq \frac{1}{3}\diam_{\sigma}(\gamma_{y_{j,i}(tr_2)x_{j+1,i}^1(tr_2)}).$$
Then Claim \ref{Cl25-09-29}(\ref{Cl25-09-29-3}) follows from Claim \ref{Cl25-09-29}(\ref{Cl25-09-29-2}).

Since $x_{j+1,i}^1(tr_2)\in \alpha_{y_{j,i}(tr_2)x_{j,i}(tr_2)}$, we may apply a similar argument as in \eqref{H25-12-09-1} (with $y_{j,i}(tr_2)=y_{j-1,q}(pr_2)$, $x_{j,i}(tr_2)=x_{j-1,q}(pr_2)$ and $x_{j+1,i}^1(tr_2)=x_{j,q}^1(pr_2)$) to obtain that
$$
\sigma(y_{j,i}(tr_2),x_{j+1,i}^1(tr_2))< \sigma(y_{j,i}(tr_2),x_{j,i}(tr_2))< \frac{2}{M_0-1}(M_0^{t+1}-M_0)r_2.
$$
and again by using Claim \ref{Cl25-09-29} (\ref{Cl25-09-29-4}), we know that for each $x\in \gamma_{y_{j,i}(tr_2)x_{j+1,i}^1(tr_2)}[y_1,x_{j+1,i}^2(tr_2)]$,
$$ 
\sigma(y_{j,i}(tr_2),x)\geq\sigma(y_{j,i}(tr_2),x_{j,i}(tr_2)) +M_0^{t+1}r_2>\Big(1+\frac{M_0-1}{2}\Big) \sigma(y_{j,i}(tr_2),x_{j,i}(tr_2)),
$$
and thus, Claim \ref{Cl25-09-29}(\ref{Cl25-09-29-4}) follows from Lemma \ref{H25-09-20}(\ref{H25-09-20-3a}). The proof of the claim is complete.
\smallskip

Now, we may complete the proof of lemma based on Claim \ref{Cl25-09-29} as follows.
We know from (\ref{H25-09-21-2}) that 
$$\sigma(x_{1,0},y_{1,0}) \leq \sigma(x_1,x_{1,0})+\sigma(x_2,y_{1,0})+\sigma(x_1,x_2)\leq 4\ell(\alpha_{xy}),$$
and by (\ref{H25-12-11-0}),
\be\label{H25-12-11-1} \ell(\alpha_{x_{1,0}y_{1,0}})\leq \sigma(x_{1,0},y_{1,0})+e^{-\theta_1}d_{\Omega}(y_{1,0})<5\ell(\alpha_{xy}).\ee
For $j=1$ and $t\in\{1,\cdots, M\}$, it follows from the conclusion \eqref{CL-09-27-1} of the lemma and \eqref{H25-12-11-1} that for each $x\in\gamma_{xy}[x_{1,1}(tr_2),x_0]$,
$$\sigma(y_{1,0},x)\geq\sigma(y_{1,0},x_{1,1}(tr_2))=M_0^t r_2>M_0^t\ell(\alpha_{x_{1,0}y_{1,0}})\geq M_0\ell(\alpha_{x_{1,0}y_{1,0}}),$$
and for each $x\in\gamma_{xy}[y_{1,1}(tr_2),x_0]$,
$$\sigma(x_{1,0},x)\geq\sigma(x_{1,0},y_{1,1}(tr_2))=M_0^t r_2>M_0^t\ell(\alpha_{x_{1,0}y_{1,0}})\geq M_0\ell(\alpha_{x_{1,0}y_{1,0}}).$$
Then we obtain from Lemma \ref{H25-09-20}(\ref{H25-09-20-3a}) that for each $x\in\gamma_{y_{1,0}x_{1,0}}[x_{1,1}(tr_2),y_{1,1}(tr_2)]$, it holds
\be\label{H25-12-04-1}
\sigma(x_{1,1}(tr_2),x)\leq (2^{\frac{1}{[4Q]}}-1)^{-1}\cdot 6\theta  d_{\Omega}(x).
\ee

Take $\gamma_{y_{1,1}(tr_2)x_{1,1}(tr_2)}=\gamma_{x_1x_2}[y_{1,1}(tr_2), x_{1,1}(tr_2)]$. It follows from Lemma \ref{25H-09-27}(\ref{CL-09-27-1}) that
$$\sigma(y_{1,0},x_{1,1}(tr_2))=\sigma(x_{1,0},y_{1,1}(tr_2))= M_0^t r_2,$$
which, together with \eqref{H25-09-21-1}, implies that
\beqq
\begin{aligned}
	\diam_{\sigma}(\gamma_{y_{1,1}(tr_2)x_{1,1}(tr_2)})&\geq  \sigma(x_2,x_0)-\sigma(x_2,x_{1,1}(tr_2))\\
	&\geq \frac{1}{3}\theta_0r_1-\sigma(x_2,y_{1,0})-\sigma(y_{1,0},x_{1,1}(tr_2))\\ 
	&> \frac{1}{3}(\theta_0-6M_0^t)r_1>e^{\theta_{0,1}^{8M-\frac{3}{2}}}r_1.
\end{aligned}
\eeqq

This, together with \eqref{H25-12-04-1}, shows that the condition of Claim \ref{Cl25-09-29} holds for $j=1$. Hence
Claim \ref{Cl25-09-29}\eqref{Cl25-09-29-1}$\sim$\eqref{Cl25-09-29-4} imply Lemma \ref{25H-09-27}(\ref{CL-09-27-2}). Lemma \ref{25H-09-27}(\ref{CL-09-27-3}) is clear from our choice. Hence the proof of lemma is complete.
\epf

Now, we are ready to finish the proof of the theorem.

\bpf[Proof of Theorem \ref{2017-10-11-1}]
By Lemma \ref{25H-09-27}(\ref{CL-09-27-3}) and an elementary computation, we have\be\label{H25-09-22-1} M_2=\mbox{Card}\left\{\bigcup\limits_{j=1}^{M-2}\bigcup\limits_{t=j}^{M_1}\bigcup\limits_{i=1}^{j_t}\{x_{j,i}(tr_2)\}\right\}> (3Q)^{M}.\ee

Take arbitrarily two points $w_1\not=w_2\in\left\{\bigcup\limits_{j=1}^{M-2}\bigcup\limits_{t=j}^{M_1}\bigcup\limits_{i=1}^{j_t}\{x_{j,i}(tr_2)\}\right\}$. Then it follows from  Lemma \ref{25H-09-27} that there exist two integers $p_1, p_2$ in $\{1,\cdots,M_1\}$
such that

$(i)$ for each $i\in\{1,\cdots,p_1-1\}$ and each $j\in\{1,\cdots,p_2-1\}$,
$t_i<t_{i+1}$ and $s_j<s_{j+1}.$

$(ii)$ for each $i\in\{1,\cdots,p_1\}$ and each $j\in\{1,\cdots,p_2\}$,
$t_i$ and $s_j$ are integers in $\{1,\cdots,M_1\}$,
and \beqq
\sigma(y_{1,0},w_1)=\sum_{i=1}^{p_1}M_0^{t_i}r_2\quad\mbox{and }\quad\sigma(y_{1,0},w_2)=\sum_{j=1}^{p_2}M_0^{s_j}r_2.\eeqq
Consequently, we have
\be\label{H25-09-22-6}
\sigma(w_1,w_2)\geq |\sigma(y_{1,0},w_1)-\sigma(y_{1,0},w_2)|\geq 2^{2+\frac{1}{[4Q]}}r_1.
\ee

For any $j\in\{1,\cdots, M-2\}$, $t\in\{j,\cdots, M_1-1\}$ and  $i\in\{1,\cdots,j_t\}$, we see from \eqref{H25-09-30-1} that 
$$d_{\Omega}(x_{j,i}(tr_2))\geq  \frac{2(2^{\frac{1}{[4Q]}}-1) }{3\theta}r_1>\frac{(2^{\frac{1}{[4Q]}}-1)}{2\theta}r_1,$$
and so, $$\overline{\mathbb{B}}\left(x_{j,i}(tr_2),\frac{(2^{\frac{1}{[4Q]}}-1)}{4\theta}r_1\right)\subset \Omega.$$

For each $q\in\{1,\cdots, M_2\}$, set $B_q=\mathbb{B}\left(x_{j,i}(tr_2),\frac{(2^{\frac{1}{[4Q]}}-1)}{4\theta}r_1\right)$. Then we know from \eqref{H25-09-22-6} that for each $q_1\not=q_2\in\{1,\cdots, M_2\}$,
\be\label{H25-09-22-7}\overline{B}_{q_1}\cap \overline{B}_{q_2}=\emptyset.\ee

Moreover, for each $q\in\{1,\cdots,M_2\}$ and $x\in \overline{B}_q$, we infer from Lemma \ref{25H-09-27}(\ref{CL-09-27-2}) that
$$\sigma(y_{0,1},x)<4\theta_{0,2}r_2+\frac{(2^{\frac{1}{[4Q]}}-1)}{4\theta}r_1= \left(\frac{(2^{\frac{1}{[4Q]}}-1)}{4\theta}+16\theta_{0,2}\right) r_1.$$
Thus for each $q\in\{1,\cdots,M_1\}$, it holds
$$\overline{B}_q\subset \mathbb{B}(y_{0,1},(1+16\theta_{0,2}) r_1).$$

Notice that by Lemma \ref{qs-5}, there are at most $Q^{\log_2{(2^{\frac{1}{[4Q]}}-1)^{-1}}\cdot9\theta\theta_{0,2}}$ balls $B_q$ such that
they are disjoint from each other in the ball $\mathbb{B}(y_{0,1},(1+16\theta_{0,2}) r_1)$. This, however, contradicts with \eqref{H25-09-22-1} and \eqref{H25-09-22-7}. The proof of Theorem \ref{2017-10-11-1} is thus complete.
\epf

\subsection{Proof of Theorem \ref{positive-answer} }
Based on Theorem \ref{2017-10-11-1}, we are ready to prove Theorem \ref{positive-answer}. To be more precise, for any $x,y\in \Omega$, $\gamma_{xy}\in \Lambda_{xy}(\Omega)$ and $\alpha_{xy}\in G_{xy}(\Omega)$, we shall prove that
$$\ell(\gamma_{xy})\leq\theta_1\ell(\alpha_{xy})$$
with $\theta_1=(9\theta_0)^2 Q^{2\log_2 56\theta_0}$. For notational simplicity, we write $\gamma=\gamma_{xy}$ and $\alpha=\alpha_{xy}$.

We again use a contradiction argument. Suppose, on the contrary, that
\be\label{10-23-1}
\ell(\gamma)>\theta_1\ell(\alpha).
\ee
Then we have
\begin{equation}\label{2017-08-27-14}
	d(x,y)\geq\frac{3}{4}\max\{d_{\Omega}(x),d_{\Omega}(y)\}.
\end{equation}
Indeed, if not, then Lemma \ref{Lemma-2.1} with $a=\frac{4}{3}$ gives
$$\ell(\gamma)\leq \frac{40}{9}e^{\frac{10}{3}}d(x,y)\leq \theta_1\ell(\alpha),$$
which contradicts with \eqref{10-23-1}.

By \eqref{2017-08-27-14} and Lemma \ref{qs-6}, there exist a finite sequence of balls $\big\{B_\rho=\mathbb{B}(x_\rho,r_\rho)\big\}_{\rho=1}^{t}$ in $\Omega$ and a finite sequence of points $\{x_\rho\}_{\rho=1}^{t+1}$ in $\gamma$ which satisfy all assertions in Lemma \ref{qs-6}. Here, $x_1=x$ and $x_{t+1}=y\in B_{t}$ ( with the possibility that $x_{t+1}=x_t$). Then $\{x_\rho\}_{\rho=1}^{t+1}$ forms a partition of $\gamma$. Furthermore, it follows from Lemma \ref{qs-6}\eqref{qs-6(7)} that
\be\label{2017-08-27-15}
\ell(\gamma)=\sum_{\rho=1}^{t}\ell(\gamma[x_\rho,x_{\rho+1}])\leq \frac{11}{9}\sum_{\rho=1}^{t}d_{\Omega}(x_\rho).
\ee

Let $x_0\in \gamma$ be such that
$$d_{\Omega}(x_0)\geq \frac{1}{2}\sup_{z\in\gamma}\{d_{\Omega}(z)\}.$$
By Theorem \ref{2017-10-11-1}, we have
\be\label{2017-10-12-1}
\sigma(x,x_0)\leq \diam_{\sigma}(\gamma)\leq \theta_0\sigma(x,y).
\ee
To layer the elements in the partition $\{x_\rho\}_{\rho=1}^{t+1}$, we set
\beq\label{20-08-10-3}
S_0=\max_{1\leq \rho\leq t+1}\{d_{\Omega}(x_\rho)\}\;\;\mbox{ and } \;\;T_0=\min_{1\leq \rho\leq t+1}\{d_{\Omega}(x_\rho)\}.
\eeq
Then there must exist an integer $t_2\geq 0$ such that
\be\label{10-23-2} 2^{t_2}\,T_0\leq S_0<2^{t_2+1}\,T_0.\ee

For each $i\in\{0,\ldots,t_2\}$, we define the $i$th layer $A_i$ of the partition of $\gamma$ as
\be\label{17-08-27-01} A_i=\big\{u^i_{j}\in\{x_1,\ldots,x_{t+1}\}:\; 2^{i}\,T_0\leq d_{\Omega}(u^i_{j})<2^{i+1}\,T_0\big\},\ee
and then, set $q_i=\Card\{A_i\}$, with the usual convention that $q_i=0$ if $A_i=\emptyset$, where $``\Card"$ means cardinality.
We define
\beq\label{eq:def for lambdas}
\lambda_1=[Q^{\log_2 56\theta_0}]+1 \quad  \text{and}  \quad \lambda_2=4[Q^{2\log_2 56\theta_0}]+25[Q^{\log_2 56\theta_0}]+36.
\eeq

\smallskip

\textbf{Case 1}: $\max\big\{q_i:\; i\in \{0,1,\ldots,t_2\}\big\}\leq\lambda_2$.
\smallskip

In this case, we have
\beqq
\begin{aligned}
	\ell(\gamma)&\stackrel{\eqref{2017-08-27-15}}{\leq} \frac{11}{9}\sum_{i=0}^{t_2}\sum_{j=0}^{q_i}d_{\Omega}(u_{i,j})
	\leq \frac{11}{9}\lambda_2\sum_{i=0}^{t_2}2^{i+1}T_0
	\stackrel{\eqref{10-23-2}+\eqref{17-08-27-01}}{\leq}  \frac{44}{9}\lambda_2 d_{\Omega}(x_0)\\
	&\leq \frac{44}{9}\left(\sigma(x,x_0)+d_{\Omega}(x) \right)\stackrel{\eqref{2017-08-27-14}+\eqref{2017-10-12-1}}{\leq} \frac{44}{9}\left(\big(\theta_0+\frac{4}{3} \big)\ell(\alpha)\right)< \theta_1\ell(\alpha),
\end{aligned}
\eeqq
which clearly contradicts with \eqref{10-23-1}.

\smallskip

\textbf{Case 2}: $\max\big\{q_i:\; i\in \{0,1,\ldots,t_2\}\big\}>\lambda_2$.
\smallskip

In this case, we need the following result, whose lengthy proof will be postponed to Subsection \ref{subsec:key proposition}.
\bprop\label{2017010-20-1}
Suppose that $\max\big\{q_i:\; i\in \{0,1,\ldots,t_2\}\big\}>\lambda_2$. Then there are partitions $P_{\gamma}=\{u_j\}_{j=0}^{s_0+1}\subset \gamma$ and $P_{\alpha}=\{w_j\}_{j=0}^{s_0+1}\subset \alpha$, where $u_0=w_0=x$ and $u_{s_0+1}=w_{s_0+1}=y$, such that the following conclusions hold:
\ben
\item
$s_0\geq \lambda_1$.
\item\label{17-10-11-9}
For each $j\in\{0,1,\ldots,s_0\}$ and every $i\in\{0,\ldots,t_2\}$,
$${\Card}\big\{\gamma[u_j,u_{j+1}]\cap A_i\big\}\leq \lambda_2.$$
\item\label{17-10-11-7}
For each $j\in\{1,\ldots,s_0-1\}$,
$$\sigma(u_j,u_{j+1})\geq 3\theta_0\max\{d_{\Omega}(u_j),d_{\Omega}(u_{j+1})\}$$ and $$\sigma(w_j,w_{j+1})\geq \frac{3}{2}\theta_0\max\{d_{\Omega}(u_j),d_{\Omega}(u_{j+1})\}.$$
\item\label{2017-10-19-1}
For each $j\in\{1,\ldots,s_0\}$, $$ w_j\in \alpha[w_{j-1},w_{j+1}]\;\;\mbox{ and }\;\;\sigma(u_j,w_j)\leq \theta_0 d_{\Omega}(u_j).$$
\item\label{2017-10-13-1}
$\sigma(w_0,w_1)\geq 3\theta_0 d_{\Omega}(u_1)\;\mbox{ and }\; \sigma(w_{s_0},w_{s_0+1})\geq 3\theta_0 d_{\Omega}(u_{s_0})$.
\een
\eprop

Here and hereafter, the elements of a sequence of points on $\gamma=\gamma_{xy}$ or $\alpha=\alpha_{xy}$ are consecutively listed along the direction from $x$ to $y$.

With the aid of Proposition \ref{2017010-20-1}, we are ready to complete the proof of Theorem \ref{positive-answer}.
\bpf[Proof of Theorem \ref{positive-answer}]
Since $\max\big\{q_i:\; i\in \{0,1,\ldots,t_2\}\big\}>\lambda_2$, by Proposition \ref{2017010-20-1}, there are partitions $P_{\gamma}=\{u_j\}_{j=0}^{s_0+1}\subset\gamma$ and $P_{\alpha}=\{w_j\}_{j=0}^{s_0+1}\subset \alpha$ such that all conclusions of Proposition \ref{2017010-20-1} are satisfied, where $v_0=w_0=x$ and $v_{s_0+1}=w_{s_0+1}=y$.

Next, we shall prove that  for each $j\in\{0,1,\ldots,s_0\}$, it holds
\beq\label{black}
\ell(\gamma[u_j,u_{j+1}])\leq \frac{88}{27\theta_0}\lambda_2(3\theta_0^2+1)\ell(\alpha[w_j,w_{j+1}]).
\eeq

For each $j\in\{0,1,\ldots,s_0\}$, Theorem \ref{2017-10-11-1} implies that
\beq\label{eq:for length of gamma}
\diam_{\sigma}(\gamma_{xy}[u_j,u_{j+1}])\leq \theta_0\sigma(u_j,u_{j+1}).
\eeq
Set
$$d_{\Omega}(v_j)=\max\left\{d_{\Omega}(u): u\in \gamma[u_j, u_{j+1}]\bigcap \Big(\bigcup_{i=0}^{t_2}A_i\Big)\right\}.$$
Then there is an integer $t_3\geq 0$ such that
\be\label{20-08-10-1} 2^{t_3}\,T_0\leq d_{\Omega}(v_j) <2^{t_3+1}\,T_0,\ee
where $T_0$ is given by \eqref{20-08-10-3}. By \eqref{10-23-2}, $t_3\leq t_2+1$. Still, the proof of \eqref{black} requires the following claim.
\smallskip

\textbf{Claim}: For each $j\in\{0,1,\ldots,s_0\}$, $$\ell(\gamma_{xy}[u_j,u_{j+1}])\leq \frac{44}{9}\lambda_2 d_{\Omega}(v_j).$$
\smallskip

For the proof of this claim, let
$$\{v_{j,s}\}_{s=0}^{p}=\gamma_{xy}[u_j,u_{j+1}]\cap \{x_\rho\}_{\rho=1}^{t+1},$$
where $v_{j,0}=u_j$ and $v_{j,p+1}=u_{j+1}$. Set
\beq\label{20-08-10-2}
A_{j}^i=\big\{u^i_{j,m}\in\{v_{j,0},\ldots, v_{j,p}:\; 2^{i}\,T_0\leq d_{\Omega}(u^i_{j,m})<2^{i+1}\,T_0\big\}
\eeq
and ${\Card}\{A^i_{j}\}=q_{j,i}$.
Then by Proposition \ref{2017010-20-1}\eqref{17-10-11-9}, $q_{j,i}\leq \lambda_2$, and thus,
\beqq
\begin{aligned}
	\ell(\gamma[u_j,u_{j+1}])&= \sum_{s=0}^{p}\ell(\gamma[v_{j,s},v_{j,s+1}])\stackrel{\text{Lemma } \ref{qs-6} \eqref{qs-6(7)}}{\leq} \frac{11}{9}\sum_{s=0}^{p}d_{\Omega}(v_{j,s})\leq \frac{11}{9}\sum_{i=0}^{t_3}\sum_{m=0}^{q_{j,i}}d_{\Omega}(u_{j,m}^i) \\
	&\stackrel{\eqref{20-08-10-2}}{\leq} \frac{11}{9}\lambda_2\sum_{i=0}^{t_3}2^{i+1}T_0 \stackrel{\eqref{20-08-10-1}}{\leq}\frac{44}{9}\lambda_2 d_{\Omega}(v_j).
\end{aligned}
\eeqq
This completes the proof of the claim.
\medskip

Let us continue the proof of \eqref{black} based on the above claim.
Since $v_j\in \gamma[u_{j},u_{j+1}]$, we have
\begin{align*}
	d_{\Omega}(v_j)&\leq \min\big\{\sigma(u_j,v_j)+d_{\Omega}(u_j),\sigma(u_{j+1},v_j)+d_{\Omega}(u_{j+1})\big\}
	\\ &\leq \diam_{\sigma}(u_{j},u_{j+1})+\min\big\{d_{\Omega}(u_j),d_{\Omega}(u_{j+1})\big\}.
\end{align*}
This, combining with the above claim and \eqref{eq:for length of gamma}, gives
\be\label{17-08-28-07}
\ell(\gamma[u_j,u_{j+1}])\leq
\frac{44}{9}\lambda_2 \big(\theta_0\sigma(u_j,u_{j+1})+\min\{d_{\Omega}(u_j),d_{\Omega}(u_{j+1})\}\big).
\ee

We split the arguments into the following two cases.

\textbf{Case A}: $j\in\{1,\ldots,s_0-1\}$.
\smallskip

By the assertions \eqref{17-10-11-7} and \eqref{2017-10-19-1} of Proposition \ref{2017010-20-1}, we have
$$\sigma(w_j,w_{j+1})\geq \sigma(u_j,u_{j+1})-\sigma(u_j,w_j)-\sigma(u_{j+1},w_{j+1})\geq \frac{1}{3}\sigma(u_j,u_{j+1}),$$
which, together with the assertions \eqref{17-10-11-7} and \eqref{2017-10-13-1} of Proposition \ref{2017010-20-1}, shows that
$$\ell(\gamma[u_j,u_{j+1}])\stackrel{\eqref{17-08-28-07}}{\leq}\frac{44}{9}\lambda_2\Big(3\theta_0+\frac{2}{\theta_0}\Big)\sigma(w_j,w_{j+1})<\frac{88}{27\theta_0}\lambda_2(3\theta_0^2+1)\ell(\alpha[w_j,w_{j+1}]).$$

\smallskip

\textbf{Case B}: $j\in\{0,s_0\}$.
\smallskip

We only consider the case $j=0$, as the proof for the other case is similar.

If $d_{\Omega}(u_1)\leq \frac{1}{2\theta_0}\sigma(u_0,u_1)$, then $$d_{\Omega}(u_0)\leq d_{\Omega}(u_1)+\sigma(u_0,u_1)\leq \frac{2\theta_0+1}{2\theta_0}\sigma(u_0,u_1).$$

Since $w_0=u_0=x$ and the assumption $d_{\Omega}(u_1)\leq \frac{1}{2\theta_0}\sigma(u_0,u_1)$ of the case, it follows from Proposition \ref{2017010-20-1}(\ref{2017-10-19-1}) that
$$\sigma(w_0,w_1)\geq \sigma(u_0,u_1)-\sigma(u_1,w_1)\geq \frac{1}{2}\sigma(u_0,u_1).$$
This gives $$d_{\Omega}(u_1)\leq\frac{1}{\theta_0}\sigma(w_0,w_1),$$ and thus, we have
$$\ell(\gamma[u_0,u_1])\stackrel{\eqref{17-08-28-07}}{\leq} \frac{44}{9}\lambda_2\Big(2\theta_0+\frac{1}{\theta_0}\Big)\sigma(w_0,w_1)\leq \frac{44}{9\theta_0}\lambda_2(2\theta_0^2+1)\ell(\alpha[w_0,w_1]).$$

For the remaining case, that is, $d_{\Omega}(u_1)> \frac{1}{2\theta_0}\sigma(u_0,u_1)$, again, by Proposition \ref{2017010-20-1}(\ref{2017-10-13-1}) and \eqref{17-08-28-07}, we obtain
$$\ell(\gamma[u_0,u_1])\leq \frac{44}{9}\lambda_2 (2\theta_0^2+1)d_{\Omega}(u_1)\leq \frac{44}{27\theta_0}\lambda_2 (2\theta_0^2+1)\ell(\alpha[w_0,w_1]).$$

In either case, we have proved \eqref{black}.
\medskip

Now, it follows from \eqref{black} that
$$ \ell(\gamma)=\sum_{j=1}^{s_0+1}\ell(\gamma[u_{j-1},u_j])\leq \frac{88}{27\theta_0}\lambda_2(3\theta_0^2+1)\ell(\alpha)<\theta_1\ell(\alpha),$$
which again contradicts with \eqref{10-23-1}. Thus the proof of Theorem \ref{positive-answer} is complete.
\epf
\medskip

	\subsection{Proof of Proposition \ref{2017010-20-1}}\label{subsec:key proposition}
	In this subsection, we present the proof of Proposition \ref{2017010-20-1}. It requires a couple of auxiliary lemmas. Let us recall that the $i$-th layer $A_i$ of the partition of $\gamma$ is defined in \eqref{17-08-27-01}, and the constants $\lambda_1$ and $\lambda_2$ are defined in \eqref{eq:def for lambdas}.
	
	\blem\label{Mod-again-01} For each $i\in\{0,\ldots,t_2\}$ and any $u\in\gamma$, it holds
	$${\Card}\{\mathbb{B}(u,  3\cdot2^{i+2}\theta_0 T_0)\cap A_i\}\leq \lambda_1.$$
	\elem
	
	\bpf Let $v\in A_i\cap \mathbb{B}(u, 3\cdot2^{i+2}\theta_0 T_0)$. Then for any $w\in \overline{\mathbb{B}}(v, 4^{-1}d_{\Omega}(v))$,
	by (\ref{17-08-27-01}), we have $$d(w,u)\leq d(w,v)+d(v,u)< (6\theta_0+4^{-1})2^{i+1}T_0<7\cdot 2^{i+1}\theta_0 T_0.$$
	This implies that
	$$\overline{\mathbb{B}}(v, 4^{-1}d_{\Omega}(v))\subset \mathbb{B}(u, 7\cdot2^{i+1}\theta_0 T_0).$$
	
	Taking into consideration of Lemma \ref{qs-6}(\ref{qs-6(4)}) and \eqref{17-08-27-01}, the lemma follows directly from Lemma \ref{qs-5} (applied with $R=7\cdot2^{i+1}\theta_0 T_0$, $a=56\theta_0$ and $r=\frac{R}{a}=2^{i-2}T_0$).
	\epf

	The next lemma gives a useful partition of the $i$th layer $A_i$ of $\gamma$.
	
	\blem\label{Mod-again-03}
	For $v_1\in\gamma$ and $v_2\in\gamma[v_1,x_{t+1}]$, suppose that there exists some $i\in\{0,\ldots,t_2\}$ such that
	\ben
	\item[(i)]
	the set $E_i=\gamma[v_1,v_2]\cap A_i=\{u_r\}_{r=1}^q$ with $q>\lambda_2$.
	\item[(ii)]
	there are $w_1\in\alpha$ and $w_2\in\alpha[w_1,y]$ such that for any $u\in A_{i}$, $$\max\{\sigma(v_1,w_1),\sigma(v_2,w_2)\}\leq \theta_0 d_{\Omega}(u).$$
	\een
	Then there exists a partition $\{u_{p_j}\}_{j=1}^{\iota}\subset E_i$  such that
	\ben
	\item\label{2017-10-08-0}
	$\iota\geq \lambda_1$.
	\item\label{20-08-4-1}
	for each $j\in\{0,1,\ldots,\iota\}$,
	$${\Card}\big\{\gamma[u_{p_j},u_{p_{j+1}}]\cap E_i \big\}\leq \lambda_1+1,$$ where $u_{p_0}=v_1$ and $u_{p_{\iota+1}}=v_2$.
	\item\label{2017-10-08-1-0}
	$\min\{\sigma(v_1,u_{p_1}),\sigma(v_2,u_{_{p_{\iota}}})\}\geq 6\theta_0 d_{\Omega}(u_{p_1}).$
	\item\label{2017-10-08-1-10}
	$\min\{\sigma(u_{p_j},u_{p_{j+1}}):\; j\in \{1, \ldots, \iota-1\}\}\geq 6\theta_0 d_{\Omega}(u_{p_j}).$
	\item\label{2017-10-18-1}
	for each $j\in\{1,\ldots,\iota\}$ and $\alpha_{v_2w_2}\in G_{v_2w_2}(\Omega)$,
	$$\sigma(u_{p_j},\alpha_{v_2w_2})>\theta_0 d_{\Omega}(u_{p_j}).$$
	\een\elem
	\bpf
	Since $E_i\subset A_i$, Lemma \ref{Mod-again-01} gives
	$${\rm Card}\big\{\mathbb{B}(v_1,3\cdot2^{i+2}\theta_0 T_0)\cap E_i \big\}\leq \lambda_1.$$
	It follows that there are at least $k_1\geq \lambda_2-\lambda_1$ points $\{u_r\}_{r=1}^{k_1}$ in $E_i$, which are not contained in $\mathbb{B}(v_1,3\cdot2^{i+2}\theta_0 T_0)$.
	Let
	$$
	s_1:=\min\{r\in \{1,\cdots,k_1\}:u_r\notin B(v_1,3\cdot 2^{i+2}\theta_0 T_0)\}.
	$$  
	Then we have
	$${\rm Card}\big\{ \gamma[v_1,u_{s_1}]\cap E_i\big\}\leq \lambda_1+1$$
	and
	\beqq
	\sigma(v_1,u_{s_1})\geq d(v_1,u_{s_1})\geq  3\cdot2^{i+2}\theta_0 T_0\stackrel{\eqref{17-08-27-01}}{\geq} 6\theta_0 d_{\Omega}(u_{s_1}).
	\eeqq
	
	Applying Lemma \ref{Mod-again-01} again, we obtain $\Card\big\{\mathbb{B}(u_{s_1},3\cdot2^{i+2}\theta_0 T_0)\cap E_i \big\}\leq \lambda_1$, and thus,
	there are at least $k_2\geq\lambda_2-2\lambda_1-1$ points $\{u_r\}_{r=s}^{k_2+s-1}$ in $E_i$ with $s>s_1$,
	which are not contained in $\mathbb{B}(u_{s_1},3\cdot2^{i+2}\theta_0 T_0)$.
	Let
	$$
	s_2:=\min\{r\in \{s,\cdots,k_2+s-1\}:u_r\notin B(u_{s_1},3\cdot 2^{i+2}\theta_0 T_0)\}.
	$$
	Then
	$${\Card}\big\{ \gamma[u_{s_1},u_{s_2}]\cap E_i\big\}\leq \lambda_1+1$$
	and
	\beqq
	\sigma(u_{s_1},u_{s_2})\geq d(u_{s_1},u_{s_2})\geq 3\cdot2^{i+2}\theta_0 T_0 \stackrel{\eqref{17-08-27-01}}{\geq} 6\theta_0 d_{\Omega}(u_{s_1}).
	\eeqq
	
	Repeating this procedure, we may find a finite sequence of points $\{u_{s_r}\}_{r=1}^{\iota_1}\subset E_i$ such that
	$$ \iota_1\geq \frac{q}{\lambda_1+1}\geq 4[Q^{\log_2 56\theta_0}]+8\;\;\mbox{ and }\;\;{\rm Card}\big\{\gamma[u_{s_{\iota_1}+1},v_2]\cap E_i \big\}\leq \lambda_1+1.$$
	Applying Lemma \ref{Mod-again-01}, once again, we infer that
	$$ {\rm Card}\big\{\mathbb{B}(v_2,3\cdot2^{i+2}\theta_0 T_0)\cap E_i \big\}\leq \lambda_1\;\;\mbox{ and }\;\;
	{\rm Card}\big\{\mathbb{B}(v_2,3\cdot2^{i+2}\theta_0 T_0)\cap \{u_{s_r}\}_{r=1}^{\iota_1} \big\}\leq \lambda_1.$$
	This shows that there are $\iota$ points in $\{u_{s_r}\}_{r=1}^{\iota_1}$, which are not contained in $\mathbb{B}(v_2,3\cdot2^{i+2}\theta_0 T_0)$.
	Denote these $\iota$ points by $\{u_{p_j}\}_{j=1}^{\iota}$. Then it follows from the preceeding construction that
	
	\begin{enumerate}
		\item\label{2017-10-11-2}
		$\iota\geq \iota_1-\lambda_1-1 \geq 2[Q^{\log_2 56\theta_0}]+5\geq \lambda_1.$
		\item
		for each $j\in \{0,1,\ldots, \iota\},$ $${\rm Card}\big\{\gamma[u_{p_j},u_{p_{j+1}}]\cap E_i \big\}\leq \lambda_1+1,
		$$  where $u_{p_0}=v_1$ and $u_{p_{\iota+1}}=v_2$.
		\item
		$\min\{\sigma(v_1,u_{p_1}),\sigma(v_2,u_{_{p_{\iota}}})\}\geq 6\theta_0 d_{\Omega}(u_{p_1}).$
		\item
		$\min\{\sigma(u_{p_j},u_{p_{j+1}}):\; j\in \{1, \ldots, \iota-1\}\}\geq 6\theta_0 d_{\Omega}(u_{p_j}).$
	\end{enumerate}
	
	To finish the proof of lemma, it remains to verity the statement \eqref{2017-10-18-1}.  Suppose, on the contrary, that there exists $j\in\{1,\ldots,\iota\}$ such that
	$$\sigma(u_{p_j},\alpha_{v_2w_2})\leq \theta_0  d_{\Omega}(u_{p_j}).$$
	
	On the one hand, we know from the assumption $\rm{(ii)}$ of lemma that
	\[
	\begin{aligned}
		\sigma(u_{p_{j}},v_2)&\leq \sigma(u_{p_{j}},w_2)+\sigma(w_2,v_2)\leq \sigma(u_{p_j},\alpha_{v_2w_2})+\ell(\alpha_{v_2w_2})+\sigma(w_2,v_2)\\
		&\leq \theta_0 d_{\Omega}(u_{p_j})+(2+e^{-\theta_1})\sigma(w_2,v_2)<4\theta_0 d_{\Omega}(u_{p_j}).
	\end{aligned}
	\]
	On the other hand, for each $j\in \{1, \ldots, \iota\}$,
	$$\mathbb{B}(v_2,3\cdot2^{i+2}\theta_0 T_0) \cap  \{u_{p_j}\}_{j=1}^{\iota} =\emptyset,$$
	and thus,
	$$\sigma(u_{p_{j}},v_2)\geq d(u_{p_{j}},v_2)\geq 3\cdot2^{i+2}\theta_0 T_0\stackrel{\eqref{17-08-27-01}}{\geq} 6\theta_0 d_{\Omega}(u_{p_j}).$$
	This is a contradiction, and hence, the proof of lemma is complete.
	\epf
	
	Corresponding to the partition of $\gamma$ in Lemma \ref{Mod-again-03}, we have an associated partition on $\alpha$.
	\blem\label{Mod-again-10-16-1} Under the assumptions of Lemma \ref{Mod-again-03}, there exists a partition $\{w_{p_j}\}_{j=1}^{\iota}\subset \alpha[w_1,w_2]$ such that
	\ben
	\item\label{2017-10-15-1}
	for each $j\in\{1,\ldots,\iota\}$,
	$$w_{p_j}\in\alpha[w_{p_{j-1}},w_{p_{j+1}}]\;\;\mbox{ and }\;\;\sigma(u_{p_j},w_{p_j})\leq \theta_0 d_{\Omega}(u_{p_j}),$$ where $w_{p_{0}}=w_1$ and $w_{p_{\iota+1}}=w_2$.
	\item\label{Mod-again-03-03}
	for each $j\in\{0,1,\ldots, \iota\}$, $$\sigma(w_{p_{j}},w_{p_{j+1}})\geq 3\theta_0 d_{\Omega}(u_{p_1}).$$
	\een
	\elem
	\bpf (1) Let $\{u_{p_j}\}_{j=1}^{\iota}\subset \gamma[v_1,v_2]$ be the point sequence obtained in Lemma \ref{Mod-again-03}. Fix $\alpha_{v_1w_1}\in G_{v_1w_1}(\Omega)$ and  $\alpha_{v_2w_2}\in G_{v_2w_2}(\Omega)$.
	
	We first consider the point $u_{p_1}$. Since $(\Omega,\sigma)$ satisfies the $\theta$-ball separation condition, there exists a point $w_{p_1}\in \alpha_{v_1w_1}\cup\alpha[w_1,w_2]\cup \alpha_{v_2w_2}$  such that
	\be\label{20-08-5-1}
	\sigma(u_{p_1},w_{p_1})\leq \theta d_{\Omega}(u_{p_1})\leq \theta_0 d_{\Omega}(u_{p_1}).
	\ee
	It follows from Lemma \ref{Mod-again-03}(\ref{2017-10-18-1}) that
	$w_{p_1}\in\alpha_{v_1w_1}\cup\alpha[w_1,w_2]$.
	
	Suppose that $w_{p_1}\in\alpha_{v_1w_1}$. Since by the assumption
	$\rm{(ii)}$ of Lemma \ref{Mod-again-03}, we obtain that
	$$
	\sigma(v_1,u_{p_1})\leq\sigma(v_1,w_{p_1})+\sigma(w_{p_1},u_{p_1})\stackrel{\eqref{20-08-5-1}}{\leq}\ell(\alpha_{v_1w_1})+\theta_0 d_{\Omega}(u_{p_1})\leq  (2+e^{-\theta_1})\theta_0 d_{\Omega}(u_{p_1}),$$
	which contradicts with Lemma \ref{Mod-again-03}\eqref{2017-10-08-1-0}.
	This shows that $w_{p_1}\in\alpha[w_1,w_2]$.
	
	Fix $\alpha_{u_{p_1}w_{p_1}}\in \Gamma_{u_{p_1}w_{p_1}}(\Omega)$ and $\alpha_{u_{p_2}w_{p_2}}\in \Gamma_{u_{p_2}w_{p_2}}(\Omega)$. Let us consider the point $u_{p_2}$. Since $(\Omega,\sigma)$ satisfies the $\theta$-ball separation condition, there exists a point
	$w_{p_2}\in \alpha_{u_{p_1}w_{p_1}}\cup\alpha[w_{p_1},w_2]\cup \alpha_{v_2w_2}$  such that
	\be\label{20-08-6-1}
	\sigma(u_{p_2},w_{p_2})\leq \theta d_{\Omega}(u_{p_2})\leq \theta_0 d_{\Omega}(u_{p_2}).
	\ee
	It follows from Lemma \ref{Mod-again-03}(\ref{2017-10-18-1}) that
	$w_{p_2}\in\alpha_{u_{p_1}w_{p_1}}\cup\alpha[w_{p_1},w_2]$.
	
	Suppose that $w_{p_2}\in \alpha_{u_{p_1}w_{p_1}}$. Then we have
	\beqq
	\begin{aligned}
		\sigma(u_{p_1},u_{p_2}) &\leq  \sigma(u_{p_1},w_{p_2})+\sigma(w_{p_2},u_{p_2})\stackrel{\eqref{20-08-6-1}}{\leq}  \ell(\alpha_{u_{p_1}w_{p_1}})+\theta_0 d_{\Omega}(u_{p_2})\\
		&\leq (1+e^{-\theta_1})\sigma(u_{p_1},w_{p_1})+\theta_0 d_{\Omega}(u_{p_2})\stackrel{\eqref{17-08-27-01}+\eqref{20-08-5-1}}{\leq} (3+2e^{-\theta_1})\theta_0 d_{\Omega}(u_{p_2}),
	\end{aligned}
	\eeqq
	which contradicts with Lemma \ref{Mod-again-03}\eqref{2017-10-08-1-10}.
	This shows that $w_{p_2}\in\alpha[w_{p_1},w_2]$.
	
	Repeating this procedure, we may find a finite sequence of points $\{w_{p_j}\}_{j=1}^{\iota}\subset \alpha[w_1,w_2]$ which satisfies the first statement of Lemma~\ref{Mod-again-10-16-1}.
	
	(2) Since by the assumption
	$$\max\{\sigma(v_1,w_1),\sigma(v_2,w_2)\}\leq\theta_0 d_{\Omega}(u_{p_1}),$$
	we obtain from Lemma \ref{Mod-again-03}\eqref{2017-10-08-1-0} and the conclusion \eqref{2017-10-15-1} of the lemma that
	\beqq
	\sigma(w_1,w_{p_1}) \geq \sigma(v_1,u_{p_1})-\sigma(v_1,w_1)-\sigma(u_{p_1},w_{p_1}) \geq  4\theta_0 d_{\Omega}(u_{p_1}),
	\eeqq
	and  that
	\beqq
	\sigma(w_2,w_{p_{\iota}}) \geq  \sigma(v_2,u_{p_{\iota}})-\sigma(u_{p_{\iota}},w_{p_{\iota}})-\sigma(v_2,w_2)
	\stackrel{\eqref{17-08-27-01}}{\geq}  3\theta_0 d_{\Omega}(u_{p_1}).
	\eeqq
	Moreover, for each $j\in\{2,\ldots,\iota\}$, it follows from Lemmas \ref{Mod-again-03}\eqref{2017-10-08-1-10} and  \ref{Mod-again-10-16-1}\eqref{2017-10-15-1} that
	\beqq
	\sigma(w_{p_{j-1}},w_{p_j}) \geq  \sigma(u_{p_{j-1}},u_{p_j})-\sigma(u_{p_j},w_{p_j})-\sigma(u_{p_{j-1}},w_{p_{j-1}})\stackrel{\eqref{17-08-27-01}}{\geq} 3\theta_0 d_{\Omega}(u_{p_1}).
	\eeqq
	This proves the second statement, and hence, the proof of lemma is complete.
	\epf
	
	Now, we are ready to prove Proposition \ref{2017010-20-1}.
	
	\bpf[Proof of Proposition \ref{2017010-20-1}]
	Set $$r_1=\min\big\{i\in\{0,1,\ldots,t_2\}:\; q_i>\lambda_2\big\},$$
	$v_1=w_1=x$ and $v_2=w_2=y$. Then by Lemmas \ref{Mod-again-03} and \ref{Mod-again-10-16-1}, there are partitions $\{u^{1}_{\nu_1}\}_{\nu_1=0}^{\iota_1+1}$ of $\gamma$ and $\{w^{1}_{\nu_1}\}_{\nu_1=0}^{\iota_1+1}$ of $\alpha$ such that all conclusions in Lemmas \ref{Mod-again-03} and \ref{Mod-again-10-16-1} are satsified, where $u_{0}^1=x$, $u_{\nu_1+1}^1=y$, $w^{1}_{0}=x$, $w^{1}_{\nu_1+1}=y$, $\{u^{1}_{\nu_1}\}_{\nu_1=1}^{\iota_1}\subset A_{r_1}$ and $\iota_1\geq \lambda_1$.
	
	If for each $\nu_1\in\{0,1, \ldots, \iota_1\}$ and for every $i\in  \{r_1, \ldots, t_2\}$,
	$${\Card}\big\{\gamma[u^{1}_{\nu_1}, u^{1}_{\nu_1+1}]\cap A_i\big\}\leq \lambda_2,$$
	then for all $i\in \{0,1, \ldots, t_2\}$, we have
	\beq\label{25-12-2-1}
	{\Card} \big\{ \gamma[u^{1}_{\nu_1},u^{1}_{\nu_1+1}]\cap A_{i} \big\}\leq \lambda_2.
	\eeq
	
	Let
	$$P_{\gamma}^{1}=\left\{u^{1}_{\nu_1}\right\}_{\nu_1=0}^{\iota_1+1}\;\;\mbox{ and }\;\; P_{\alpha}^{1}=\left\{w^{1}_{\nu_1}\right\}_{\nu_1=0}^{\iota_1+1}.$$
	
	Since $P_{\gamma}^{1}$ and $P_{\alpha}^{1}$ satisfy
	all conclusions in Lemmas \ref{Mod-again-03} and \ref{Mod-again-10-16-1}, and since $\iota_1\geq \lambda_1$, by letting $s_0=\iota_1$, we see from \eqref{25-12-2-1} that $P_{\gamma}^{1}$ and $P_{\alpha}^{1}$
	are our desired partitions.
	
	Otherwise, there are some $\nu_1\in\{0,1, \ldots, \iota_1\}$ and an $i\in \{r_1+1, \ldots, t_2\}$ such that
	$${\Card}\big\{\gamma[u^{1}_{\nu_1}, u^{1}_{\nu_1+1}]\cap A_i\big\}> \lambda_2.$$
	Set
	$$r_2=\min\big\{i\in\{r_1+1,\ldots,t_2\}:\; q_i>\lambda_2\big\}.$$ Clearly, \beqq\label{20-08-7-1}r_2>r_1.\eeqq
	
	Assume that there are $K\geq 1$ sub-curves $\gamma[u^{1}_{\nu_{1,1}}$, $u^{1}_{\nu_{1,1}+1}]$, $\cdots$, $\gamma[u^{1}_{\nu_{1,K}},$ $u^{1}_{\nu_{1,K}+1}]$ of $\gamma$ such that for all $k\in \{1, \ldots, K\}$,
	$${\Card} \big\{ \gamma[u^{1}_{\nu_{1,k}}, u^{1}_{\nu_{1,k}+1}]\cap A_{r_2} \big\}> \lambda_2,$$
	and for all remaining sub-curves $\gamma[u^{1}_{\nu_1}, u^{1}_{\nu_1+1}]$ of $\gamma$, it holds
	$$
	{\Card} \big\{ \gamma[u^{1}_{\nu_1}, u^{1}_{\nu_1+1}]\cap A_{r_2} \big\}\leq \lambda_2.
	$$
	
	For each $\nu_1\in\{0,1, \ldots, \iota_1\}$ and all $u\in A_{r_2}$, by the choice of $w^{r_1}_{\nu_1}$ in Lemma \ref{Mod-again-10-16-1} (1), we have
	\beqq
	\sigma(u^{r_1}_{\nu_1},w^{r_1}_{\nu_1})  \leq  \theta_0 d_{\Omega}(u^{r_1}_{\nu_1})\stackrel{\eqref{17-08-27-01}}{\leq} 2^{r_1+1}\theta_0 T_0
	\leq 2^{r_2}\theta_0 T_0\stackrel{\eqref{17-08-27-01}}{\leq}  \theta_0 d_{\Omega}(u).
	\eeqq
	This implies that the assumptions of Lemmas \ref{Mod-again-03} and \ref{Mod-again-10-16-1} are satisfied, and thus, it follows that for each $k\in \{1, \ldots, K\}$, there exist partitions of $\gamma[u^{1}_{\nu_{1,k}}, u^{1}_{\nu_{1,k}+1}]$ and $\alpha[w^{1}_{\nu_{1,k}}, w^{1}_{\nu_{1,k}+1}]$, respectively, such that all conclusions in Lemmas \ref{Mod-again-03} and \ref{Mod-again-10-16-1} are satisfied.
	In this way, we get a subdivision of $\gamma$ (resp. $\alpha$), which is denoted by
	$$P_{\gamma}^{2}=\left\{u^{2}_{\nu_2}\right\}_{\nu_2=0}^{\iota_2+1}\;\;\quad  \left(\mbox{resp.}\;\;P_{\alpha}^{2}=\left\{w_{\nu_2}^{2}\right\}_{\nu_2=0}^{\iota_2+1}\right),$$
	where $\iota_2$ is a constant with $\iota_{2}\geq \iota_1\geq \lambda_1$, such that for each $\nu_2\in \{0,1, \ldots, \iota_2\}$ and for every $i\leq r_2$, it holds
	$${\Card} \big\{ \gamma[u^{2}_{\nu_2},u^{2}_{\nu_2+1}]\cap A_{i} \big\}\leq \lambda_2.$$
	
	If for each $\nu_2\in\{0,1, \ldots, \iota_2\}$ and for all $i\in \{r_2, \ldots, t_2\}$,
	$${\Card}\big\{\gamma[u^{2}_{\nu_2}, u^{2}_{\nu_2+1}]\cap A_i\big\}\leq \lambda_2,$$
	then for all $i\in \{0,1, \ldots, t_2\}$, we have
	\beq\label{25-12-2-2}{\Card} \big\{ \gamma[u^{2}_{\nu_2},u^{2}_{\nu_2+1}]\cap A_{i} \big\}\leq \lambda_2.\eeq
	
	Let
	$$P_{\gamma}^{2}=\left\{u^{2}_{\nu_1}\right\}_{\nu_2=0}^{\iota_2+1}\;\;\mbox{ and }\;\; P_{\alpha}^{2}=\left\{w^{2}_{\nu_2}\right\}_{\nu_2=0}^{\iota_2+1}.$$
	
	Since $P_{\gamma}^{2}$ and $P_{\alpha}^{2}$ satisfy
	all conclusions in Lemmas \ref{Mod-again-03} and \ref{Mod-again-10-16-1}, and since $\iota_2\geq \lambda_1$, by letting $s_0=\iota_2$, we know from \eqref{25-12-2-2} that $P_{\gamma}^{2}$ and $P_{\alpha}^{2}$ are the desired partitions.
	
	Otherwise, there are some $\nu_2\in\{0,1, \ldots, \iota_2\}$ and an $i\in \{r_2+1, \ldots, t_2\}$ such that
	$${\Card}\big\{\gamma[u^{2}_{\nu_2}, u^{2}_{\nu_2+1}]\cap A_i\big\}\geq \lambda_2.$$
	Set
	$$r_3=\min\big\{i\in\{r_2+1,\ldots,t_2\}:\; q_i>\lambda_2\big\}.$$ Clearly, \beqq\label{25-11-26-1}r_3>r_2.\eeqq
	$\ldots$.
	
	By repeating the above procedures for at most $k_0$ times, where $k_0\leq [\frac{\ell(\alpha)}{T_0}]+1$, we obtain the following partitions of $\gamma$ and $\alpha$:
	$$P_{\gamma}^{{k_0}}=\left\{u^{{k_0}}_{\nu_{k_0}}\right\}_{\nu_{k_0}=0}^{\iota_{k_0}+1}\quad \text{and}\quad P_{\alpha}^{{k_0}}=\left\{w^{{k_0}}_{\nu_{k_0}}\right\}_{\nu_{k_0}=0}^{\iota_{k_0}+1}$$
	such that 	\begin{itemize}
		\item
		$\iota_{k_0}$ is a constant with $\iota_{k_0}\geq \lambda_1$.
		\item
		$P_{\gamma}^{{k_0}}$ and $P_{\alpha}^{{k_0}}$ satisfy all conclusions in Lemmas \ref{Mod-again-03} and \ref{Mod-again-10-16-1}.
		\item for each $\nu_{k_0}\in \{0,1, \ldots, \iota_{k_0}\}$ and any $i\in\{0,1,\ldots,t_2\}$, it holds
		$${\Card} \big\{ \gamma[u^{{k_0}}_{\nu_{k_0}},u^{{k_0}}_{\nu_{k_0}+1}]\cap A_i \big\}\leq \lambda_2.$$
	\end{itemize}
	
	Clearly, $P_{\gamma}^{{k_0}}$ and $P_{\alpha}^{{k_0}}$ are our desired partitions, and hence, the proof of Proposition \ref{2017010-20-1} is complete.
	\epf
	

	\section{Ball separation condition with Gehring-Hayman inequality implies Gromov hyperbolicity}\label{sec-3}
	In this section, we shall prove Theorem \ref{thm:sufficient for Gromov hyperbolic}. To be more precise, suppose that $\Omega$ satisfies both the $C$-Gehring-Hayman inequality and the $C$-ball separation condition with $C\geq 100$. Then we shall prove that $(\Omega, k)$ is $\delta$-Gromov hyperbolic with $\delta=50 C^6(3+C)^2$.
	
	For any  $x_1$, $x_2$ and $x_3$ in $\Omega$, let  $\gamma_{x_1 x_2} \in \Lambda_{x_1 x_2}(\Omega)$,  $\gamma_{x_1 x_3} \in \Lambda_{x_1 x_3}(\Omega)$  and  $\gamma_{x_2 x_3} \in \Lambda_{x_2 x_3}(\Omega)$. Then we shall show the geodesic triangle  $\Delta_{x_1 x_2 x_3}$  has the  $\delta$-thin property.
	
	Fix an arbitrary point $x_0 \in \gamma_{x_1x_2}$, and let $y_0 \in \gamma_{x_1x_3} \cup \gamma_{x_2 x_3}$ be such that
	$$
	\sigma(x_0,y_0)\leq \inf\limits_{y\in\gamma_{x_1x_3} \cup \gamma_{x_2 x_3}}\sigma(x_0,y)+\frac{1}{32}\inf\limits_{x\in\gamma_{x_1x_3} \cup \gamma_{x_2 x_3}}d_{\Omega}(x).
	$$
	Without loss of generality, we may assume that  $y_0 \in \gamma_{x_{1} x_{3}}$.
	
	Since $(\Omega, \sigma)$ satisfies the $C$-ball separation condition, we have $$B_{\sigma}\big(x_0,Cd_{\Omega}(x_0)\big)\cap \big(\gamma_{x_1x_3}\cup \gamma_{x_2x_3}\big)\neq \emptyset,$$
	and so,
	\be\label{eq-20-1}
	\sigma(x_0, y_0) \leq \inf\limits_{y\in\gamma_{x_1x_3} \cup \gamma_{x_2 x_3}}\sigma(x_0,y)+\frac{1}{32}\inf\limits_{x\in\gamma_{x_1x_3} \cup \gamma_{x_2 x_3}}d_{\Omega}(x)\leq C d_{\Omega}(x_0)+\frac{1}{32}d_{\Omega}(y_0).
	\ee
	
	We shall use a contradiction argument to prove
	\beqq
	k_{\Omega}(x_0,\gamma_{x_1 x_3}\cup \gamma_{x_2 x_3})\leq  50 C^6(3+C)^2.
	\eeqq
	To this end, suppose, on the contrary, that
	\be\label{H25-26-0}
	k_{\Omega}(x_0,\gamma_{x_1 x_3}\cup \gamma_{x_2 x_3})> 50 C^6(3+C)^2.
	\ee
	We divide the proof into a few steps. In the first step, we shall prove the following assertion.
	
	\medskip
	\textbf{Step 1.} For each $\gamma_{x_0y_0}\in \Lambda_{x_0y_0}$ and any $x \in \gamma_{x_0 y_0}$, we have
	\be\label{eq-20-2}
	\sigma(y_0,x) \leq 4 C^{3}(1+C) d_{\Omega}(x).
	\ee
	
	If $d_{\Omega}(y_0)\geq 2\sigma(y_0,x)$, then
	$$d_{\Omega}(x)\geq d_{\Omega}(y_0)-\sigma(y_0,x)\geq \sigma(y_0,x),$$ and so, (\ref{eq-20-2}) holds.
	
	Next, we consider the remaining case, that is, $d_{\Omega}(y_0)< 2\sigma(y_0,x)$.
	Since $(\Omega,d)$ satisfies the $C$-Gehring-Hayman inequality, it holds
	\be\label{eq-20-3}
	\sigma(y_0,x) \leq \ell(\gamma_{x_0 y_0}) \leq C \sigma(x_0,y_0).
	\ee
	Let $\alpha=\alpha_x$ be a rectifiable curve connecting $x_0$ and $y_0$ in $\Omega$ such that
	$$\ell(\alpha)\leq \sigma(x_0,y_0)+d_{\Omega}(x).$$
	
	Since $(\Omega,\sigma)$ satisfies the $C$-ball separation condition, we have
	$$B_{\sigma}\big(x,Cd_{\Omega}(x)\big)\cap \alpha\neq \emptyset,$$
	and so, there exists some point  $y \in \alpha$  such that
	\be\label{eq-20-4}
	\sigma(x, y) \leq Cd_{\Omega}(x).
	\ee
	
	Applying the $C$-ball separation condition again, we obtain that there exist two points
	$y_1 \in \gamma_{x_1 x_3}[x_1,y_0] \cup \gamma_{x_1 x_2}[x_1, x_{0}]$ and $y_2\in \gamma_{x_1x_3}[y_0,x_3]\cup \gamma_{x_2 x_3}\cup \gamma_{x_1 x_2}[x_{0}, x_2]$
	such that
	\be\label{H-15-25-1}
	\max \left\{\sigma\left(x, y_1\right), \sigma\left(x, y_{2}\right)\right\} \leq C d_{\Omega}(x).
	\ee
	Note that $\gamma_{x_1 x_3}=\gamma_{x_1 x_3}[x_1,y_0]\cup \gamma_{x_1x_3}[y_0,x_3]$. In below, we present the proof of  \eqref{eq-20-2}  in two cases.
	\medskip
	
	\textbf{Case 1-1.} Either  $y_1\in \gamma_{x_1 x_3}[x_1,y_0]$ or  $y_2\in \gamma_{x_1 x_3}[y_0,x_3] \cup \gamma_{x_2 x_3}$ happens.
	\medskip
	
	In this case, we take $z=y_1$ if $y_1\in \gamma_{x_1 x_3}[x_1,y_0]$, and $z=y_2$ otherwise.
	Then we have
	\beqq
	\begin{aligned}
		\sigma(x_0,y)+\sigma(y,y_0)&\leq \ell(\alpha)\leq \sigma(x_0,y_0)+d_{\Omega}(x) \leq\sigma(x_0,z)+d_{\Omega}(x)+\frac{1}{32}d_{\Omega}(y_0)\\
		&\leq\sigma(x_0,y)+\sigma(y,z)+d_{\Omega}(x)+\frac{1}{32}d_{\Omega}(y_0) \\
		&\leq \sigma(x_0,y)+\sigma(x,y)+\sigma(x,z)+d_{\Omega}(x)+\frac{1}{32}d_{\Omega}(y_0),
	\end{aligned}
	\eeqq
	which implies
	$$
	\sigma(y_0, y) \leq \sigma(x, y)+\sigma(x, z)+d_{\Omega}(x)+\frac{1}{32}d_{\Omega}(y_0) \stackrel{\eqref{eq-20-4}+\eqref{H-15-25-1} }{\leq} 2Cd_{\Omega}(x)+d_{\Omega}(x)+\frac{1}{32}d_{\Omega}(y_0).
	$$
	Consequently, we obtain from the condition $d_{\Omega}(y_0)< 2\sigma(y_0,x)$ that
	$$\sigma(y_0, x) \leq \sigma(y_0, y)+\sigma(x, y) \stackrel{\eqref{eq-20-4} }{\leq} 3 C d_{\Omega}(x)+d_{\Omega}(x)+\frac{1}{32}d_{\Omega}(y_0)\leq 3 C d_{\Omega}(x)+d_{\Omega}(x)+\frac{1}{16}\sigma(y_0,x).$$
	This gives \eqref{eq-20-2}.
	\medskip

	\textbf{Case 1-2.} $y_1 \in \gamma_{x_1 x_2}[x_1, x_{0}]$  and  $y_{2} \in \gamma_{x_1 x_2}[x_{0}, x_2]$.
	\medskip
	
	In this case, note that
	\beqq
	\begin{aligned}
		\ell(\gamma_{x_1x_2}[y_1, y_{2}]) & \geq  \sigma(x_{0}, x)-\sigma(x, y_1)+\sigma(x_{0}, x)-\sigma(x, y_{2}) \\
		& \stackrel{\eqref{H-15-25-1}}{\geq} 2 \sigma(x, x_{0})-2 C d_{\Omega}(x).
	\end{aligned}
	\eeqq
	
	Since  $(\Omega,d)$ satisfies the $C$-Gehring-Hayman inequality, we further have
	$$\ell(\gamma_{x_1x_2}[y_1, y_{2}]) \leq C \sigma(y_1,y_{2})\leq C\left(\sigma(x, y_1)+\sigma(x, y_{2})\right)\stackrel{\eqref{H-15-25-1}}{\leq} 2C^2d_{\Omega}(x).$$
	Combining the above two estimates gives
	\be\label{eq-20-5}
	\sigma(x_{0}, x) \leq C(1+C)d_{\Omega}(x).
	\ee
	
	If  $\sigma(x_{0}, x) \leq \frac{1}{2} d_{\Omega}(x_{0})$, then  $d_{\Omega}(x) \geq d_{\Omega}(x_{0})-\sigma(x_{0}, x) \geq \frac{1}{2} d_{\Omega}(x_{0})$, and thus, we have
	$$
	\sigma(y_0, x) \stackrel{\eqref{eq-20-3}}{\leq} C \sigma(x_{0}, y_0) \stackrel{\eqref{eq-20-1}}{\leq} 2C^{2} d_{\Omega}(x_{0}) \leq 4 C^{2} d_{\Omega}(x).
	$$

	If  $\sigma(x_{0}, x)>\frac{1}{2} d_{\Omega}(x_{0})$, then we have
	$$
	\sigma(y_0, x) \stackrel{\eqref{eq-20-3}}{\leq} C \sigma(x_{0}, y_0) \stackrel{\eqref{eq-20-1}}{\leq} 2C^{2} d_{\Omega}(x_{0}) \leq 4 C^{2} \sigma(x_{0}, x) \stackrel{\eqref{eq-20-5}}{\leq} 4 C^{3}(1+C) d_{\Omega}(x).
	$$
	
	{In either case,} \eqref{eq-20-2} is proved.
	
	In the second step, we shall prove the following estimate.
	\medskip
	
	\textbf{Step 2.} The following estimate holds:
	\be\label{eq-20-star}
	d_{\Omega}(y_0)\leq\frac{1}{16(3+C)^{4}}\sigma(x_0,y_0).
	\ee
	
	\medskip
	
	Suppose, on the contrary, that \eqref{eq-20-star} fails, that is, $$d_{\Omega}(y_0)>\frac{1}{16(3+C)^{4}}\sigma(x_0,y_0).$$ Then we may obtain a contraction to \eqref{H25-26-0} as follows. Since $(\Omega,d)$ satisfies the $C$-Gehring-Hayman inequality, we have for each $x\in\gamma_{x_0y_0}$,
	\beqq
	\ell(\gamma_{x_0 y_0}[y_0, x]) \leq C \sigma(y_0, x) \stackrel{\eqref{eq-20-2}}{\leq} 4 C^4(1+C) d_{\Omega}(x).
	\eeqq
	This gives, in particular, that for each  $w \in \gamma_{x_0 y_0}$ and $x\in\gamma_{x_0 y_0}[w,x_0]$, it leads to
	$$\ell(\gamma_{x_0 y_0}[w, x]) \leq \ell(\gamma_{x_0 y_0}[y_0, x]) \leq 4 C^4(1+C) d_{\Omega}(x).$$
	Then we know from Lemma \ref{lem-3-1} that for each  $w \in \gamma_{x_0 y_0}$, it holds
	\be\label{eq-20-7}
	k_{\Omega}(w, x_0) \leq 16C^4(1+C) \log \left(1+\frac{\ell(\gamma_{x_0y_0})}{d_{\Omega}(w)}\right) \stackrel{\eqref{eq-20-3}}{\leq} 16C^4(1+C) \log \left(1+\frac{C\sigma(x_0,y_0)}{d_{\Omega}(w)}\right).
	\ee
	
	Since $d_{\Omega}(y_0)>\frac{1}{16(3+C)^{4}}\sigma(x_0,y_0)$, we may apply \eqref{eq-20-7} with $w=y_0$ to obtain
	$$k_{\Omega}(x_0,y_0)<80C^4(3+C) \log 16(3+C)<50 C^6(3+C)^2,$$
	which clearly contradicts with \eqref{H25-26-0}. This proves \eqref{eq-20-star}.
	
	The following claim is very useful for the discussions in the next step.
	
	\bcl\label{claim:find z1 and z2}
	There exist two points  $z_{1} \in \gamma_{x_{1} x_{2}}[x_{1}, x_0]$  and  $z_{2} \in \gamma_{x_{1}x_{2}}[x_0, x_{2}]$  such that
	\be\label{eq-20-8}
	\sigma(y_0, z_{1}) \leq \frac{3}{5(3+C)^{2}} \sigma(x_0, y_0)\ \text{ and }\  \sigma(y_0, z_{2}) \leq \frac{3}{5 (3+C)^2} \sigma(x_0, y_0).
	\ee
	\ecl
	
	We shall only prove the existence of $z_1$ in the claim, as the proof of the other case is similar. Let $w_1 \in \gamma_{x_0 y_0}$ be such that
	\be\label{eq-20-9}
	\sigma(y_0, w_1)=\frac{1}{8 (3+C)^{ 4}} \sigma(x_0, y_0).
	\ee  
	
	If $\sigma(y_0, x_{1}) \leq 2 \sigma(y_0, w_1)$, then for $z_{1}=x_{1}$, \eqref{eq-20-9} gives
	$$\sigma(y_0,z_1)=\sigma(y_0, x_{1})\leq \frac{1}{4(3+C)^4}\sigma(x_0,y_0) \leq \frac{3}{5(3+C)^{2}} \sigma(x_0, y_0).$$
	In this case, the proof of claim is complete.
	
	For the remaining case, that is, $\sigma(y_0, x_{1}) > 2 \sigma(y_0, w_1)$, let
	$u_1\in\gamma_{x_{1}x_{3}}[y_0,x_{1}]$ be the last point along the direction from $y_0$ to $x_{1}$ such that
	\be\label{eq-20-10}
	\sigma(y_0,u_1)=2\sigma(y_0,w_1)=\frac{1}{4(3+C)^4}\sigma(x_0,y_0).
	\ee
	Then the $C$-Gehring-Hayman inequality ensures that
	for each  $u \in \gamma_{u_1 w_1}$,
	\be\label{eq-20-11}
	\sigma(u_1, u) \leq \ell(\gamma_{u_1 w_1}) \leq C\sigma(u_1, w_1) \leq C\left(\sigma(u_1, y_0)+\sigma(y_0, w_1)\right)\stackrel{\eqref{eq-20-9}+ \eqref{eq-20-10}}{\leq}\frac{3\sigma(x_0,y_0)}{8(3+C)^3},
	\ee
	and the $C$-ball separation condition guarantees that for each $u \in \gamma_{u_{1} w_1}$, there exists some point $v \in \gamma_{x_1x_2}[x_1,x_0]\cup \gamma_{x_1 x_3}[x_1, u_{1}] \cup \gamma_{y_0 x_0}[w_1, x_0]$  such that
	\be\label{eq-20-14}
	\sigma(u, v) \leq C d_{\Omega}(u).
	\ee
	
	\textbf{Case 2-1.} $v \in \gamma_{x_1x_2}[x_1,x_0]$.
	\medskip
	
	In this case, note first that
	\beqq
	\begin{aligned}
		\sigma(y_0, v) &\leq \sigma(y_0, u_1)+\sigma(u_1, u)+\sigma(u, v)\stackrel{\eqref{eq-20-14}}{\leq} \sigma(y_0, u_1)+\sigma(u_1, u)+C d_{\Omega}(u)\\
		&\leq \sigma(y_0, u_1)+\sigma(u_1, u)+C (d_{\Omega}(u_1)+\sigma(u,u_1))  \\
		&\leq \sigma(y_0, u_1)+(1+C)\sigma(u_1,u)+C(d_{\Omega}(y_0)+\sigma(y_0,u_1))\\
		&\stackrel{\eqref{eq-20-10}}{=}\frac{1+C}{4(3+C)^{4}} \sigma(x_0, y_0)+(1+C)\sigma(u_1,u)+Cd_{\Omega}(y_0).
	\end{aligned}
	\eeqq
	Then  \eqref{eq-20-star} and \eqref{eq-20-11} ensure that
	$$\sigma(y_0, v)\leq \frac{6C^2+29C+22}{16(3+C)^4}\sigma(x_0, y_0)<\frac{3}{5(3+C)^{2}} \sigma(x_0, y_0),$$ from which
	the claim follows by taking $v=z_1$ above.
	\medskip
	
	\textbf{Case 2-2.} $v \in \gamma_{x_1 x_3}[x_1, u_{1}] \cup \gamma_{y_0 x_0}[w_1, x_0]$.
	\medskip
	
	In this case, we claim that the following estimate holds:
	\be\label{eq-20-15}
	\sigma\left(u_{1}, u\right)<12 C^4(3+C)^3d_{\Omega}(u).
	\ee
	
	\textbf{Case 2-2-1.} $v\in \gamma_{x_0 y_0}[w_1, x_0]$.
	\smallskip
	
	In this case, by the $C$-Gehring-Hayman inequality, we have
	\beqq
	\begin{aligned}
		\sigma(y_0, v) \geq  C^{-1} \ell(\gamma_{y_0 x_0}[y_0, v]) \geq C^{-1} \sigma(y_0, w_1)\stackrel{\eqref{eq-20-9}}{>} \frac{1}{8C(3+C)^{4}} \sigma(x_0, y_0),
	\end{aligned}
	\eeqq
	and so,
	\beqq
	\begin{aligned}
		\sigma(u_1, u) & \stackrel{\eqref{eq-20-11}}{\leq} \frac{3}{8(3+C)^3} \sigma(x_0, y_0) \leq 3C(3+C) \sigma(y_0, v)  \stackrel{\eqref{eq-20-2}}{<} 12C^4(3+C)^2 d_{\Omega}(v) \\
		&\leq 12C^4(3+C)^2 \left(d_{\Omega}(u)+\sigma(u, v)\right) \stackrel{\eqref{eq-20-14}}{\leq} 12C^4(3+C)^3  d_{\Omega}(u).
	\end{aligned}
	\eeqq
	This proves \eqref{eq-20-15}.
	
	\medskip
	
	\textbf{Case 2-2-2.} $v\in \gamma_{x_1 x_3}[x_1, u_{1}]$.
	\smallskip
	
	In this case, let $\alpha=\alpha_u$ be a rectifiable curve connecting $u_1$ and $y_0$ in $\Omega$ such that
	$\ell(\alpha)\leq \sigma(u_1,y_0)+d_{\Omega}(u)$. Then by the $C$-ball separation condition, there exists some point
	$v_1\in\gamma_{y_0 x_0}[w_1, y_0]\cup \alpha$ such that
	\be\label{H25-26-4}
	\sigma(u, v_1) \leq Cd_{\Omega}(u).
	\ee
	
	Moreover, the choice of $u_1$ implies that $$\sigma(y_0,u_1)\leq \sigma(y_0,v).$$
	Then we have
	\beq\label{H25-26-7}
	\begin{aligned}
		\frac{1}{4 (3+C)^4}\sigma(x_{0}, y_{0})&\stackrel{\eqref{eq-20-10}}{=} \sigma(y_{0}, u_1)\leq \sigma(y_{0}, v)\leq \sigma(y_{0}, u)+\sigma(u,v)
		\\  &\stackrel{\eqref{eq-20-14}}{\leq} \sigma(y_{0}, u)+Cd_{\Omega}(u).
	\end{aligned}
	\eeq
	
	We first consider the case $v_1 \in \gamma_{x_{0} y_{0}}[y_{0}, w_{1}]$.
	\begin{itemize}
		\item  If $\sigma(y_{0}, v_1)\geq \frac{1}{8 (3+C)^{5}} \sigma(x_{0}, y_{0})$, then we may argue as in Case 2-2-1 to derive \eqref{eq-20-15}.
		
		\item If $\sigma(y_{0}, v_1)<\frac{1}{8(3+C)^{5}} \sigma(x_{0}, y_{0})$, then
		\beqq
		\sigma(y_{0}, u) \leq \sigma(y_{0}, v_1)+\sigma(u, v_1)\stackrel{\eqref{H25-26-4}}{<}\frac{1}{8 (3+C)^{5}}\sigma(x_{0}, y_{0})+Cd_{\Omega}(u),
		\eeqq
		which, together with (\ref{eq-20-11}) and (\ref{H25-26-7}), shows that
		$$\sigma(u_1,u)<6C(3+C)d_{\Omega}(u),$$ which implies \eqref{eq-20-15}.
	\end{itemize}

	Next, we consider the case $v_1\in\alpha$. In this case, we have
	\beqq
	\begin{aligned}
		\sigma(y_{0}, u) \nonumber&\leq \sigma(y_{0}, v_1)+\sigma(u, v_1)\stackrel{\eqref{H25-26-4}}{\leq} \ell(\alpha)-\ell(\alpha[u_1,v_1])+Cd_{\Omega}(u)
		\\ &\leq \sigma(y_0,u_1)+d_{\Omega}(u)-\sigma(u_1,u)+\sigma(u,v_1)+Cd_{\Omega}(u)
		\\ &\stackrel{\eqref{H25-26-4}}{\leq} \sigma(y_0,u_1)+d_{\Omega}(u)-\sigma(u_1,u)+2Cd_{\Omega}(u)\\
		&\stackrel{\eqref{H25-26-7}}{\leq} \sigma(y_{0}, u)+d_{\Omega}(u)-\sigma(u_1,u)+3Cd_{\Omega}(u).
	\end{aligned}
	\eeqq
	Thus it follows that
	$$\sigma(u_1,u)\leq 3Cd_{\Omega}(u)+d_{\Omega}(u)=(1+3C)d_{\Omega}(u).$$
	The above estimate gives \eqref{eq-20-15}.\medskip
	
	Let us continue the proof based on \eqref{eq-20-15}.
	Note that \eqref{eq-20-7} and the $C$-Gehring-Hayman inequality imply
	$$k_{\Omega}(w_1,x_0)\leq 8C^4(1+C)\log\Big(1+\frac{\ell(\gamma_{x_0y_0}[w_1,x_0])}{d_{\Omega}(w_1)})\leq 8C^4(1+C)\log\Big(1+\frac{C\sigma(y_0,x_0)}{d_{\Omega}(w_1)}\Big).$$
	By \eqref{eq-20-2} and \eqref{eq-20-9}, we have
	\[
	d_{\Omega}(w_1)\geq \frac{\sigma(y_0,w_1)}{4C^3(1+C)}=\frac{\sigma(x_0,y_0)}{32C^3(3+C)^4}.
	\]
	The above two estimates yield
	$$k_{\Omega}(w_1,x_0)\leq 150C^4(1+C)\log 2(3+C),$$
	which, together with \eqref{H25-26-0}, guarantees
	\be\label{H25-26-9}
	k_{\Omega}(u_1,w_1)\geq k_{\Omega}(x_0,u_1)-k_{\Omega}(w_1,x_0)> 49C^6(3+C)^2.
	\ee
	
	Moreover, by Lemma \ref{lem-3-1}, \eqref{eq-20-11}  and \eqref{eq-20-15}, we have
	$$
	k_{\Omega}(u_{1}, w_1) \leq 48 C^4(3+C)^3\log \left(1+\frac{\ell(\gamma_{u_1w_1})}{d_{\Omega}(u_{1})}\right) \leq 48 C^4 (3+C)^3\log\Big(1+\frac{3\sigma(x_0,y_0)}{8(3+C)^3d_{\Omega}(u_{1})}\Big).
	$$
	This, combining with \eqref{H25-26-9}, ensures that
	\be\label{eq-20-16}
	d_{\Omega}(u_{1})<\frac{1}{48C^3(3+C)^{6}} \sigma(x_0, y_0).
	\ee
	
	Since $(\Omega,\sigma)$ satisfies the $C$-ball separation condition, there exists some point $w \in \gamma_{x_1x_2}[x_1, x_0] \cup \gamma_{y_0 x_0}$  such that
	$$
	\sigma(u_1, w) \leq C d_{\Omega}(u_1),
	$$
	and so, by \eqref{eq-20-16}, $$\sigma(u_{1}, w) \leq \frac{1}{48C^2(3+C)^{6}} \sigma(x_0, y_0).$$
	Then we have
	$$
	\sigma(y_0, w) \geq \sigma(u_{1}, y_0)-\sigma(u_{1}, w)\stackrel{\eqref{eq-20-10} }{>}\frac{1}{5(3+C)^{4}} \sigma(x_0, y_0).
	$$
	
	Suppose that $w\in \gamma_{y_0 x_0}$. Then the above estimate gives
	$$\frac{1}{20C^3(3+C)^5} \sigma(x_0, y_0)<\frac{1}{4C^3(1+C)}\sigma(y_0, w)\stackrel{ \eqref{eq-20-2}}{\leq} d_{\Omega}(w)\leq d_{\Omega}(u_1)+\sigma(u_1, w)\leq (1+C)d_{\Omega}(u_1),$$
	which clearly contradicts with \eqref{eq-20-16}. This ensures that $w \in \gamma_{x_1x_2}[x_1, x_0] $.
	
	Since the similar argument used in Case 2-1 ensures that $\sigma(y_0,w)\leq \frac{3}{5(3+C)^2}\sigma(x_0,y_0)$, the first inequality in Claim \ref{claim:find z1 and z2} follows from letting $w=z_1$. This completes the proof of the existence of $z_1$, and hence, the claim is proved.
	
	\medskip
	
	\textbf{Step 3.} Contradiction.
	\medskip

	By \eqref{eq-20-8}, we have
	\begin{align*}
		\ell(\gamma_{x_1 x_2}[z_1, z_2])& \geq \sigma(z_1, x_0)+\sigma(z_2, x_0) \\
		& \geq \sigma(y_0, x_0)-\sigma(y_0, z_1)+\sigma(y_0, x_0)-\sigma(y_0, z_2) \\
		& \geq\left(2-\frac{6}{5(3+C)^{2}}\right) \sigma(y_0, x_0).
	\end{align*}
	
	On the other hand, the $C$-Gehring-Hayman inequality gives
	\begin{align*}
		\ell(\gamma_{x_1 x_2}[z_1, z_2])  \leq C \sigma(z_1, z_2) \leq C(\sigma(z_1, y_0)+\sigma(y_0, z_2))  \leq \frac{6C}{5(3+C)^2} \sigma(y_0, x_0),
	\end{align*}
	which clearly contradicts with the previous estimate. This indicates that \eqref{H25-26-0} can not hold, and thus, the proof of Theorem \ref{thm:sufficient for Gromov hyperbolic} is complete.

	\section{Geometric applications}\label{sec:6}

	In this section, we give two geometric applications of our obtained results. First of all, as an application of Theorem \ref{positive-answer}, we show that the ball separation condition, together with a geometric {\it LLC-2} condition, completely characterizes the inner uniformity. Recall that a domain $\Omega$ in a metric space $X=(X,d)$ is $c_0$-{\it LLC-2} if each pair of points $a,b\in \Omega\backslash \overline{\mathbb{B}(x,r)}$ can be joined in $\Omega\backslash \overline{\mathbb{B}(x,r/c_0)}$, where $c_0\geq 1$.
	\begin{thm}[Geometric characterization of inner uniformity]\label{thm:ball separation+LLC2 implies GH metric space}
		Let $Q>1$, and let $X=(X,d)$ be a $Q$-doubling length space and $\Omega\subset X$ a proper subdomain such that $(\Omega,k)$ is geodesic. Then
		the following conclusions hold:
		\begin{enumerate}
			\item 	If $\Omega$ is $c_0$-LLC-2 with the $C$-ball separation condition, then it is {$C_1$}-inner uniform with $C_1=C_1(C,Q,C_0)$.
			
			\item If $\Omega$ is $C$-inner uniform, then it is {$C_2$}-LLC-2 and satisfies the $C_2$-ball separation condition with $C_2=C_2(C)$.
		\end{enumerate}
	\end{thm}
	
	We remark that it follows from the proof of Theorem \ref{thm:ball separation+LLC2 implies GH metric space} below that we can take $C_1=\max\{(2Q)^{4(36C\cdot Q^{5C})^{(8QC)^8}},({16}ec_0C)^2([Q^{2\log_2 70ec_0C}]+1)\}$ and $C_2=130C^4e^{32C^4}$. Also,
	Theorem \ref{thm:ball separation+LLC2 implies GH metric space} is almost sharp in the sense that the $Q$-doubling assumption cannot be removed, as there is an {\it LLC-2} Gromov hyperbolic domain $\Omega$ in an infinite dimensional Banach space, which fails to be inner uniform; see \cite[Remark 3.16]{Vai2004}.
	
	As an immediate consequence of Theorem \ref{thm:ball separation+LLC2 implies GH metric space}, we obtain that in a locally compact $Q$-doubling length space $X$, a proper subdomain $\Omega\subset X$ is inner uniform if and only if it is John and satisfies the ball separation condition, with explicit dependence on the relevant coefficients. {This shows that an inner uniform domain differs from a John domain exactly a geometric ball separation condition. }
	
	Note that under the assumption of Theorem \ref{thm:ball separation+LLC2 implies GH metric space}$(1)$, by Theorem \ref{positive-answer}, $\Omega$ satisfies the $c_1$-Gehring-Hayman inequality with $c_1=(2Q)^{4(36C\cdot Q^{5C})^{(8QC)^8}}$. Thus it is left to verify the double cone condition, i.e. Definition \ref{def:uniform domain} (1). For this, we have the following stronger result, which shows that each quasihyperbolic geodesic satisfies the double cone condition.
	
	\blem\label{Theorem-1.8}
	For $x_1, x_2\in \Omega$, let $\gamma=\gamma_{x_1x_2}\in\Lambda_{x_1x_2}(\Omega)$. Then for each $x\in \gamma$, we have
	$$\min\left\{\ell(\gamma[x_1,x]),\ell(\gamma[x_2,x])\right\}\leq c_2 d_{\Omega}(x),$$
	where $c_2=({16}ec_0C)^2([Q^{2\log_2 70ec_0C}]+1)$.
	\elem

	\bpf 
	Let $x_0 \in \gamma$ be such that
	\beqq
	\min \left\{\diam( \gamma[x_1, x_0]), \diam( \gamma[x_2, x_0])\right\} \geq \frac{1}{3} \diam(\gamma).
	\eeqq

	We first prove that for each $i \in \{1,2\}$ and $x\in  \gamma[x_i, x_0]$,
	\be\label{H25-0729-1}
	\diam( \gamma[x_i ,x])\leq 10c_0 Cd_{\Omega}(x).
	\ee
	Without loss of generality, we only consider the case $x\in \gamma[x_1,x_0]$
	and shall prove that 
	\beq\label{333}
	\diam(\gamma[x_1, x]) \leq 10c_0 Cd_{\Omega}(x).
	\eeq

	Let $z_1\in \gamma[x_1,x]\cap \left(\Omega \backslash\mathbb{B}\left(x,\frac{1}{9}\diam( \gamma[x_1, x])\right)\right)$. It follows from $\diam(\gamma[x_2,x])\geq \frac{1}{3}\diam(\gamma)$ that  there exists $z_2 \in \gamma[x_2,x]$ that such that $z_2 \in \Omega \backslash\mathbb{B}\left(x,\frac{1}{9}\diam( \gamma[x_1, x])\right)$.
	
	Since $\Omega$ is $c_0$-{\it LLC-2}, there exists some rectifiable curve
	$$\beta \subset \Omega\backslash \mathbb{B}\left(x,\frac{1}{9c_0}\diam( \gamma[x_1, x])\right)$$ joining $z_1$ and $z_2$,
	and since $(\Omega,\sigma)$ satisfies the $C$-ball separation condition, we have
	\beqq
	\mathbb{B}_{\sigma}(x,Cd_{\Omega}(x)) \cap \beta \neq \emptyset,
	\eeqq
	from which it follows that
	\beqq
	\diam(\gamma[x_1, x]) \leq 9c_0 Cd_{\Omega}(x).
	\eeqq
	This establishes \eqref{333}, and thus, completes the proof of \eqref{H25-0729-1}.

	Next, our aim is to prove that for each $i=1,2$ and $x \in \gamma[x_i, x_0]$, it holds
	\beq\label{444}
	\ell (\gamma[x_i, x])\leq (16ec_0C)^2([Q^{2\log_2 70ec_0C}]+1) d_{\Omega}(x).
	\eeq
	
	Without loss of generality, by (\ref{H25-0729-1}), we only consider the case $x\in \gamma[x_1,x_0]$. Let $y_1 \in \gamma[x_1, x]$ be such that $$d_{\Omega}(y_1)\geq \frac{1}{2}\sup_{y \in \gamma[x_1, x]}\{d_{\Omega}(y)\}.$$
	Then we have
	\beq\label{555}
	d_{\Omega}(y_1) \leq d_{\Omega}(x)+ d(y_1 ,x)\leq d_{\Omega}(x)+\diam(\gamma[x_1, x]) \stackrel{\eqref{333}}{\leq} (1+10c_0 C)d_{\Omega}(x).
	\eeq
	Note that by \eqref{333}, for each $z\in \gamma[x_1 ,x])$, it holds
	\beqq
	d(x_1,z) \leq \diam( \gamma[x_1, z]) \leq 10c_0 C d_{\Omega}(z).
	\eeqq
	Thus we infer from Lemma \ref{qs-8} (with $\beta=\gamma[x_1, x]$, $\mu_1=10c_0C$ and $x_0=y_1$) that
	\beqq
	\ell(\gamma[x_1 ,x]) \leq {25}e^2c_0C([Q^{2\log_2 70ec_0C}]+1)d_{\Omega}(y_1) \stackrel{\eqref{555}}{\leq} ({16}ec_0C)^2([Q^{2\log_2 70ec_0C}]+1)d_{\Omega}(x).
	\eeqq
	This establishes \eqref{444}. The proof of Lemma \ref{Theorem-1.8} is thus complete.
	\epf
	\medskip

	\begin{proof}[Proof of Theorem \ref{thm:ball separation+LLC2 implies GH metric space}]
		
		(1) This follows from Lemma \ref{Theorem-1.8} and Theorem \ref{positive-answer}.
		One can take
		$$C_1=\max\{c_1,c_2\}=\max\{(2Q)^{4(36C\cdot Q^{5C})^{(8QC)^8}},({16}ec_0C)^2([Q^{2\log_2 70ec_0C}]+1)\}.$$

		(2)  Since $\Omega$ is $C$-inner uniform, it is $(2C+1)$-{\it LLC-2} by \cite[Lemma 3.7]{Vai7}. Furthermore, it follows from \cite[Lemma 2.35]{Vai7} (or the proof of \cite[Lemma 2.13]{BHK}) that every quasihyperbolic geodesic $\gamma$ in a $C$-inner uniform domain $\Omega$ is a {$C_2$}-inner uniform curve with $C_2=128C^4e^{32C^4}$,
		and thus, the domain $\Omega$ satisfies the $C_2$-ball separation condition.
		
	\end{proof}
	
	Secondly, as a direct application of Theorem \ref{thm:new main metric space}, we obtain the following corollary, which provides an affirmative answer to an open question \cite[Question 1.9]{ZP-2024-Pisa} in the setting of locally compact $Q$-doubling length spaces. This result is also new even if $X=\IR^n$. Recall that a proper subdomain $\Omega$ in a locally compact length space $(X,d)$ is called {\it quasihyperbolic $c$-John} if there is $c\geq1$ such that every quasihyperbolic geodesic $\gamma\subset \Omega$ is a double $c$-cone curve. 
	\begin{thm}\label{thm:Zhou problem}
		Let $X$ be a locally compact $Q$-doubling length space and $\Omega\subset X$ a proper length subdomain. If $\Omega$ is quasihyperbolic $c$-John, then $(\Omega,k)$ is $\delta$-Gromov hyperbolic, quantitatively.
	\end{thm}
	\begin{proof}
		Since $\Omega$ is length, the inner metric $\sigma$ coincides with $d$.
		As $\Omega$ is quasihyperbolic $c$-John, by definition, $(\Omega,\sigma)$ satisfies the $c$-ball separation condition,
		and thus, it follows from Theorem \ref{thm:new main metric space} that $(\Omega,k)$ is $\delta$-Gromov hyperbolic with $\delta=\delta(c,Q)$.
	\end{proof}
	
	\appendix
	
	\section{A few technical lemmas}\label{appendix-A}
	In this appendix, we collect some technical results that are needed for the discussions in Section \ref{sec-4-3}.
	\blem\label{lem-23-4.2}
	Suppose that  $\gamma_{x z_{1,1}y_{1}} \in P_{\alpha}^{\gamma_{xz_{1}y}}(3C)$, $y_2\in\alpha[y_1,y]$ and $\gamma_{y_1z_{2,1}y_2} \in O_{\alpha[y_1,y_2]}^{\gamma_{x z_1 y}}(2C)$. For  $y_{3} \in \alpha[y_{1}, y_{2}]$, fix $\gamma_{y_1y_3}\in \Lambda_{y_1y_3}(\Omega)$ and $\gamma_{xy_3}\in \Lambda_{xy_3}(\Omega)$. If $k_{\Omega}(z_{2,1}, z_{3})>24C$ and  $2C\leq k_{\Omega}(z_{3}, \gamma_{y_{1}y_{3}}) \leq 7C$ for some $z_{3} \in \gamma_{y_{1} y_{2}}[y_{1}, z_{2,1}]$, then $2C \leq k_{\Omega}(z_{1}, \gamma_{x y_{3}})\leq {3}C$.
	\elem
	\bpf
To get the upper bound for the quantity $k_{\Omega}(z_{1}, \gamma_{x y_{3}})$ in the lemma, let $z_{5,3} \in \gamma_{x y_{3}}$ be such that $k_{\Omega}(z_{1}, z_{5,3})=k_{\Omega}(z_{1}, \gamma_{x y_{3}})$. Since $\gamma_{x z_{1,1}y_{1}} \in P_{\alpha}^{\gamma_{xz_{1}y}}(3C)$, Definition \ref{def:class P}(3) implies that 
	\be\label{H25-0505-2}
	k_{\Omega}(z_{1}, \gamma_{xy_3})=k_{\Omega}(z_{1}, z_{5,3})\leq 3C.
	\ee
	
In the following, we show the lower bound for $k_{\Omega}(z_{1}, \gamma_{x y_{3}})$ in the lemma, that is,
\be\label{25-12-15-1}
k_{\Omega}(z_{1}, \gamma_{x y_{3}})\geq 2C.
\ee

We start the proof with some preparation.	First of all, we prove
	\be\label{eq-23-1-1}
	k_{\Omega}(z_{1}, \gamma_{y_{1}y_{3}}) >21C.
	\ee
To show such a lower bound for $k_{\Omega}(z_{1}, \gamma_{y_{1}y_{3}})$,	
let $z_{1,3} \in \gamma_{y_{1} y_{3}}$ be such that
	\be\label{H25-27-1}
	k_{\Omega}(z_{2,1}, \gamma_{y_{1}y_{3}})=k_{\Omega}({z_{2,1}}, z_{1,3}),
	\ee
	and then, fix $\gamma_{z_{2,1}z_{1,3}}\in\Lambda_{z_{2,1}z_{1,3}}(\Omega)$.  Since  $k_{\Omega}(z_3,\gamma_{y_1y_3}[y_1,z_{1,3}])\geq k_{\Omega}(z_{3}, \gamma_{y_{1}y_{3}})$, it follows from the assumption $2C\leq k_{\Omega}(z_{3}, \gamma_{y_{1}y_{3}}) \leq 7C$ in the lemma that $k_{\Omega}(z_3,\gamma_{y_1y_3}[y_1,z_{1,3}])\geq 2 C$. Then the assumptions of $(\Omega,k)$ being $C$-Gromov hyperbolic ensures that there exists some point $z_{3,1} \in \gamma_{z_{2,1}z_{1,3}}$ such that
	$$k_{\Omega}(z_{3}, z_{3,1}) \leq C,$$
	which, together with the assumption $k_{\Omega}(z_{2,1}, z_{3})>24C$, shows that
	\begin{equation*}
		k_{\Omega}(z_{2,1}, z_{1,3})\geq k_{\Omega}(z_{2,1}, z_{3,1}) \geq k_{\Omega}(z_{2,1}, z_{3})-k_{\Omega}(z_{3}, z_{3,1}) >23 C.
	\end{equation*}
	
Since $\gamma_{y_1z_{2,1}y_2} \in O_{\alpha[y_1,y_2]}^{\gamma_{x z_1 y}}(2C)$, we know
\be\label{eq-23-2}
k_{\Omega}(z_{1},z_{2,1}) \leq 2C.
\ee	
This, together with \eqref{H25-27-1}, leads to
	$$
	k_{\Omega}(z_{1}, \gamma_{y_{1}y_{3}}) \geq k_{\Omega}(z_{2,1}, \gamma_{y_{1}y_{3}})-k_{\Omega}(z_{1}, {z_{2,1}}) >21C,
	$$
	which gives \eqref{eq-23-1-1}.
	

	
	Next, we find a point $w_1\in\gamma_{y_1y_2}[z_3,z_{2,1}]$ such that
	\be\label{25H-0724-7}
	2C<k_{\Omega}(w_1,\gamma_{y_1y_3})\leq 7C.
	\ee
	
Observe that, in the lemma, it is assumed that $2C\leq k_{\Omega}(z_{3}, \gamma_{y_{1}y_{3}}) \leq 7C$. Based on this observation, to prove the existence of the point $w_1$, we consider two possibilities.
	If $k_{\Omega}(z_3,\gamma_{y_1y_3})>2C$, then $w_1=z_3$ satisfies \eqref{25H-0724-7}.
	For the other possibility $k_{\Omega}(z_3,\gamma_{y_1y_3})=2C$, we take $w_{1,1}\in \gamma_{y_1y_2}[z_3,z_{2,1}]$ be such that \be\label{25H-0724-8}
	k_{\Omega}(z_3,w_{1,1})=5C.
	\ee
	
	We claim
	\be\label{25H-0724-9}k_{\Omega}(w_{1,1},\gamma_{y_1y_3})\geq 3C.\ee
	Otherwise, there exists some point $w_{1,2}\in\gamma_{y_1y_3}$ such that
	$$k_{\Omega}(w_{1,1},w_{1,2})<3C.$$
	Then we infer from \eqref{25H-0724-8} that $$k_{\Omega}(z_3,w_{1,1})\geq 2C+k_{\Omega}(w_{1,1},w_{1,2}).$$
This implies that Lemma \ref{lem-2-3.0} is applicable to  the points $y_1, w_{1,1}, w_{1,2}$ and $z_3$,  
 and then, it follows that $k_{\Omega}(z_3,\gamma_{y_1y_3})\leq C$. This contradicts with the assumption $2C\leq k_{\Omega}(z_{3}, \gamma_{y_{1}y_{3}}) \leq 7C$ in the lemma. This contradiction proves (\ref{25H-0724-9}).
	
	By (\ref{25H-0724-9}), we may select $w_1\in\gamma_{y_1y_2}[z_3,w_{1,1}]$ such that $k_{\Omega}(w_1,\gamma_{y_1y_3})>2C$. 
Obviously,  (\ref{25H-0724-8}) implies that 
\be\label{25-12-16-1}
k_{\Omega}(z_3,w_1)\leq 5C.
\ee 	
	Then the assumption $k_{\Omega}(z_3,\gamma_{y_1y_3})=2C$ in this possibility guarantees that such a point $w_1$ is what we need. 

The preparation still needs the following assertion: There is some point $z_{2,3} \in \gamma_{x y_{1}}[z_{1,1}, y_{1}]$ such that
\be\label{H25-26-12}
k_{\Omega}(w_1, z_{2,3}) \leq C,
\ee	 
and for any $z\in\gamma_{xy_1}[z_{2,3},y_1]$,
\be\label{H25-0505-1}
k_{\Omega}(z_{1}, z)>6C.
\ee

First, we find a point in $\gamma_{x y_{1}}[z_{1,1}, y_{1}]$ by applying Lemma \ref{lem-2-3.0}.
	Since  $\gamma_{x z_{1,1}y_{1}} \in P_{\alpha}^{\gamma_{xz_{1}y}}(3C)$, Lemma \ref{Lemma4-1.1} shows
\be\label{H25-26-10}
k_{\Omega}(z_{1},z_{1,1})<5C,\ee	
and thus, it follows from \eqref{eq-23-2} that
\be\label{H25-26-11}
k_{\Omega}(z_{1,1},z_{2,1})\leq k_{\Omega}(z_1,z_{1,1})+k_{\Omega}(z_1,z_{2,1})\leq 7C.
\ee
Then the assumption $k_{\Omega}(z_{2,1}, z_{3})>24C$ in the lemma and \eqref{25-12-16-1} ensure that
\be\label{25H-0724-10} k_{\Omega}(w_1,z_{2,1})\geq k_{\Omega}(z_{2,1}, z_{3})-  k_{\Omega}(z_3,w_1) \geq 19C\geq 12C+k_{\Omega}(z_{1,1},z_{2,1}).\ee
This illustrates that Lemma \ref{lem-2-3.0} is applicable to
 the points
$y_1, z_{2,1}, z_{1,1}$ and $w_1$, and thus, it follows that there is some point $z_{2,3} \in \gamma_{x y_{1}}[z_{1,1}, y_{1}]$ such that
$
k_{\Omega}(w_1, z_{2,3}) \leq C.
$

To finish the proof of the assertion, let $z\in\gamma_{xy_1}[z_{2,3},y_1]$. Then

\beq\label{25-12-15-2}
k_{\Omega}(z_{1}, z) &\geq & k_{\Omega}(z_{1,1}, z)- k_{\Omega}(z_{1}, z_{1,1})\geq   k_{\Omega}(z_{1,1}, z_{2,3})-k_{\Omega}(z_{1}, z_{1,1})\\ \nonumber & {\geq} & k_{\Omega}(z_{2,1},w_1)-k_{\Omega}(z_{2,3},w_1)-k_{\Omega}(z_{1,1},z_{2,1})-k_{\Omega}(z_{1}, z_{1,1})
\eeq

By substituting \eqref{H25-26-12} and \eqref{H25-26-10}$-$\eqref{25H-0724-10} to \eqref{25-12-15-2}, we see from the arbitrariness of $z$ in $\gamma_{xy_1}[z_{2,3},y_1]$ that \eqref{H25-0505-1} is true.
	
We end the preparation with the choice of a suitable point from $\gamma_{x y_{3}}$. 		
Note that \eqref{25H-0724-7} and \eqref{H25-26-12} ensure that
$k_{\Omega}(z_{2,3}, \gamma_{y_1,y_3}) \geq k_{\Omega}(w_1, \gamma_{y_1,y_3})-k_{\Omega}(w_1, z_{2,3})>C$. Then the assumption of $(\Omega,k)$ being $C$-Gromov hyperbolic implies that there exists some point $z_{4,3} \in \gamma_{x y_{3}}$ such that
\be\label{eq-23-5}
k_{\Omega}(z_{2,3}, z_{4,3}) \leq C.
\ee
	
Now, we are ready to prove \eqref{25-12-15-1}. Let us divide the discussions into two cases.
	\smallskip
	
	\textbf{Case I}: $z_{5,3} \in \gamma_{x y_{3}}[x, z_{4,3}]$.
	\smallskip

	In this case, since $(\Omega,k)$ is $C$-Gromov hyperbolic, we may find some point  $z_{6,3} \in \gamma_{xy_1}[x, z_{2,3}] \cup \gamma_{z_{2,3} z_{4,3}}$ such that
	\beq\label{202511-17-1}
	k_{\Omega}(z_{5,3}, z_{6,3}) \leq C.
	\eeq
	Then \eqref{H25-0505-2} gives
	\be\label{H25-0505-3}k_{\Omega}(z_1,z_{6,3})\leq k_{\Omega}(z_1,z_{5,3})+k_{\Omega}(z_{5,3}, z_{6,3})\leq 4C.\ee
	This implies that $z_{6.3} \in \gamma_{x y_{1}}[x, z_{2,3}]$. Indeed, if not, then $k_{\Omega}(z_{2,3},z_{6,3})\leq k_{\Omega}(z_{2,3},z_{4,3})$,
	and thus, \eqref{H25-0505-1} and \eqref{eq-23-5} lead to
	$$k_{\Omega}(z_1,z_{6,3})\geq k_{\Omega}(z_1,z_{2,3})-k_{\Omega}(z_{2,3},z_{6,3})\geq
	5C,$$
	which contradicts with \eqref{H25-0505-3}.
	This fact implies that
	$
	k_{\Omega}(z_{1}, z_{5,3}) \geq k_{\Omega}(z_{1}, z_{6,3})-k_{\Omega}(z_{5,3}, z_{6,3}).
	$
	Then by the assumption $\gamma_{x z_{1,1}y_{1}} \in P_{\alpha}^{\gamma_{xz_{1}y}}(3C)$ and Lemma \ref{Lemma4-1.1}, 
	 we know from \eqref{202511-17-1} that
	\beqq\label{H25-27-6}
	k_{\Omega}(z_{1}, z_{5,3}) \geq   k_{\Omega}(z_{1},z_{1,1}) -k_{\Omega}(z_{5,3}, z_{6,3})\geq 2C.
	\eeqq
	
	
	\textbf{Case II}: $z_{5,3} \in \gamma_{x y_{3}}[z_{4,3}, y_{3}]$.
	\smallskip
	
	In this case, since $(\Omega,k)$ is $C$-Gromov hyperbolic, there exists some $z_{7,3}\in\gamma_{xy_1}\cup \gamma_{y_1y_3}$
	such that \beqq\label{H25-27-7}k_{\Omega}(z_{5,3},z_{7,3})\leq C.\eeqq
Then it follows from  (\ref{H25-0505-2}) that
	$$k_{\Omega}(z_1,z_{7,3})\leq k_{\Omega}(z_1,z_{5,3})+k_{\Omega}(z_{5,3},z_{7,3})\leq 4C,$$
	and so, (\ref{eq-23-1-1}) yields  $z_{7,3}\notin\gamma_{y_1y_3}$. Moreover, \eqref{H25-0505-1} yields  $z_{7,3}\notin\gamma_{xy_1}[z_{2,3},y_1]$.
	Thus $z_{7,3}\in\gamma_{xy_1}[x,z_{2,3}]$. Based on this fact, a similar discussion as in \eqref{H25-27-6}
	shows that $$k_{\Omega}(z_1,z_{5,3})\geq 2C.$$

	The proof of the lemma is thus complete.
	\epf
	
	The following are two more technical results.
	
	\blem\label{lem-22-3.1} Fix  $\gamma_{xx_1y}\in \Lambda_{xy}(\Omega)$, $\alpha\in \Gamma_{xy}(\Omega)$ and $\gamma_{xz_1z}\in Q_{\alpha}^{\gamma_{xx_1y}}$. Suppose that there are  $y_1\in\gamma_{xy}[x_1,y]$ and $w\in \alpha[z, y]$ such that $\gamma_{xw_1w}\in Q_{\alpha}^{\gamma_{xy_1y}}$.
	If $k_{\Omega}(x_1,y_1)\geq 30C$, then the following assertions hold:
	\ben
	\item\label{Lemma4-1}
	For any $\gamma_{zw}\in\Lambda_{zw}(\Omega)$, there exists some point $w_0\in\gamma_{zw}$ such that $\gamma_{zw_0w}\in O_{\alpha[z,w]}^{\gamma_{xx_1y}}(2C)$.
	\item\label{Lemma4-2}  For each  $u\in \gamma_{xw_1w}[w_1,w]$ and $v\in \gamma_{xz_1z}[z_1,z]$, $k_{\Omega}(u,v) \geq3C$.
	\een
	\elem
	\bpf First of all, fix $\gamma_{zw}\in \Lambda_{zw}(\Omega)$. Since $\gamma_{xw_1w}\in Q_{\alpha}^{\gamma_{xy_1y}}$, $k_{\Omega}(w_1,y_1)\leq 7C$, and thus,
	$$k_{\Omega}(x_1,y_1)\geq 30C\geq 23C+k_{\Omega}(w_1,y_1),$$
	which, together with Lemma \ref{lem-2-3.0} (with $x=x$, $w=x_1$, $y=y_1$ and $z=w_1$), yields that
	there exists some point	$y_2\in\gamma_{xw}[x,w_1]$ such that
	\be\label{H25-03-2}k_{\Omega}(x_1, y_2)\leq C.\ee
	
	Since the assumption $\gamma_{xw_1w}\in Q_{\alpha}^{\gamma_{xy_1y}}$ in the lemma ensures that $k_{\Omega}(y_1, \gamma_{xw_1w})\geq 2C$. Then the fact
	of $(\Omega,k)$ being $C$-Gromov hyperbolic implies that there exists some point $w_2\in \gamma_{wy}$ $(\in \Lambda_{wy}(\Omega))$  such that
	\be\label{25H-0724-1} k_{\Omega}(y_1,w_2)\leq C.\ee
	Hence \be\label{25H-0724-2} k_{\Omega}(w_1,w_2)\leq k_{\Omega}(w_1,y_1)+k_{\Omega}(y_1,w_2)\leq 8 C\ee
	and \be\label{25H-0724-0} k_{\Omega}(w_1,y_2)\geq k_{\Omega}(x_1,y_1)-k_{\Omega}(x_1,y_2)-k_{\Omega}(y_1,w_1)\geq 21C.\ee

	Now, we are ready to prove Lemma \ref{lem-22-3.1}(\ref{Lemma4-1}). Since $\gamma_{xz_1z}\in Q_{\alpha}^{\gamma_{xx_1y}}$, $k_{\Omega}(x_1,z_1)\geq 2C$, and thus, there exists some point $x_2\in\gamma_{zy}$ $(\in \Lambda_{zy}(\Omega))$
	such that \be\label{25H-0724-3} k_{\Omega}(x_1,x_2)\leq C.\ee
	
	As $(\Omega,k)$ is $C$-Gromov hyperbolic, there exists some point $w_0\in \gamma_{zw}\cup \gamma_{wy}$ such that
	\be\label{H25-0506-1}k_{\Omega}(x_2,w_0)\leq C,\ee
	and so,
	\be\label{25H-0724-4}k_{\Omega}(x_1,w_0)\leq k_{\Omega}(x_1, x_2)+k_{\Omega}(x_2,w_0)\stackrel{ \eqref{25H-0724-3}\; + \;\eqref{H25-0506-1}}{\leq} 2C.\ee
	
	We claim that $w_0\in \gamma_{zw}$. Otherwise,
	\be\label{25H-0724-5}k_{\Omega}(w_0,w_2)\geq k_{\Omega}(x_1,y_1)-k_{\Omega}(x_1,w_0)-k_{\Omega}(y_1,w_2)
	\stackrel{ \eqref{25H-0724-1}\; + \;\eqref{25H-0724-4}}{\geq} 17C.\ee
	
	If $w_0\in \gamma_{wy}[w,w_2]$, then (\ref{25H-0724-2}) and (\ref{25H-0724-5}) imply
	that $$k_{\Omega}(w_0,w_2)\geq 9C+k_{\Omega}(w_1,w_2).$$
	Applying Lemma \ref{lem-2-3.0} (with $x=w,y=w_2$, $z=w_1$ and $w=w_0$), we obtain that there exists some $w_3\in \gamma_{xw}[w,w_1]$ such that
	$$k_{\Omega}(w_0,w_3)\leq C,$$
	which, together with (\ref{H25-03-2}), (\ref{25H-0724-3})  and (\ref{H25-0506-1}), shows
	$$k_{\Omega}(y_2,w_3)\leq k_{\Omega}(y_2,x_1)+k_{\Omega}(x_1,x_2)+k_{\Omega}(x_2,w_0)+k_{\Omega}(w_0,w_3)\leq 4C,$$
	which contradicts with
	$$k_{\Omega}(y_2,w_3)\geq k_{\Omega}(w_1,y_2) \stackrel{ \eqref{25H-0724-0}}{\geq} 21C.$$
	
	If $w_0\in \gamma_{wy}[w_2,y]$, then a similar discussion as above will lead to a contradiction. Hence $w_0\in \gamma_{zw}$, and thus, \eqref{25H-0724-4} implies that the lemma is true.
	
	Next, we prove Lemma \ref{lem-22-3.1}(\ref{Lemma4-2}).
	Since $k_{\Omega}(x_1, y_1) \geq 30C$, we get from  \eqref{H25-03-2} and the assumption $\gamma_{xw_1w}\in Q_{\alpha}^{\gamma_{xy_1y}}$ that for each $u\in \gamma_{xw_1w}[w_1,w]$, it holds
	\be\label{H25-0506-2}
	k_{\Omega}(u,w_2) \geq k_{\Omega}(w_2, w_1) \geq k_{\Omega}(x_1, y_1)-k_{\Omega}(x_1, w_2)-k_{\Omega}(y_1, w_1) \geq 22C.
	\ee
	Then by \eqref{H25-0506-1}, $$k_{\Omega}(u,w_2) \geq 22C\geq 21C+k_{\Omega}(w_2,w_0).$$
	Applying Lemma \ref{lem-2-3.0} (with $x=w,y=w_2$ and $z=w_0$), we obtain that there exists some point  $u_1 \in \gamma_{xw}[w_0,w]$ such that $$k_{\Omega}(u,u_1) \leq C.$$ Then by (\ref{H25-0506-1}) and (\ref{H25-0506-2}),
	\be\label{H25-0506-3}k_{\Omega}(u_1,w_0)\geq k_{\Omega}(u,w_2) -k_{\Omega}(u,u_1)-k_{\Omega}(w_2,w_0)\geq 20C.\ee
	Since $\gamma_{xz_1z}\in Q_{\alpha}^{\gamma_{xx_1y}}$, it follows from Lemma \ref{lem-22-3.1}(\ref{Lemma4-1}) that
	$$k_{\Omega}(z_1,w_0)\leq k_{\Omega}(x_1,z_1)+k_{\Omega}(x_1,w_0)\leq 9C.$$
	
	If $k_{\Omega}(z_1,v)\geq 11C$, then $$k_{\Omega}(z_1,v)\geq 2C+k_{\Omega}(z_1,w_0).$$
	Applying Lemma \ref{lem-2-3.0} (with $x=z$, $y=z_1$ and $z=w_0$), we obtain that there exists some point  $v_1 \in \gamma_{zw}[z,w_0]$ such that $$k_{\Omega}(v,v_1) \leq C,$$
	which, together with (\ref{H25-0506-3}), shows that
	$$k_{\Omega}(u,v)\geq k_{\Omega}(u_1,v_1)-k_{\Omega}(v,v_1)-k_{\Omega}(u,u_1)\geq k_{\Omega}(u_1,w_0)-2C\geq 18C.$$
	
	If $k_{\Omega}(z_1,v)< 11C$, then we obtain from \eqref{H25-03-2} and \eqref{H25-0506-2} that
	$$
	\begin{aligned}
		k_{\Omega}(u,v)&\geq k_{\Omega}(u,w_2)-k_{\Omega}(w_2,z_1)-k_{\Omega}(z_1,v)\geq 11C-k_{\Omega}(w_2,z_1)\\
		&\geq 11C-k_{\Omega}(w_2,x_1)-k_{\Omega}(x_1,z_1)\geq 3C.
	\end{aligned}
	$$
	This completes the proof of Lemma \ref{lem-22-3.1}(\ref{Lemma4-2}).
	\epf
	\medskip

	\blem\label{lem-24-5.3}
	Fix $\gamma_{xz_1y}\in \Lambda_{xy}(\Omega)$ and $\alpha\in\Gamma_{xy}(\Omega)$. Suppose that $y_{1} \in \alpha$, $y_{2} \in \alpha[y_{1}, y]$ and $y_{3} \in \alpha[y_{1}, y_{2}]$. Then the following two assertions hold:
	
	\begin{enumerate}
		\item \label{lem-24-5.3-1} If $\gamma_{xz_{1,1}y_{1}}\in P_{\alpha}^{\gamma_{xz_{1}y}}(3C)$ and $\gamma_{y_{1}z_{1,2}y_{2}}\in O_{\alpha[y_{1},y_{2}]}^{\gamma_{xz_1y}}(2C)$, then for each $z \in \gamma_{xy_{1}}[y_{1}, z_{1,1}]$ with $k_{\Omega}(z_{1,1}, z) \geq 10C$,
		$$k_{\Omega}(z,\gamma_{y_{1}y_{2}}[z_{1,2},y_{1}])\leq C.$$
		
		\item\label{lem-24-5.3-2} If $\gamma_{xz_{1,1}y_{1}} \in Q_{\alpha}^{\gamma_{xz_{1}y}}$, $\gamma_{x z_{2,1} y_{3}} \in Q_{\alpha}^{\gamma_{x z_{1}y}}$ and $k_{\Omega}(z_{1}, \gamma_{y_{1}y_{3}}) \geq 11 C$, then for each  $z\in \gamma_{y_{1} y_{3}}$,
		$$k_{\Omega}(z, \gamma_{x y_{1}}[y_{1}, z_{1,1}] \cup \gamma_{x y_{3}}[z_{2,1}, y_{3}]) \leq C.$$
		
	\end{enumerate}
	\elem
	
	\bpf
	(1). Since $\gamma_{xz_{1,1}y_{1}}\in P_{\alpha}^{\gamma_{xz_{1}y}}(3C)$, Lemma \ref{Lemma4-1.1} gives $k_{\Omega}(z_{1}, z_{1,1})<5C$. Since $\gamma_{y_{1}z_{1,2}y_{2}}\in O_{\alpha[y_{1}y_{2}]}^{\gamma_{xz_1y}}(2C)$, we obtain
	$$
	k_{\Omega}(z_{1,1}, z_{1,2}) \leq k_{\Omega}(z_{1}, z_{1,1})+k_{\Omega}(z_{1}, z_{1,2}) <7C,
	$$
	and thus, $$k_{\Omega}(z_{1,1}, z) \geq 10C \geq 3C+k_{\Omega}(z_{1,1}, z_{1,2}).$$
	Applying Lemma \ref{lem-2-3.0} (with $x=y_1$, $y=z_{1,1}$,  $z=z_{1,2}$ and $w=z$), we infer that there exists some point $w \in \gamma_{y_{1}y_2}[z_{1,2}, y_{1}]$ such that
	$$
	k_{\Omega}(z, w) \leq C,
	$$
	from which Lemma \ref{lem-24-5.3}\eqref{lem-24-5.3-1} follows.
	
	(2). Since $(\Omega,k)$ is $C$-Gromov hyperbolic, for each $z\in \gamma_{y_{1}y_{3}}$, there exists some $w\in \gamma_{x y_{1}}\cup \gamma_{x y_{3}}$ such that
	\be\label{H25-31-1}k_{\Omega}(z, w)\leq C.\ee
	We shall prove
	\be\label{eq-24-1}
	w\in \gamma_{x y_{3}}[z_{2,1}, y_{3}]\cup \gamma_{x y_{1}}[z_{1,1}, y_{1}]
	\ee
	via a contradiction argument.
	
	Suppose that \eqref{eq-24-1} fails. Then $w \in \gamma_{xy_{3}}[x, z_{2,1}]\cup \gamma_{x y_{1}}[x,z_{1,1}]$. As the discussion is similar for $w \in \gamma_{x y_{1}}[x,z_{1,1}]$, without loss of generality, we may assume $w \in \gamma_{xy_{3}}[x, z_{2,1}]$. Since $\gamma_{x z_{2,1} y_{3}} \in Q_{\alpha}^{\gamma_{x z_{1}y}}$, $k_{\Omega}(z_{1}, z_{2,1}) \leq 7C$, and thus, by the triangle inequality and \eqref{H25-31-1}, we have 
	$$
	\begin{aligned}
		k_{\Omega}(z_{2,1},w)&\geq k_{\Omega}(z_{2,1},z)-k_{\Omega}(z,w)\geq k_{\Omega}(z_1,z)-k_{\Omega}(z_{2,1},z_1)-C\\
		&\geq k_{\Omega}(z_1,\gamma_{y_1y_3})-8C\geq 3C\geq 2C+k_{\Omega}(z,w).
	\end{aligned}
	$$
	Then we may apply Lemma \ref{lem-2-3.0} (with $x=y_3$, $w=z_{2,1}$, $y=w$ and $z=z$) to obtain a point $u\in \gamma_{y_{1}y_{3}}[z, y_{3}]$ with
	$$
	k_{\Omega}(z_{2,1}, u) \leq C.
	$$
	Consequently, we get
	$$k_{\Omega}(z_{1}, \gamma_{y_{1}y_{3}})\leq k_{\Omega}(z_{1}, z_{2,1})+k_{\Omega}(z_{2,1}, u) \leq 8C,$$
	which clearly contradicts with the assumption $k_{\Omega}(z_{1}, \gamma_{y_{1}y_{3}}) \geq 11 C$.
	Thus \eqref{eq-24-1} holds, and hence, the proof of lemma is complete.
	\epf
	
	
	%
	%
	
	\medskip
	\textbf{Acknowledgment.} C.-Y. Guo is supported by the Young Scientist Program of the Ministry of Science and Technology of China (No.~2021YFA1002200), the NSF of China (No.~12311530037), the Taishan Scholar Project and the Jiangsu Provincial Scientific Research Center of Applied Mathematics (Grant No.~BK20233002). M. Huang and X. Wang are partly supported by NSF of China (No.~12371071 and No.~12571081).
	
	We would like to thank Prof.~Pekka Koskela for numerous helpful discussions and comments, in particular, for bringing Question \ref{Ques:metric version} to our attention. We also thank Prof. Xiangdong Xie and Dr. Abhishek Pandey for their valuable comments and suggestions.

\end{document}